\documentclass[a4paper,10pt]{article}

\usepackage[T1]{fontenc}
\usepackage[all]{xy}
\usepackage{amsmath,amssymb,amsthm,geometry}
\title{Coverings of Laura Algebras: the Standard Case}

\author{Ibrahim Assem
\\
\and
Juan Carlos Bustamante
\\
\and
Patrick Le Meur
}

\date{\today}

\theoremstyle{plain}
\newtheorem{prop}{Proposition}[section]
\newtheorem{lem}[prop]{Lemma}
\newtheorem{cor}[prop]{Corollary}
\newtheorem{Thm}{Theorem}[]
\newtheorem{thm}[prop]{Theorem}

\theoremstyle{definition}
\newtheorem{df}[prop]{Definition}
\newtheorem{pb}{Problem}[]
\theoremstyle{remark}
\newtheorem{ex}[prop]{Example}
\newtheorem{exs}[prop]{Examples}
\newtheorem{rem}[prop]{Remark}

\makeatletter
\renewcommand\section{\@startsection {section}{1}{0mm}%
                                   {-3.5ex \@plus -1ex \@minus -.2ex}%
                                   {2.3ex \@plus.2ex}%
                                   {\normalfont\large\bfseries}}

\renewcommand\subsection{\@startsection {subsection}{2}{0mm}%
  {-3.5ex \@plus -1ex \@minus -.2ex}%
  {0.5ex \@plus.2ex}%
  {\normalfont\normalsize\bfseries}
}

\makeatother

\def\sq{\null\hfill$\blacksquare$\\}
\def\ts#1{\text{\normalfont\textsf{#1}}}

\begin{document}

\maketitle

\abstract{
In this paper, we study the covering theory of laura algebras. We
prove that if a connected laura algebra is standard (that is, has
a standard connecting component), then it has 
Galois coverings associated to
the coverings of the connecting component. As a consequence,  the
first Hochschild cohomology group of a standard laura algebra 
vanishes if and only if it has no proper Galois coverings.
}
\section*{Introduction}
$\ $

Introduced in $1945$, the Hochschild cohomology groups are subtle and
interesting invariants of associative algebras. The lower dimensional
groups have simple interpretations: for instance, the $0^{th}$ group
is the centre of the algebra, the $1$st group can be thought of as the
group of outer derivations of the algebra, while the $2$nd  and $3$rd
groups are related to the rigidity properties of the algebra. In
\cite[\S 3, Pb. 1]{S93a}, Skowro\'nski has related the vanishing
of the first Hochschild cohomology group $\ts{HH}^1(A)$ of an algebra
$A$ (with coefficients in the bimodule $_AA_A$) to the simple
connectedness of $A$. Recall that a basic and connected finite
dimensional algebra over an algebraically closed field $k$ is
\emph{simply connected} if it has no proper Galois covering or,
equivalently, if the fundamental group (in the sense of
\cite{martinezvilla_delapena}) of any presentation is trivial. In
particular, Skowro\'nski posed the following problem: for which algebras $A$ do
we have $\ts{HH}^1(A)=0$ if and only if $A$ is simply connected? This
problem has been the subject of several investigations: notably this
equivalence holds true for algebras derived equivalent to hereditary algebras
\cite{lemeur2}, weakly shod algebras \cite{ws} (see also
\cite{al}), large classes of selfinjective algebras \cite{R02} and
 schurian cluster-tilted algebras \cite{AR09}. It was proved in
 \cite{bl} that, for
a representation-finite algebra, the first Hochschild cohomology group
vanishes if and only if its Auslander-Reiten quiver is simply
connected.
Note that if $A$ is a representation-finite triangular algebra, then
its Auslander-Reiten quiver is simply connected if and only if $A$
has no proper Galois covering, that is, $A$ is simply connected.

Here, we study this conjecture for laura algebras. These are defined
as follows. Let $\ts{mod}\, A$ be the category of finitely generated
right $A$-modules, and $\ts{ind}\,A$ be a full subcategory consisting
of exactly one representative from each isomorphism class of
indecomposable $A$-modules. The \emph{left part} $\mathcal L_A$ of
$\ts{mod}\,A$ is the full subcategory of $\ts{ind}\,A$ consisting of
those modules whose predecessors have projective dimension at most
one, and the \emph{right part} $\mathcal R_A$ is defined dually. These
classes were introduced in \cite{hrs} in order to study the module
categories of quasi-tilted algebras. Following \cite{ac,sk2}, we say
that $A$ is \emph{laura} provided $\ts{ind}\,A\backslash\left(\mathcal
  L_A\cup \mathcal R_A\right)$ has only finitely many objects. Part of
the importance of laura algebras comes from the fact that this class
contains (and generalises) the classes of  representation-finite
algebras,  tilted, quasi-tilted and weakly shod algebras. Laura
algebras have appeared naturally in the study of Auslander-Reiten
components:  an Auslander-Reiten component is called
\emph{quasi-directed} if it is generalised standard and almost all its
modules are directed. It was shown in \cite{ac} that a laura algebra
which is  not quasi-tilted has a unique faithful
convex quasi-directed Auslander-Reiten component (which is also the unique
non-semiregular component). Conversely,  any convex quasi-directed component occurs in this
way \cite{S04a}. The techniques used for the study of laura algebras
were applied in \cite{LS06} to obtain useful results on the infinite radical of the
module category. Their representation dimension is at most three and
this is a class of algebras with possibly infinite global dimension
which satisfies the finitistic dimension conjecture \cite{apt}. Also,
laura algebras have been characterised in terms
of the Gabriel-Rojter measure as announced by Lanzilotta in the
ICRA~XI in Mexico, $2004$ (see also \cite{ACLST05}). For further
properties of laura algebras we refer the reader to
\cite{ac,ac04,act,alr,sk2}. Here we concentrate on the conjecture
that a laura algebra $A$ is simply connected if and only if
$\ts{HH}^1(A)=0$.

Our approach, already used in \cite{ws,lemeur2}, uses
coverings. Covering theory was introduced by Gabriel and his school
(see, for instance, \cite{bog,gabriel,r}) and consists in replacing an
algebra by a locally bounded category, called its covering, which is
sometimes easier to study. We recall that a tilted algebra is
characterised by the existence of at least one, and at most two,
\emph{connecting components} (it has two if and only if it is
concealed, in which case the connecting components are
postprojective and preinjective) see \cite{ass}. If $A$ is a laura
not quasi-tilted algebra, then its unique faithful quasi-directed
component is also called a connecting component (see
\cite{ac}). Hence, by laura algebra \emph{with connecting component}, we mean
a connected laura algebra which is either tilted or not
quasi-tilted. We call a laura algebra with connecting component
\emph{standard} provided its connecting components are all standard (it is
known from \cite{R88} that the connecting components of concealed
algebras are standard).
 This generalises the
 notion of standard representation-finite algebra (see
 \cite{bog}). Several classes of laura algebras are standard, notably
  tilted algebras or  weakly shod algebras. Our first main
 theorem says that if $\widetilde{\Gamma}\to \Gamma$ is a Galois
 covering of the connecting component such that there exists a
 well-behaved covering functor $k(\widetilde{\Gamma})\to \ts{ind}\,\Gamma$
 then it induces a  covering of the algebra.
\begin{Thm}
  \label{thm2}
Let $A$ be a laura algebra with connecting component $\Gamma$ and
$\pi\colon\widetilde{\Gamma}\to\Gamma$ be a Galois covering with group
$G$ with respect to which there exists a well-behaved
covering functor $p\colon k(\widetilde{\Gamma})\to \ts{ind}\,\Gamma$. Then
there exists a covering functor  $F\colon \widetilde A\to A$ whose fibres are in bijection
with  $G$. If moreover $A$ is standard, then $F$ is a Galois covering with
group $G$.
\end{Thm}
Note that if $\Gamma$ is standard then there always exists a
well-behaved covering functor $p$.

In order to prove Theorem~\ref{thm2}, we consider a more
general situation. We first consider an Auslander-Reiten
component, which contains a left section (in the sense
of \cite{assem}) and  show that to a Galois covering of this
component such that there is a corresponding well-behaved functor
corresponds a covering of its support algebra with
nice properties, see Theorem~\ref{prop2} below. Applying this result to the connecting component of a
laura algebra yields the required covering.

Because of the theorem,  if $A$ is standard, then we are able
to work with Galois
coverings which are notably easier to handle than covering
functors. We prove that if $A$ is  standard laura, then any
Galois covering of the connecting component induces a Galois
covering of $A$, with the same group. This allows us to prove
our second main theorem, which settles the conjecture for
standard laura algebras.
\begin{Thm}
  \label{thm1}
Let $A$ be a standard laura algebra, and $\Gamma$ its connecting component(s). The following are
equivalent:
\begin{enumerate}
\item[(a)] $A$ has no proper Galois covering, that is, $A$ is simply connected.
\item[(b)] $\ts{HH}^1(A)=0$.
\item[(c)] $\Gamma$ is simply connected.
\item[(d)] The orbit graph $\mathcal O(\Gamma)$  is a tree.
\end{enumerate}
Moreover, if these conditions are verified, then $A$ is weakly shod.
\end{Thm}

If one drops the standard condition, then the above theorem may
fail. Indeed, there are examples of non-standard representation-finite
algebras which have no proper Galois covering and with non-zero first Hochschild
cohomology group (see \cite{bog,bl}, or below). However, some
implications are still true in Theorem~\ref{thm1} without
assuming standardness. Indeed, we always have: (c) and (d) are
equivalent and (c) implies (a) and (b).

Our paper is organised as follows. After a short preliminary section,
we prove a few preparatory lemmata on covering
functors in Section~\ref{sec:s2}. In Section~\ref{sec:s3}, we give examples of
standard laura algebras. Section~\ref{sec:s4} is devoted to properties
of tilting modules which are in the image of the push-down functor
associated to a covering functor. In Section~\ref{sec:s5bis} we study
the coverings of Auslander-Reiten components having left sections. The proof of Theorem~\ref{thm2}
occupies Section~\ref{sec:s5}. We concentrate on Galois
coverings in Section~\ref{sec:s6}, and prove
Theorem~\ref{thm1} in Section~\ref{sec:s7}.

\section{Preliminaries}
\label{sec:s1}
\paragraph{Categories and modules}$\ $

Throughout this paper, $k$ denotes a fixed algebraically closed
field. All our categories are locally bounded $k$-categories, in the
sense of \cite[2.1]{bog}. We assume that all locally bounded $k$-categories
are small and
all functors are $k$-linear (the categories of finite dimensional
modules and their bounded derived categories are skeletally small).

Let $F\colon\mathcal E\to \mathcal B$ be a $k$-linear functor and $G$ be a
group acting on $\mathcal E$ and $\mathcal B$ by automorphisms. Then $F$ is called
\emph{$G$-equivariant} if $F\circ 
g=g\circ F$ for every $g\in G$.

A basic finite dimensional algebra $A$ can be
considered equivalently as a locally bounded $k$-category as follows:
Fix a complete set $\{e_1,\ldots,e_n\}$ of primitive orthogonal
idempotents, then the object set of $A$ is the set
$\{e_1,\ldots,e_n\}$ and the morphisms space from $e_i$ to $e_j$ is
$e_jAe_i$. The composition of morphisms is induced by
the multiplication in $A$.

Let $\mathcal C$ be a locally bounded $k$-category. We denote by $\mathcal C_o$
its object class. A \emph{right $\mathcal
  C$-module} $M$ is a $k$-linear functor $M\colon \mathcal C^{op}\to \ts{MOD}\, 
k$, where $\ts{MOD}\, k$ is the category of $k$-vector spaces. 
We write
$\ts{MOD}\, \mathcal C$ for the category of $\mathcal C$-modules and $\ts{mod}\, \mathcal C$ for the
full subcategory of the \emph{finite dimensional} $\mathcal C$-modules,
that is, those modules $M$ such that $\sum\limits_{x\in \mathcal C_o}\dim\
M(x)<\infty$. If $\mathcal A$ is a subcategory of $\ts{MOD}\,\mathcal C$, we use the
notation $X\in\mathcal A$ to express that $X$ is an object in $\mathcal A$.
For every $x\in \mathcal C_o$,
the indecomposable projective $\mathcal C$-module associated to $x$ is $\mathcal
C(-,x)$. The standard duality $\ts{Hom}_k(-,k)$ is denoted by $D$. 
Let $M$ be a $\mathcal C$-module. If $\mathcal B$ is a full
subcategory of $\mathcal C$, then $M_{|\mathcal B}$ is the induced $\mathcal
B$-module. If $\mathcal X$ is a subcategory of $\ts{mod}\, \mathcal C$, then the $\mathcal
X$-module $\ts{Hom}_{\mathcal C}(-,M)_{|\mathcal
  X}$ is denoted by $\ts{Hom}_{\mathcal C}(\mathcal X,M)$.  Also, $\ts{Hom}_{\mathcal
  C}(M,\mathcal C)$ denotes the $\mathcal
C^{op}$-module $\ts{Hom}_{\mathcal C}(M,\bigoplus\limits_{x\in \mathcal C_o}\mathcal
C(-,x))$ (if $A=\mathcal C$ is a finite dimensional algebra, this is just
the left $A$-module $\ts{Hom}_A(M,A)$).

We let $\ts{ind}\, \mathcal C$ be a full subcategory of $\ts{mod}\, \mathcal
C$ consisting of a complete set of representatives of the isomorphism
classes of  indecomposable $\mathcal C$-modules. We write $\ts{proj}\,\mathcal C$ and
$\ts{inj}\,\mathcal C$ for the full subcategories of $\ts{ind}\,\mathcal C$ of projective and
injective modules, respectively. Whenever we speak about
an indecomposable $\mathcal C$-module, we always mean that it belongs to
$\ts{ind}\, \mathcal C$. 

For a full subcategory $\mathcal A$ of $\ts{mod}\, \mathcal C$, we denote by $\ts{add}\, \mathcal A$
 the full subcategory of $\ts{mod}\, \mathcal C$ with objects the direct sums
of summands of modules in $\mathcal A$. If $M$ is a module, then $\ts{add}\, M$
denotes $\ts{add}\,\{M\}$.

The Auslander-Reiten translations in $\ts{mod}\, \mathcal C$ are denoted by
$\tau_{\mathcal C}=D\,\ts{Tr}$ and $\tau_{\mathcal C}^{-1}=\ts{Tr}\,D$. The Auslander-Reiten
quiver of $\mathcal C$ is denoted by $\Gamma(\ts{mod}\, \mathcal C)$. For a component $\Gamma$
of $\Gamma(\ts{mod}\, \mathcal C)$, we denote by $\mathcal O(\Gamma)$ its orbit graph (see
\cite[4.2]{bog}, or Section~\ref{sec:s7}
below). The component $\Gamma$  is \emph{non-semiregular} if it contains both an
injective and a projective module. It is \emph{faithful} if its
\emph{annihilator} $\ts{Ann}\,\Gamma=\bigcap\limits_{X\in\ \Gamma}\ts{Ann}\, X$ is
zero. Following \cite{sk}, a component  $\Gamma$ is
\emph{generalised standard} if $\ts{rad}^{\infty}(X,Y)=0$ for every
$X,Y\in \Gamma$. Denoting by $k(\Gamma)$ the mesh category of $\Gamma$
(see \cite[2.5]{bog}), $\Gamma$ is \emph{standard} if there exists an
isomorphism of
$k$-categories $k(\Gamma)\xrightarrow{\sim}\ts{ind}\,\Gamma$
which extends the
identity  on vertices, and which maps meshes to almost split
sequences. 
Let $\pi\colon\widetilde{\Gamma}\to \Gamma$ be a morphism of translation quivers. Let $\mathcal
X$ be a full convex subquiver of $\widetilde{\Gamma}$. We let $k(\mathcal X)$ be the full
subcategory of $k(\widetilde{\Gamma})$ with objects the vertices in $\mathcal X$. Following
\cite[3.1]{bog}, a functor
$p\colon k(\mathcal X)\to \ts{ind}\, \Gamma$ is \emph{well-behaved} (with
respect to $\pi$) if it satisfies:
\begin{enumerate}
\item $p(X)=\pi(X)$ for every $X\in\mathcal X$.
\item Let $X\in\mathcal X$. Let $(u_i\colon Z_i\to X)_{i=1,\ldots,t}$ be all
  the arrows in $\mathcal X$ ending at $X$ (or $(v_j\colon X\to Y_j)_{j=1,\ldots,s}$ be all
  the arrows in $\mathcal X$  starting from $X$), then the morphism $\begin{bmatrix}
      p(u_1)&\ldots&p(u_t)
    \end{bmatrix}\colon\bigoplus\limits_{i=1}^tp(Z_i)\to p(X)$ (or $\begin{bmatrix}
    p(v_1)&\ldots &p(v_s)
  \end{bmatrix}^t\colon p(X)\to\bigoplus\limits_{j=1}^s Y_j$,
  respectively) is irreducible.
% \item Let $X\in\widetilde{\Gamma}$ be non projective. Let $(u_i\colon \tau X\to
%   E_i)_{i=1,\ldots,s}$ (or $(v_i\colon E_i\to X)_{i=1,\ldots,s}$) be all
%   the arrows in $\widetilde{\Gamma}$ starting from $\tau X$ (or ending at $X$,
%   respectively). If the mesh in $\widetilde{\Gamma}$ ending at $X$ is contained in $\mathcal
%   X$, then  $\sum\limits_{i=1}^sp(v_i)p(u_i)=0$.
\end{enumerate}
Condition $2$  above imply that if a
mesh in $\widetilde{\Gamma}$ is contained in $\mathcal X$, then 
$p$ maps this mesh to an almost split sequence.

For notions and results on modules, we refer the reader to
\cite{ass}. For coverings and fundamental groups of
translation quivers, we refer the reader to \cite[\S 1]{bog}. Note
that the translation quivers we consider are not valued translation
quivers and may have multiple arrows

\paragraph{Paths}$\ $

Let $\mathcal C$ be a locally bounded $k$-category. Let $X,Y$ be in $\ts{ind}\, \mathcal
C$. Following the convention used in \cite{hrs}, a \emph{path}
$X\rightsquigarrow Y$ from $X$ to $Y$ in $\ts{ind}\, \mathcal 
C$ is a sequence of non-zero morphisms:
\begin{equation}
(\star)\ \ X=X_0\xrightarrow{f_1}X_1\to\cdots\to
X_{t-1}\xrightarrow{f_t}X_t=Y\ (t\geqslant 0)\notag
\end{equation}
where $X_i\in \ts{ind}\, \mathcal C$ for all $i$. We then say that $X$ is a
\emph{predecessor} of $Y$ and that $Y$ is a \emph{successor} of
$X$. A path from $X$ to $X$ involving at least one non-isomorphism is
a \emph{cycle}. A module $X\in \ts{ind}\, \mathcal C$ which lies on no cycle is
\emph{directed}. If each $f_i$ in $(\star)$ is irreducible, we say
that $(\star)$ is a \emph{path of irreducible morphisms} or a
\emph{path in}
 $\Gamma(\ts{mod}\, \mathcal C)$. A path $(\star)$ of irreducible morphisms is
\emph{sectional} if $\tau_{\mathcal C}X_{i+1}\neq X_{i-1}$ for all
$i$ with $0<i<t$.

 An indecomposable module $M\in\mathcal L_A$ is
\emph{$\ts{Ext}$-injective} in $\ts{add}\, \mathcal L_A$ if $\ts{Ext}_A^1(-,M)_{|\mathcal
  L_A}=0$ (see \cite{as}). This is the case if and only if
$\tau_A^{-1}M\not\in\mathcal L_A$.

The endomorphism algebra of the direct sum of the indecomposable
projective modules lying in $\mathcal L_A$ is called the \emph{left support of
  $A$}. If $A$ is laura with connecting component, then its left support is a product of tilted
algebras (see \cite[4.4, 5.1]{act}).

 An algebra $A$ is \emph{weakly shod} if the length of any
path in $\ts{ind}\, A$ from an injective to a projective is bounded
\cite{cl}. Also, $A$ is \emph{quasi-tilted} if its global
dimension $\ts{gl.dim}\, A$ is at most two and $\ts{ind}\, A=\mathcal L_A\cup \mathcal R_A$,
see \cite{hrs}.

\section{Covering functors}$\ $
\label{sec:s2}

A $k$-linear functor  $F\colon\mathcal E\to \mathcal B$ is a
\emph{covering functor} if (see \cite[3.1]{bog}):
\begin{enumerate}
\item $F^{-1}(x)\neq\emptyset$ for every
$x\in \mathcal B_o$.
\item  For every $x,y\in \mathcal E_o$, the two following $k$-linear maps
  are bijective:
\begin{equation}
\bigoplus\limits_{F(y')=F(y)}\mathcal E(x,y')\to \ \mathcal B(F(x),F(y)),\ 
\text{and} \ \ \bigoplus\limits_{F(x')=F(x)}\mathcal E(x',y)\to \
 \mathcal B(F(x),F(y))\ .\notag
\end{equation}
\end{enumerate}
Following \cite[\S 3]{gabriel}, $F$ is a \emph{Galois covering
  with group $G$} if there exists a group morphism $G\to \ts{Aut}(\mathcal
E)$ such that $G$ acts freely on $\mathcal E_o$,  $F\circ g=F$ for
every $g\in G$
and  the functor $\mathcal E/G\to \mathcal B$ induced by $F$  is an
isomorphism. We refer the reader to
\cite[3.1]{gabriel} for the definition of $\mathcal E/G$.  Galois coverings are
 covering functors.

If $F\colon \mathcal E\to \mathcal B$ is a covering functor, then $F$ defines an
adjoint pair
$(F_{\lambda},F_.)$ of  functors $F_{\lambda}\colon \ts{MOD}\, \mathcal E\to \ts{MOD}\, \mathcal B$
and $F_.\colon \ts{MOD}\, \mathcal B\to \ts{MOD}\, \mathcal E$ (see \cite[3.2]{bog}). The functor
$F_.$ is the \emph{pull-up functor} and
$F_{\lambda}$ is the \emph{push-down}. We recall
their construction: If $M\in\ts{MOD}\, \mathcal B$, then $F_.M=M\circ F^{op}$; if
$M\in \ts{MOD}\,\mathcal E$, then $F_{\lambda}M$ is the $\mathcal B$-module such that
$F_{\lambda}M(x)=\bigoplus\limits_{F(\widetilde x)=x}M(\widetilde x)$, for every
$x\in \mathcal B_o$. Both $F_{\lambda}$ and $F_.$ are exact.

Let $F\colon\mathcal E\to\mathcal
B$ be a covering functor between locally bounded $k$-categories.
We prove a few facts relative to $F$. Some are
easy to prove in case $F$ is a Galois covering. However, in
 general, the proofs are more complicated. This
can be explained by the following fact: $F^{op}\colon\mathcal
E^{op}\to\mathcal B^{op}$ is also a covering functor, and
$D\,F^{op}_{\lambda}\simeq F_{\lambda}\,D$ if $F$ is Galois. However,
this isomorphism no longer exists in the general case of covering
functors (see \cite[3.4]{bog}, for instance).

As a motivation for the results in this section, we start with
the following construction. 
We recall that the universal covering of a translation quiver $\Gamma$
was introduced in \cite[1.2]{bog} using a homotopy relation denoted as
$H$. We define $\hat{H}$ to be the smallest equivalence relation
containing $H$ and satisfying the following additional relation: Let $\alpha$
and $\beta$ be two arrows in $\Gamma$ having the same source and the
same target, then $\alpha$ and $\beta$ are equivalent for
$\hat{H}$. Using the construction of \cite[1.3]{bog} with respect to
the relation $\hat{H}$ we construct a covering of $\Gamma$ which we
call the \emph{generic covering of $\Gamma$}. It is an immediate
consequence of this definition and of \cite[1.3]{bog} that the generic
covering is a Galois covering and is a quotient of the universal
covering. They coincide if
$\Gamma$ has no multiple arrows (for example, if $\Gamma$ is the
Auslander-Reiten quiver of a representation-finite algebra).

The following proposition is mainly due to 
Riedtmann (see \cite[2.2]{r}).
\begin{prop}
\label{prop:well_behaved}
Let $A$ be a basic finite dimensional algebra.  Let $\Gamma$ be a
component of $\Gamma(\ts{mod}\, A)$. Let
$\pi\colon\widetilde{\Gamma}\to\Gamma$ be
  the generic covering. Then there
  exists a well-behaved functor $p\colon k(\widetilde{\Gamma})\to \ts{ind}\, \Gamma$. If,
  moreover, $\Gamma$ is generalised standard, then
  $p$ is a covering functor.
\end{prop}
\noindent{\textbf{Proof:}} The functor $p$ was constructed
in \cite[2.2]{r} for the stable part of the
Auslander-Reiten quiver of a self-injective representation-finite algebra. The
covering property was proved in \cite[2.3]{r}
under the same setting. The construction of $p$ was
generalised to any Auslander-Reiten component in
\cite[3.1]{bog}. It is
easily seen that the arguments given in \cite[2.3]{r} to prove that
$p$ is a covering functor apply to the
 case of generalised standard components.\sq

Note that if $\Gamma$ is a standard Auslander-Reiten component, then,
by definition, there exists a well-behaved functor
$k(\Gamma)\to\ts{ind}\,\Gamma$. In particular, any  covering of
translation quivers $p\colon
\Gamma'\to \Gamma$ gives rise to a well-behaved covering functor
$k(\Gamma')\to\ts{ind}\,\Gamma$ by composing the functors $k(p)\colon
k(\Gamma')\to k(\Gamma)$ and $k(\Gamma)\to\ts{ind}\,\Gamma$.

The results of this section will be applied to covering functors as in
\ref{prop:well_behaved}.
We now turn to the general situation where $F\colon\mathcal E\to \mathcal B$ is a
covering functor between locally bounded $k$-categories.

Since $F_{\lambda}$ and $F_.$ are exact, we still have an adjunction
at the level of derived categories. Here and in the sequel, $\mathcal D(\ts{MOD}\,
\mathcal E)$ and $\mathcal D^b(\ts{mod}\,\mathcal E)$ denote the derived category of $\mathcal
E$-modules and the bounded derived category of finite dimensional $\mathcal
E$-modules, respectively. The following lemma is
immediate. For a background on derived categories, we refer the reader
to \cite[Chap. III]{GM03}.
\begin{lem}
\label{lem:adjder}
  $F_{\lambda}$ and $F_.$ induce an adjoint pair $(F_{\lambda},F_.)$ of exact functors:
  \begin{equation}
    \xymatrix{
\mathcal D(\ts{MOD}\, \mathcal E) \ar@<+1ex>[r]^{F_{\lambda}}&
\mathcal D(\ts{MOD}\, \mathcal B) \ar@<+1ex>[l]^{F_.}&.
}\notag
  \end{equation}
Moreover $F_{\lambda}(\mathcal D^b(\ts{mod}\, \mathcal E))\subseteq\mathcal D^b(\ts{mod}\, \mathcal B)$.\hfill$\blacksquare$
\end{lem}

Let $x\in \mathcal E_o$. By condition $2$ in the
definition of a covering functor, $F$ induces a canonical isomorphism
$F_{\lambda}(\mathcal E(-,x))\xrightarrow{\sim} \mathcal B(-,F(x))$ of $\mathcal B$-modules (see
\cite[3.2]{bog}). In the sequel, we always identify these two
modules by means of this isomorphism. Using this identification we get
the following  result.
\begin{lem}
\label{lem:bijadj}
  Let $M\in\mathcal D(\ts{MOD}\,\mathcal E)$. Then $F_{\lambda}$ induces two linear
  maps for every $x\in \mathcal B_o$:
  \begin{align}
    & \varphi_M\colon \bigoplus\limits_{F(\widetilde x)=x} \mathcal D(\ts{MOD}\,\mathcal E)(M,\mathcal E(-,\widetilde
    x))\to \mathcal D(\ts{MOD}\, \mathcal B)(F_{\lambda} M,\mathcal 
    B(-,x))\ ,\notag\\
and\ & \psi_M\colon\bigoplus\limits_{F(\widetilde x)=x} \mathcal D(\ts{MOD}\,\mathcal E)(\mathcal E(-,\widetilde x),M)\to \mathcal D(\ts{MOD}\, B)(\mathcal
    B(-,x),F_{\lambda} M)\ .\notag
  \end{align}
These maps are functorial in $M$, and are bijective if $M$ is
quasi-isomorphic to a bounded complex of finite dimensional projective
modules (for example, if $\ts{gl.dim}\, \mathcal E<\infty$ and $M\in\ts{mod}\,\mathcal E$).
\end{lem}
\noindent{\textbf{Proof:}} 
 If $M=P[l]$ where $l\neq 0$ and $P$ is
  a projective $\mathcal E$-module, then $\varphi_M$ is bijective (because
  $\ts{Ext}_{\mathcal E}^{-l}(P,\mathcal E(-,\widetilde x))=0$). Also, 
if $M$ is an indecomposable projective
  $\mathcal E$-module, then $\varphi_M$ is bijective (because $F$ is a covering functor).
Finally, if $M\to M'\to M''\to M[1]$ is a triangle in $\mathcal D(\ts{MOD}\, \mathcal E)$,
  then $\varphi_M$, $\varphi_{M'}$ and $\varphi_{M''}$ are bijective as
  soon as two of them are so.
Consequently, $\varphi_M$ is bijective if $M$ is quasi-isomorphic to a
bounded complex of finite dimensional projective $\mathcal E$-modules.
The second map is handled similarly.\sq

In general, $F_{\lambda}$ does not commute with the
Auslander-Reiten translations. However, we
have the following.
\begin{lem}
\label{lem:dimtau}
  Let $X\in \ts{ind}\, \mathcal E$ be such that $F_{\lambda}X\in \ts{ind}\, \mathcal B$ and $\ts{pd}\, X<\infty$. Then
  $\dim \tau_{\mathcal E}X=\dim \tau_{\mathcal B}F_{\lambda}X$.
\end{lem}
\noindent{\textbf{Proof:}} Let $X\in\ts{mod}\,\mathcal E$ be any module. Let $P_1\to P_0\to X\to 0$ be a minimal
projective presentation in $\ts{mod}\, \mathcal E$. By \cite[3.2]{bog}, we deduce that $F_{\lambda} P_1\to
F_{\lambda} P_0\to F_{\lambda} X\to 0$ is a minimal projective presentation in $\ts{mod}\,\mathcal
B$. So we have exact sequences in $\ts{mod}\,\mathcal E^{op}$ and $\ts{mod}\,\mathcal B^{op}$,
respectively:
\begin{align}
  & 0\to \ts{Hom}_{\mathcal E}(X,\mathcal E)\to \ts{Hom}_{\mathcal E}(P_0,\mathcal E)\to \ts{Hom}_{\mathcal E}(P_1,\mathcal E)\to \ts{Tr}_{\mathcal E}X\to
  0\ ,\notag\\
\text{and}\ & 0\to \ts{Hom}_{\mathcal B}(F_{\lambda}X,\mathcal B)\to \ts{Hom}_{\mathcal B}(F_{\lambda}P_0,\mathcal B)\to \ts{Hom}_{\mathcal B}(F_{\lambda}
P_1,\mathcal B)\to \ts{Tr}_{\mathcal B}F_{\lambda} X\to
  0\ ,\notag
\end{align}
Let $X\in\ts{mod}\,\mathcal E$ be of finite projective
dimension, thus quasi-isomorphic to a bounded complex of
finite dimensional projective $\mathcal E$-modules.
The
bijections of \ref{lem:bijadj} imply that $\dim \ts{Hom}_{\mathcal E}(X,\mathcal
E)=\sum\limits_{x\in\mathcal E_o}\dim\ts{Hom}_{\mathcal E}(X,\mathcal
E(-,x))=\sum\limits_{x\in\mathcal B_o}\dim\ts{Hom}_{\mathcal B}(F_{\lambda}X,\mathcal B(-,x))=\dim
\ts{Hom}_{\mathcal B}(F_{\lambda}X,\mathcal B)$. Using the above exact sequences, we deduce
that $\dim\ts{Tr}_{\mathcal E}X=\dim\ts{Tr}_{\mathcal B}F_{\lambda}X$. Thus $\dim\tau_{\mathcal
  E}X=\dim \tau_{\mathcal B}F_{\lambda}X$ if both $X$ and $F_{\lambda}X$ are
indecomposable.\sq

\section{Standard laura algebras}$\ $
\label{sec:s3}

We now derive sufficient conditions for a laura algebra to be
standard. Weakly shod algebras are particular cases of laura
algebras. It is proved in \cite[\S 4]{cl} that if $A$ is weakly shod and not
quasi-tilted, then $A$ can be written as a one-point extension
$A=B[M]$ such that the connecting component of $A$ can be
recovered from $M$ and from the connecting components of $B$. This
motivates the following definition.
\begin{df}
 Let $A$ be a laura algebra with connecting components. An indecomposable projective $A$-module $P$ lying
    in a connecting component $\Gamma$ is a \emph{maximal projective}
    if it has an injective predecessor and no proper projective
    successor in $\ts{ind}\, A$. Furthermore, $A$ is a \emph{maximal
      extension} of $B$ if 
    there exists a maximal projective $P=eA$ such that
    $B=(1-e)A(1-e)$ and $A=B[M]$, where $M = \ts{rad} 
    P$.
\end{df}

By definition, a maximal projective
belongs to $\mathcal R_A$. In particular, by \cite[2.2]{al}, it is directed.
The notions of \emph{minimal injective} or \emph{maximal
  coextension}  are 
dual.
If $A$ is a tilted algebra which is the endomorphism algebra
of a regular tilting module, then it has neither maximal projective,
nor minimal injective (see \cite{R88}).

\begin{prop}\label{prop:ext-standard} 
 Let $A=B[M]$ be a maximal extension. Then $B$ is a product of laura
 algebras with connecting components. Moreover, if every
 connected component of $B$ is standard, then so is $A$.
\end{prop}
\noindent{\textbf{Proof:}}
By \cite{ac04}, every connected component of $B$ is a laura algebra.  
Let $P_m\in\ts{ind}\, A$ be the maximal projective such that $\ts{rad}
P_m=M$ and denote by $\Gamma$ the component of $\Gamma(\ts{mod}\, A)$ in which $P_m$
lies. So $P_m\in\mathcal R_A\cap\Gamma$. In particular, $P_m$ is directed. Note
that every proper predecessor of $P_m$ is an indecomposable $B$-module.\\

Let us prove the first assertion. If it is false, then a connected
component $B'$ of $B$ is quasi-tilted and not tilted (and, therefore,
quasi-tilted of canonical type). 
Since $A$ is connected, at least one indecomposable summand $M'$ of
$M$ lies in $\ts{ind}\, B'$. Assume first that
$M'$ is not directed. In particular, $M'\in\Gamma$ implies that
$M'\not\in\mathcal L_A\cup\mathcal R_A$. Therefore there is a
non-sectional path $M'\rightsquigarrow P$ in $\ts{ind}\, A$ with $P$
projective. If $P=P_m$, then there exists a non-sectional path
$M'\rightsquigarrow M''$ with $M''$ an indecomposable summand of
$M=\ts{rad} P_m$. This is impossible because $P_m$ is directed (see
\cite[Thm. 1 of \S 2]{HR93}). So $P\neq P_m$. By maximality of $P_m$, the path $M'\rightsquigarrow
P$ is a non-sectional path in $\ts{ind}\, B'$ ending at a projective. So
$M'\not\in\mathcal R_{B'}$. On the other hand, $M'\not\in\mathcal L_A$ means that there
exists a non-sectional path $I\rightsquigarrow M'$ in $\ts{ind}\, A$, where $I$
is injective. By maximality of $P_m$, this is a non-sectional path in
$\ts{ind}\, B'$. For the same reason, we have $\ts{Hom}_A(P_m,I)=0$, so that $I$
is injective as a $B'$-module. So 
$M'\not\in\mathcal L_{B'}\cup\mathcal R_{B'}$. This is impossible
because $B'$ is quasi-tilted. 
Therefore $M'$ is directed. Since $B'$ is quasi-tilted of canonical
type, the component $\Gamma'$ of $\Gamma(\ts{mod}\, B')$ containing $M'$ is either
the unique postprojective or the unique preinjective component (see
\cite[Prop. 4.3]{len-skow}).
Assume that $\Gamma'$ is the unique postprojective component of $\Gamma(\ts{mod}\,
B')$. Then $\Gamma'\subseteq \mathcal L_{B'}\backslash\mathcal R_{B'}$ (see
\cite[5.2]{CS96}). In particular,
there exists a non-sectional path $M'\rightsquigarrow P$ in $\ts{ind}\, B'$
with $P$ projective. Since $P_m$ is maximal, this is also
a non-sectional path in $\ts{ind}\, A$. Since $P$ is projective and since
$M'\in\Gamma$, we deduce that $P\in\Gamma$ and that the path is refinable to a
non-sectional path in $\Gamma(\ts{mod}\, A)$ and therefore in $\Gamma(\ts{mod}\, B')$
because $P_m$ is maximal. Consequently, $M'$ lies in the postprojective
component $\Gamma'$ of $\Gamma(\ts{mod}\, B')$ and is the starting point of a
non-sectional path in $\Gamma(\ts{mod}\, B')$ ending at a projective. This is
absurd. If $\Gamma'$ is the unique preinjective component of $\Gamma(\ts{mod}\,
B')$, then, using dual arguments, we also get a contradiction. Thus,
$B'$ is either tilted or not quasi-tilted.\\

Now, we assume that every connected component of $B$ is standard,
and prove that $A$ is standard. Later, in
\ref{rem:tilted_standard}, we shall see that, if
$A$ is tilted, then  its connecting components are standard. So
assume that $A$ is not tilted. Let $\Gamma$ be the connecting component of $\Gamma(\ts{mod}\, A)$ and
$\Gamma'$ be the disjoint union of the connecting components of the
Auslander-Reiten quivers of the connected components of $B$. We compare $\Gamma$ and
$\Gamma'$. More precisely, let $\mathcal
X$ be the full subquiver of $\Gamma$ with vertices those modules which are
not successors of $P_m$. So $\mathcal X$ is a full subquiver of $\Gamma(\ts{mod}\, B)$
stable under predecessors in
$\Gamma(\ts{mod}\, B)$, and it contains $\Gamma\backslash\mathcal R_A$. We claim that $\mathcal X$ is
contained in $\Gamma'$. We prove a series of assertions.

\textbf{(a) The left supports of $A$ and  $B$ coincide.}
Indeed, we have $\mathcal L_A\cap\ts{ind}\, B\subseteq \mathcal
L_B$ (see \cite[2.1]{ac04}). On the other hand, if $P\in \ts{ind}\, B$ is a
projective not lying in $\mathcal L_A$, then there is a non-sectional path
$I\rightsquigarrow P$ in $\ts{ind}\, A$ with $I$ injective. Since $P_m$ is maximal, this is a
non-sectional path in $\ts{ind}\, B$. For the same reason,
$\ts{Hom}_A(P_m,I)=0$, so that $I$ is injective as a $B$-module. So
$P\not\in\mathcal L_B$. Thus $A$ and $B$ have the same
left support.

\textbf{(b) Let $P\neq P_m$ be a projective lying in $\Gamma$. Then
  $P\in\Gamma'$.} Indeed, if
there exists a path $I\rightsquigarrow P$ in $\Gamma$ with $I$ injective,
then the maximality of $P_m$ implies that this path lies entirely in
$\ts{ind}\, B$ and starts in an injective $B$-module. So $P\in \Gamma'$. If there is no such path,
then $P\in\mathcal L_A\cap\Gamma$. So $P$ lies in a connecting component of one of
the components of the left support of $A$, which is also the left
support of $B$. From \cite[5.4]{ac}, we deduce that $P$ lies in $\Gamma'$.

\textbf{(c) Let $X\in\mathcal X$. There exists $m\geqslant 0$ such that $\tau_A^mX\in\Gamma'$.}
By assumption on $X$, we have $\tau_BX=\tau_AX$.
Assume first that $\tau_A^mX=P$ for some $m\geqslant 0$ and some
projective $P$. So $P\neq P_m$. From (b), we get that $P\in\Gamma'$. Now
assume that $X$ is left stable and 
non-periodic. If $X\in\mathcal R_A$, there exists $l\geqslant 0$ such that
$\tau_A^lX$ is $\ts{Ext}$-projective in $\mathcal R_A$. Since $X$ is left stable,
we deduce that $\tau_A^{l+1}X\in\Gamma\backslash\mathcal R_A$. So assume
that $X\in\Gamma\backslash\mathcal R_A$.
Since $A$ is laura, there exists $m$ such that
$\tau_A^mX\in\Gamma\cap\mathcal L_A$. So $\tau_A^mX$ lies in one of the connecting
components of the left support of $A$. So $\tau_A^mX\in\Gamma'$ because the left
supports of $A$ and $B$ are equal. Finally, assume that $X$ is periodic.
Then there exists
a projective module $P \in \Gamma$, a periodic direct
summand $Y$ of $\ts{rad} P$, and a path $Y \rightsquigarrow
X$ in $\Gamma\backslash\mathcal R_A$, and therefore in $\Gamma(\ts{mod}\, B)$. Since $Y$ is
periodic, then $P\neq P_m$ (otherwise $P_m$
would be a proper successor of itself). Since $P\in\Gamma'$, we have $Y\in\Gamma'$ and therefore $X\in\Gamma'$.

\textbf{(d) $\mathcal X$ is contained in $\Gamma'$.} Indeed, we already know
that $\mathcal X$ is a full subquiver of $\Gamma(\ts{mod}\, B)$. Also, we proved that
for every $X\in\mathcal X$, there exists $m\geqslant 0$ such that $\tau_A^mX=\tau_B^mX\in\Gamma'$. So $\mathcal
X$ is contained in $\Gamma'$.

We now show that $\Gamma$ is standard. By hypothesis, there exists a
well-behaved functor $\varphi \colon k(\Gamma') \to \ts{ind}\, \Gamma'$. 
Since $\mathcal X$ is a full subquiver of $\Gamma'$ stable under predecessors in
$\Gamma(\ts{mod}\, B)$, there exists a well-behaved functor
$\psi\colon k(\mathcal Y)\to \ts{ind}\,\Gamma$ where $\mathcal Y$ is a full subquiver of
$\Gamma$ such that:
\begin{enumerate}
\item $\mathcal Y$ contains $\mathcal X$.
\item $\mathcal Y$ is stable under predecessors in $\Gamma(\ts{mod}\, A)$.
\item $\psi$ and $\varphi$ coincide on $\mathcal X$.
\item $\mathcal Y$ is maximal for these properties.
\end{enumerate}
We show that $\mathcal Y=\Gamma$. Assume that $\mathcal Y\neq\Gamma$.
 Since $\mathcal Y$ contains $\mathcal X$, it contains
$\Gamma\backslash\mathcal R_A$, so there exists a source $X$ in $\Gamma\backslash\mathcal Y$.
If $X$ is projective, then $X=P_m$. So $\psi$ is
defined on every indecomposable summand $Y$ of $\ts{rad} P_m$. 
Set $\psi(X)=P_m$. Let $\alpha_1\colon X_1\to
P_m,\ldots,\alpha_t\colon X_t\to P_m$ be the
arrows ending at $X$. Then $X_1\oplus\cdots \oplus X_t=\ts{rad} P_m$, and
let $\psi(\alpha_i)$ be the inclusion $X_i\hookrightarrow P_m$. If
$X$ is not
projective, then the mesh ending at $X$ has the following shape:
\begin{equation}
  \xymatrix@=2pt{
&&&& X_1 \ar@{->}[rrrrd]^{v_1} &&&& \\
\tau_A X\ar@{->}[rrrru]^{u_1} \ar@{->}[rrrrd]_{u_n} &&&& \vdots &&&& X\\
&&&&X_n\ar@{->}[rrrru]_{v_n}&&&&&.
}\tag{$\star$}
\end{equation}
Since $X$ is a source of $\Gamma \backslash \mathcal Y$, then $\psi$ is
already defined on the full subquiver of the mesh consisting of all
vertices  except $X$. In particular, the following map is right
minimal almost split:
\begin{equation}
  \begin{bmatrix}
    \psi(u_1)&\ldots&\psi(u_n)
  \end{bmatrix}^t\colon \tau_AX\to X_1\oplus\cdots \oplus X_n\ .\notag
\end{equation}
Let $\psi(X)=X$, and $
\begin{bmatrix}
  \psi(v_1)&\ldots&\psi(v_n)
\end{bmatrix}\colon X_1\oplus\cdots \oplus X_n\to X$ be the cokernel of the
above map, following \cite[3.1,
Ex. b]{bog}. Clearly, this construction contradicts the maximality of $\mathcal
Y$. So $\mathcal Y=\Gamma$ and there exists 
a well-behaved functor $\psi \colon k(\Gamma) \to \ts{ind}\, \Gamma$ which
is the identity on objects. The arguments in the proof of
\cite[5.1]{bog} show that this is an isomorphism. So $\Gamma$ is
standard. \sq

Since weakly shod algebras are  laura, it
makes sense to speak of weakly shod algebras with connecting
components. We have the following corollary.
\begin{cor}
 Let $A$ be a (connected) weakly shod algebra with connecting components, then $A$ is standard.
 \end{cor}
\noindent{\textbf{Proof:}} By \cite[3.3]{al}, there exists a
sequence of full convex subcategories   \begin{equation}
    C=A_0\subsetneq A_1\subsetneq \cdots\subsetneq A_m=A \notag
\end{equation} with $C$ tilted and, for each $i \geqslant 0$, the
algebra $A_{i+1}$ is a maximal extension of $A_i$. The result follows
from \ref{prop:ext-standard} and induction because $C$
is standard (see \ref{rem:tilted_standard} below).\sq

The preceding result motivates the following definition, inspired from
\cite[2.3]{al}.

\begin{df} Let $A$ be a laura algebra. We say that $A$ admits a
    \emph{maximal filtration} if there exists a sequence 
\begin{equation}
    C=A_0\subsetneq A_1\subsetneq \cdots\subsetneq A_m=A \tag{$f$}
\end{equation} of full convex subcategories with $C$
a product of representation-finite algebras and, for each $i\geqslant 0$, the algebra
$A_{i+1}$ is a maximal extension, or a maximal coextension, of $A_i$.
\end{df}

\begin{cor}
Let $A$ be a laura algebra admitting a maximal filtration  ($f$):
\begin{enumerate}
\item[(a)] If $C$ is a product of standard representation-finite
  algebras, then $A$ is standard.
\item[(b)] If the Auslander-Reiten quiver of every connected component of $C$
  is simply connected, then $A$ is standard.	
\item[(c)] If $\ts{HH}^1(A)=0$, then $A$ is standard.
\end{enumerate}
\end{cor}

\noindent{\textbf{Proof:}}
\noindent  Statement $(a)$ follows directly from \ref{prop:ext-standard}.\\
\noindent  $(b)$ This follows from \ref{prop:ext-standard} and the fact that if a
representation-finite connected algebra $C$ has $\ts{HH}^1(C)=0$, or
equivalently, if its Auslander-Reiten
quiver is simply connected, then $C$ is standard  \cite[4.2]{bl}. 

\noindent $(c)$ We use induction on the length $m$ of a maximal
filtration. If $m=0$, then $A$ is representation-finite and the result
follows from \cite[4.2]{bl}. Assume that $m\geqslant 1$ and that
the statement holds for algebras admitting maximal filtrations of
length less than $m$. Without loss of generality, we may assume that
$A=A_{m-1}[M]$ is a maximal extension. We claim that $
  \ts{Ext}^1_{A_{m-1}}(M,M)=0$. Indeed, if this is not the case, then
there exists an indecomposable summand $N$ of $M$ such that $
  \ts{Ext}^1_{A_{m-1}}(M,N)\not=0$. Write $M \simeq N \oplus N'$ and let $P$
be the indecomposable projective such that $M= \ts{rad} P$. Then
$N'$ is a submodule of $P$ and $L=P/N'$ is indecomposable. By
\cite[III.2.2, (a)]{hrs} we have $\ts{id}\, L \geqslant 2$. But
this contradicts the fact that $L\in \mathcal R_A$ because it is a successor of
the maximal projective $P$. So $\ts{Ext}^1_{A_{m-1}}(M,M)=0$.
Applying \cite[5.3]{happel}, the exact sequence 
\begin{equation}
\ts{HH}^1(A) \to \ts{HH}^1(A_{m-1}) \to \ts{Ext}_{A_{m-1}}^1(M,M)\notag
\end{equation} 
yields $\ts{HH}^1(A_{m-1})=0$. By the induction
hypothesis, $A_{m-1}$ is standard. By
\ref{prop:ext-standard}, so is $A$. \sq 

\begin{exs}
\label{ex1}
\begin{itemize}
 \item [$(a)$] Let $A$ be the radical-square zero algebra given by the
   quiver
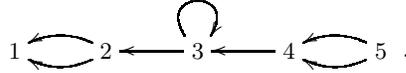
\vspace{2mm}
\begin{equation}
\xymatrix{1&2\ar@/^/[l] \ar@/_/[l]&3\ar@(ul,ur)\ar[l] & 4 \ar[l] &
  5\ar@/^/[l] \ar@/_/[l]}\ .\notag
\end{equation}
This is a laura algebra (see \cite[2.3]{ac}).
Here and in the sequel, we denote by $P_x$, $I_x$ and $S_x$ the
indecomposable projective, the indecomposable injective, and the
simple module corresponding to the vertex $x$, respectively. Clearly $P_1$ is
maximal projective and $I_5$ is minimal injective.
 Letting
$C$ be the full convex subcategory with objects
$\{2,3,4\}$ we see that 
\begin{equation}
 C\subsetneq [S_4\oplus S_4]C \subsetneq A\notag
\end{equation}
is a maximal filtration. Since $C$ is standard, so is $A$.
Its connecting component is drawn below:
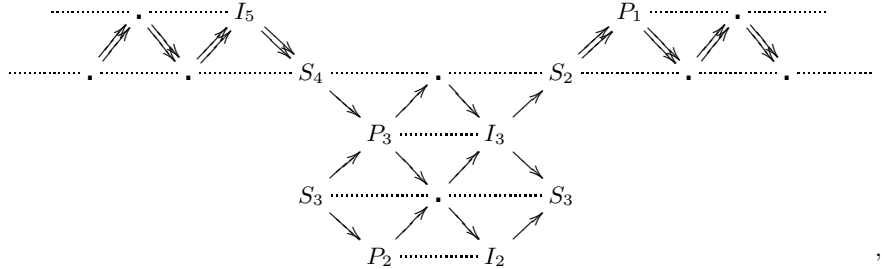
\begin{equation}
  \xymatrix@=10pt{
& \ar@{.}[rr]& &  \centerdot \ar@<1pt>[rd] \ar@<-1pt>[rd] \ar@{.}[rr]& & 
I_5 \ar@<1pt>[rd] \ar@<-1pt>[rd] & & & & & & 
P_1 \ar@<1pt>[rd] \ar@<-1pt>[rd] \ar@{.}[rr] & & 
\centerdot \ar@<1pt>[rd] \ar@<-1pt>[rd] \ar@{.}[rr]& & &\\
\ar@{.}[rr]&& \centerdot \ar@<1pt>[ru] \ar@<-1pt>[ru] \ar@{.}[rr] & & 
 \centerdot \ar@<1pt>[ru] \ar@<-1pt>[ru] \ar@{.}[rr] & &
S_4 \ar[rd] \ar@{.}[rr] &  & \centerdot  \ar[rd] \ar@{.}[rr]& 
& 
S_2 \ar@<1pt>[ru] \ar@<-1pt>[ru] \ar@{.}[rr] & &
\centerdot \ar@<1pt>[ru] \ar@<-1pt>[ru] \ar@{.}[rr] & &
 \centerdot \ar@{.}[rr]&&\\
&&&&&&& 
P_3 \ar[ru] \ar@{.}[rr] \ar[rd]& &
I_3 \ar[ru] \ar[rd] & & & & & &&\\
&&& & & & S_3 \ar[ru] \ar@{.}[rr] \ar[rd]& & 
\centerdot \ar[ru] \ar[rd] \ar@{.}[rr] & & S_3 &&\\
&&& & & & &
P_2 \ar[ru] \ar@{.}[rr] 
& & 
I_2 \ar[ru]&&&&&&&,
}\notag
\end{equation}
where the two copies of $S_3$ are identified.
 \item [$(b)$] Let $B,C$ be products of standard laura algebras, and $A$ an
   articulation of $B,C$ (in the sense of \cite{ds}). Then $A$ is
   laura with connecting components (see \cite{ds}). Using  \cite[3.9]{ds}
   it is easy to check that $A$ is standard.
\end{itemize}
\end{exs}

The section motivates the following questions.
\begin{pb}
  Which laura algebras admit maximal filtrations?
\end{pb}
\begin{pb}
  Assume that $A$ is a laura algebra which does not admit a maximal
  filtration. If $\ts{HH}^1(A)=0$, do we have that $A$ is standard?
\end{pb}
\section{Tilting modules of the first kind with respect to covering
  functors}$\ $
\label{sec:s4}

For tilting theory, we refer to \cite{ass}.
Let $B$ be a product of tilted algebras and $n$ be the
rank of its Grothendieck group.
In \cite[Cor. 4.5]{lemeur2}, it is proved that tilting modules are of the first kind with respect to any Galois
covering of $B$. More precisely, let $F\colon \widetilde B\to B$ be a Galois covering
with group $G$, where $\widetilde B$ is locally bounded. Denote by $\mathcal T$ the
class of complexes $T\in\mathcal D^b(\ts{mod}\, B)$ such that:
\begin{enumerate}
\item $T$ is multiplicity-free and has $n$
  indecomposable summands.
\item $\mathcal D^b(\ts{mod}\, B)(T,T[i])=0$ for every $i\geqslant 1$ (so $T$ is a
  \emph{silting complex} in the sense of \cite{KV88}).
\item $T$ generates the triangulated category $\mathcal D^b(\ts{mod}\, B)$.
\end{enumerate}
Any multiplicity-free tilting module lies in $\mathcal T$.
It was proved in \cite[\S 4]{lemeur2} that for any $T\in\mathcal T$ and for any
indecomposable summand
$X$ of $T$, there exists $\widetilde X\in\mathcal D^b(\ts{mod}\, \widetilde B)$ such that:
\begin{enumerate}
\item $F_{\lambda}\widetilde X\simeq X$.
\item $^g\widetilde X\not\simeq \,^h\widetilde X$ for $g\neq h$.
\item If $Y\in\mathcal D^b(\ts{mod}\, \widetilde B)$ is such that $F_{\lambda}Y\simeq X$, then
  $Y\simeq\,^g\widetilde X$ for some $g\in G$.
\end{enumerate}
Given $T\in\mathcal
T$ and  an indecomposable summand $X$ of $T$, we fix $\widetilde X\in\mathcal
D^b(\ts{mod}\, \widetilde B)$ arbitrarily such that $F_{\lambda}\widetilde X\simeq X$.

For later reference, we recall some facts. The following result was proved in
\cite[Cor. 4.5, Prop. 4.6, Lem. 4.8]{lemeur2}.
\begin{lem}
\label{lem1}
Let $F\colon \widetilde B\to B$ be a Galois covering with group $G$. Let
$T\in\ts{mod}\, B$ be a multiplicity-free tilting module. Let
  $T=T_1\oplus\cdots\oplus T_n$ be such that $T_1,\ldots, T_n$ are 
  indecomposable. For every $i$, there exists $\widetilde T_i\in \ts{ind}\,\widetilde B$
  such
  that $F_{\lambda}\widetilde T_i=T_i$. Moreover:
  \begin{enumerate}
\item[(a)] $^g\widetilde T_i\not\simeq \,^h\widetilde T_j$ for $(g,i)\neq(h,j)$.
\item[(b)] $\ts{pd}\,\widetilde T_i\leqslant 1$ for every $i$.
\item[(c)] $\ts{Ext}_{\widetilde B}^1(\,^g\widetilde
  T_i,\,^h\widetilde T_j)=0$ for every $g,h\in G,\, i,j\in\{1,\ldots,n\}$.
\item[(d)] For every indecomposable projective $\widetilde B$-module $P$, there
  exists an exact sequence $0\to P\to T^{(1)}\to T^{(2)}\to 0$ with
  $T^{(1)},T^{(2)}$ in $\ts{add}\,\{\,^g\widetilde T_i\ |\ g\in G,\ i\in\{1,\ldots,n\}\}$.
  \end{enumerate}
\null\hfill$\blacksquare$
\end{lem}
We need similar facts about covering functors which need
not be Galois. Thus we prove the following result.
\begin{prop}
  \label{prop:cov}
Let $F\colon \widetilde B\to B$ be a Galois covering with group $G$, where $\widetilde
B$ is locally bounded. With the above setting, let $p\colon\widetilde B\to B$
be a covering functor
such that $F(x)=p(x)$ for every $x\in \widetilde B_o$. Let $T\in\mathcal T$ and $X$ be an indecomposable
summand of $T$. Then:
\begin{enumerate}
\item[(a)] There exists an isomorphism $p_{\lambda}(\,^g\widetilde
  X)\xrightarrow{\sim} X$, for every $g\in G$.
\item[(b)] If $L\in\mathcal D^b(\ts{mod}\, \widetilde B)$ is such that $p_{\lambda}L\simeq X$, then
  $L\simeq\,^g\widetilde X$ for some $g\in G$.
\item[(c)] For every $L\in\mathcal D^b(\ts{mod}\, \widetilde B)$, the following maps induced by
  $p_{\lambda}$ and by the isomorphisms of (a) are linear
  bijections:
  \begin{align}
&    \varphi_{X,L}\colon \bigoplus\limits_{g\in G}\mathcal D^b(\ts{mod}\, \widetilde B)(\,^g\widetilde
    X,L)\xrightarrow{\sim} \mathcal D^b(\ts{mod}\, B)(X,p_{\lambda}L)\,,\notag \\
and\ &\psi_{X,L}\colon \bigoplus\limits_{g\in G}\mathcal D^b(\ts{mod}\, \widetilde B)(L,\,^g\widetilde
    X)\xrightarrow{\sim} \mathcal D^b(\ts{mod}\, B)(p_{\lambda}L,X)\,.\notag
  \end{align}
\end{enumerate}
\end{prop}
In order to prove the proposition, we need the following lemma. In
case $p$ is a Galois
covering, the lemma was proved in \cite[Lems. 4.2, 4.3]{lemeur2} (see also
\cite[Lems. 3.2, 3.3]{lemeur1}). For simplicity, we write
$\ts{Hom}(X,Y)$ for the space of morphisms in the derived category.
\begin{lem}
\label{lem:lift}
Let $T,T'\in\mathcal T$ be such that \ref{prop:cov} holds true for
$T$ and for $T'$.
Consider a triangle in $\mathcal D^b(\ts{mod}\, B)$:
\begin{equation}
  X\to \bigoplus\limits_{i=1}^tX_i'\to Y\to X[1]\ ,\tag{$\Delta$} 
\end{equation}
where $X\in\ts{add}\, T$ and $X_1',\ldots,X_t'$ are indecomposable
summands of $T'$. Assume that
$\ts{Hom}(Y,X'_i[1])=0$ for all $i$ (we do not assume that $Y\in\ts{add}\, T$ or
$Y\in\ts{add}\, T'$). Then for every $g\in G$, there exist $\widetilde Y\in\mathcal
D^b(\ts{mod}\,\widetilde B)$ and $g_1,\ldots,g_t\in G$ such that the triangle $\Delta$ is
isomorphic to the image under $p_{\lambda}$ of a triangle in $\mathcal D^b(\ts{mod}\,
\widetilde B)$ as follows:
\begin{equation}
  ^g\widetilde X\to \bigoplus\limits_{i=1}^t\,^{g_i}\widetilde X_i'\to \widetilde Y\to \,^g\widetilde X[1]\ .\notag
\end{equation}
Dually, consider a triangle in $\mathcal D^b(\ts{mod}\, B)$:
\begin{equation}
  Y\to \bigoplus\limits_{i=1}^tX'_i\to X\to Y[1]\ ,\tag{$\Delta'$}
\end{equation}
where $X\in\ts{add}\, T$ and $X_1',\ldots,X'_t$ are indecomposable summands
of $T'$. Assume that
$\ts{Hom}(X'_i,Y[1])=0$ for all $i$. Then for every $g\in G$, there exist $\widetilde Y\in\mathcal
D^b(\ts{mod}\, \widetilde B)$ and
$g_1,\ldots,g_t\in G$ such that the triangle $\Delta'$ is
isomorphic to the image under $p_{\lambda}$ of a triangle in $\mathcal D^b(\ts{mod}\,\widetilde
B)$ as follows:
\begin{equation}
  \widetilde Y\to \bigoplus\limits_{i=1}^t\,^{g_i}\widetilde X_i'\to \,^g\widetilde X\to \widetilde Y[1]\ .\notag
\end{equation}
\end{lem}
\noindent{\textbf{Proof:}} The proofs of \cite[Lems. 4.2,
4.3]{lemeur2} use the following key property of a Galois covering
$F\colon\widetilde B\to B$ with group $G$. Given $L,M\in\mathcal
D^b(\ts{mod}\,\widetilde B)$, 
we have linear bijections induced by $F_{\lambda}$:
\begin{equation}
\bigoplus\limits_{g\in G}\ts{Hom}(\,^gL,M)\xrightarrow{\sim}
\ts{Hom}(F_{\lambda}L,F_{\lambda}M)\ 
\text{and}\  \bigoplus\limits_{g\in G}\ts{Hom}(L,\,^gM)\xrightarrow{\sim}
\ts{Hom}(F_{\lambda}L,F_{\lambda}M)\ .\notag
\end{equation}
Of course, these bijections no longer exist for a covering functor
which is not Galois.
However, using our hypothesis that \ref{prop:cov} holds true for $T$
and for $T'$, it
is easy to check that the proofs of \cite[Lems. 4.2, 4.3]{lemeur2}
still work in the present case. Whence the lemma.\sq

\noindent{\textbf{Proof of \ref{prop:cov}:}
We proceed in several steps.

\textbf{Step $1$: If $T=B$, then \ref{prop:cov} holds true.}
 The following facts follow from the
definition of  covering functors (see also \cite[3.2]{bog}):
\begin{enumerate}
\item $Y\in\mathcal D^b(\ts{mod}\, \widetilde B)$ is a projective module if and only if
  $p_{\lambda}Y$ is a projective module.
\item $p_{\lambda}\left(\widetilde B(-,x)\right)\simeq F_{\lambda}\left(\widetilde
    B(-,x)\right)\simeq B(-,F(x))=B(-,p(x))$ for every $x\in\widetilde B_o$.
\item $^g\widetilde B(-,x)=\widetilde B(-,gx)$ for every $x\in\widetilde B_o$ and every
  $g\in G$.
\end{enumerate}
Therefore \ref{prop:cov} holds true for $T=B$.\\

Given an object $X$ in a triangulated category, we write
$\left<X\right>$ for the smallest additive full subcategory containing $X$ 
which is 
stable under direct summands and shifts (in both
directions).\\

\textbf{Step $2$: If $T,T'\in\mathcal T$ are such that $T'\in\left<T\right>$, then
  \ref{prop:cov} holds true for $T$ if and only if it does
  for $T'$.} This follows directly from the compatibility of $p_{\lambda}$
with the shift.\\

For the next step, consider the following situation.
Assume that $T,T'\in\mathcal T$ are such that:
\begin{enumerate}
\item $T=M\oplus \overline T$, where $M$ is indecomposable.
\item $T'=M'\oplus \overline T$, where $M'$ is indecomoposable.
\item There exists a non-split triangle $\Delta:\
  M\xrightarrow{u}E\xrightarrow{v} M'\to M[1]$ where $u$ is a
  left minimal $\ts{add}\, \overline T$-approximation and $v$ is a
  right minimal
  $\ts{add}\,\overline T$-approximation.
\end{enumerate}

\textbf{Step $3$: If $T,T'\in\mathcal T$ are as above, then
\ref{prop:cov} holds true for $T$
  if and only if it does for $T'$.}
 We prove that the condition is
necessary. Clearly, it suffices to prove that the assertions (a), (b),
and (c) of \ref{prop:cov} are true for $M'$. For simplicity, we identify $p_{\lambda}(\,^g\widetilde X)$ and
$X$ via the isomorphism used to define $\varphi_{X,-}$ and
$\psi_{X,-}$ for every indecomposable summand $X$ of $T$ and
$g\in G$.

Let $E=\bigoplus\limits_{i=1}^tE_i$ with the $E_i$ indecomposable. Recall from \cite[Lem. 4.4]{lemeur2}
that $\Delta$ is isomorphic to
the image under $F_{\lambda}$ of a triangle $\widetilde{\Delta}$ in $\mathcal D^b(\ts{mod}\, \widetilde B)$:
\begin{equation}
  \widetilde M
\xrightarrow{\widetilde u}
\bigoplus\limits_{i=1}^t
\,^{g_i}\widetilde E_i
\xrightarrow{\widetilde v}
\,^{g_0}\widetilde M'
\to
  \widetilde M[1]\ ,\tag{$\widetilde{\Delta}$}
\end{equation}
for some $g_0,g_1,\ldots,g_t\in G$. Moreover, $\widetilde u$ is a left
minimal $\ts{add}\,\mathcal 
X$-approximation and $\widetilde v$ is a right minimal $\ts{add}\,\mathcal
X'$-approximation, where $\mathcal X$ and $\mathcal X'$ are the following
full subcategories of $\mathcal D^b(\ts{mod}\, \widetilde B)$:
\begin{enumerate}
\item[-] $\mathcal X=\{\,^g\widetilde X\ |\ \text{$g\in G$, $X$ an indecomposable
    summand of $T$ and $^g\widetilde X\not\simeq \widetilde  M$}\}$.
\item[-] $\mathcal X'=\{\,^g\widetilde X\ |\ \text{$g\in G$, $X$ an indecomposable
    summand of $T'$ and $^g\widetilde X\not\simeq \widetilde M'$}\}$.
\end{enumerate}
Fix $g\in G$. Since \ref{prop:cov} holds true for $T$,
we apply \ref{lem:lift} to construct a triangle $\widetilde{\Delta}'\colon\
^g\widetilde M\xrightarrow{\widetilde u'}\bigoplus\limits_{i=1}^t\,^{g_i'}\widetilde
E_i\xrightarrow{\widetilde v'} Z_g\to \,^g\widetilde M[1]$ whose image under
$p_{\lambda}$ is isomorphic to $\Delta$. In particular,
$p_{\lambda}(Z_g)\simeq M'$. For simplicity, assume that $\Delta$ is
equal to the image of $\widetilde{\Delta}$ under $p_{\lambda}$, and
set $\widetilde E'=\bigoplus\limits_{i=1}^t\,^{g_i'}\widetilde E_i$. Let us
prove that $Z_g\simeq \,^{gg_0}\widetilde M'$. It suffices to prove that
$\widetilde{\Delta}'$ and $^g\widetilde{\Delta}$ are isomorphic. For this
purpose, we only need to prove that $\widetilde u'$ is a left minimal
$\ts{add}\,\,^g\mathcal X$-approximation. Let $f\colon \,^g\widetilde
M\to\,^{g'}\widetilde Y$ be non-zero,
where $Y$ is an indecomposable summand of $T$ such that $^{g'}\widetilde Y\in
\,^g\mathcal X$. Since $\varphi_{M,\widetilde M}$ is bijective
and since $\ts{End}(M)=k$, we have $Y\in \ts{add}\,\overline T$. So we
have a factorisation of $p_{\lambda}(f)$ by $u=p_{\lambda}(\widetilde u')$:
\begin{equation}
  \xymatrix{
M \ar[r]^{u} \ar[rd]_{p_{\lambda}(f)} & E \ar[d]^{f'}& \\
& Y & .
}\notag
\end{equation}
Since $\psi_{Y,\widetilde E_i}$ is bijective for every $i$, we have
$f'=\sum\limits_{h\in G}p_{\lambda}(f'_h)$, where $(f'_h)_h\in
\bigoplus\limits_{h\in G}\ts{Hom}(\,\widetilde E,\,^h\widetilde Y)$. So
$p_{\lambda}(f-f'_{g'}\widetilde u')-\sum\limits_{h\neq g'}p_{\lambda}(f'_h\widetilde u')=0$. Using
\ref{prop:cov}, we get $f=f'_{g'}\widetilde u'$. Hence
$\widetilde u'$ is a left $\ts{add}\,\,^g\mathcal X$-approximation. On the other hand,
$\widetilde u'$ is left minimal because $u=p_{\lambda}(\widetilde u')$ is
left minimal and $p_{\lambda}$ is exact. As explained above,
these facts imply that $Z_g\simeq\,^{gg_0}\widetilde M'$. So
$p_{\lambda}(\,^g\widetilde M')\simeq M'$, for every $g\in G$.

Let $Y\in\mathcal D^b(\ts{mod}\, \widetilde B)$. Using the triangles $^g\widetilde{\Delta}$ ($g\in G$) and
using that \ref{prop:cov} holds true for $T$, the maps
$\varphi_{M',Y}$ and $\psi_{M',Y}$ are bijective (recall that
$\ts{Hom}$-functors are cohomological).

Finally, if $Y\in\mathcal D^b(\ts{mod}\, \widetilde B)$, and if $f\colon p_{\lambda}Y\to M'$ is an
isomorphism, then $f=\sum\limits_{g\in G}p_{\lambda}(f_g)$
with $(f_g)_g\in \bigoplus\limits_{g\in G}\ts{Hom}(Y,\,^g\widetilde M')$. Since
$p_{\lambda}Y$ and $M'$ are indecomposable, there exists $g_1\in G$
such that $p_{\lambda}(f_{g_1})$ is an isomorphism. Since $p_{\lambda}$
is exact, we deduce that $f_{g_1}\colon Y\to\,^{g_1}\widetilde M'$ is an
isomorphism. This finishes the proof of the assertion: \ref{prop:cov} holds true for
$T'$ if it holds true for $T$. The converse implication is proved using
similar arguments.\\

\textbf{Step $4$: If $T\in\mathcal T$, then \ref{prop:cov} holds true.}
 This follows directly from the three preceding steps, and from
 \cite[Prop. 3.7]{lemeur2}.
\sq

\begin{ex}
Let $B=kQ$ be the path algebra of the Kronecker quiver
$\xymatrix{1\ar@<2pt>[r]^a\ar@<-2pt>[r]_b&2}$. There is a Galois covering
$F\colon \widetilde B\to B$ with group $\mathbb{Z}/2\mathbb{Z}=\{1,\sigma\}$, where
$\widetilde B=k\widetilde Q$ is the path algebra of the following quiver:
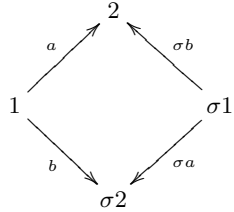
\begin{equation}
  \xymatrix{
&2&\\
1\ar[ru]^a \ar[rd]_b&&\sigma 1\ar[lu]_{\sigma b} \ar[ld]^{\sigma a}\\
&\sigma 2&
}\notag
\end{equation}
and where $F$ is the functor such that $F(\sigma^i \alpha)=\alpha$ for
every arrow $\alpha$ and every $i\in\{0,1\}$. On the other hand, there
is a covering functor $p\colon \widetilde B\to B$ such that $p(b)=p(\sigma
b)=b$, $p(a)=a$ and $p(\sigma a)=a+b$.  The $B$-module $T=e_2B\oplus\tau_B^{-1}(e_1B)$
 is tilting. One checks easily that $F_{\lambda}(e_2\widetilde B)=e_2B$,
$F_{\lambda}(\tau_{\widetilde B}^{-1}(e_1\widetilde B))=\tau_B^{-1}(e_1 B)$ and that
$p_{\lambda}(e_2\widetilde B)\simeq e_2B$,
$p_{\lambda}(\tau_{\widetilde B}^{-1}(e_1\widetilde B))\simeq\tau_B^{-1}(e_1 B)$.
\end{ex}

\section{Coverings of left sections}$\ $
\label{sec:s5bis}

Let $A$ be a basic finite dimensional $k$-algebra, $\Gamma$ a
component of $\Gamma(\ts{mod}\, A)$,
$\pi\colon\widetilde{\Gamma}\to\Gamma$  
a Galois covering of translation quivers with group $G$ such that
there exists a well-behaved functor $p\colon
k(\widetilde{\Gamma})\to\ts{ind}\,\Gamma$. A \emph{left section} (see
\cite[2.1]{assem}) in $\Gamma$ is a full subquiver $\Sigma$ such that: $\Sigma$ is
acyclic; it is convex in
$\Gamma$; and, for any $x\in \Gamma$, predecessor in $\Gamma$ of some $y\in\Sigma$,  there exists a unique $n\geqslant 0$ such
that $\tau^{-n}x \in \Sigma$. Assume that $\Sigma$ is a left section in
$\Gamma$ and let $B=A/\ts{Ann}\,\Sigma$. In this section, we construct a covering
functor $F\colon \widetilde B\to B$ associated to $p$ and  a functor
$\varphi\colon k(\widetilde{\Gamma})\to\ts{mod}\,\widetilde B$. Both $F$ and $\varphi$ are essential
in the proofs of Theorems~\ref{thm2} and \ref{thm1}.

By \cite[Thm. A]{assem}, the algebra $B$ is a full convex subcategory
of $A$ and a product of tilted
algebras and the components of $\Sigma$ form complete slices in the
connecting components of the connected components of $B$. 
Recall from
\cite{R88} that a connected  algebra $B'$ is tilted if and only if its
Auslander-Reiten quiver contains a so-called complete slice $\Sigma'$, that is, a
class of indecomposable $B'$-modules such that: ($1$) 
$U=\bigoplus\limits_{X\in\Sigma'}X$  is sincere (that is,
$\ts{Hom}_{B'}(P,U)\neq 0$ for any projective $B'$-module $P$); ($2$)
$\Sigma'$ is convex in $\ts{ind}\,B'$; ($3$) If $0\to L\to M\to
N\to 0$ is an almost split sequence, then at most one of $L$ and $N$ lies
in $\Sigma'$. Moreover, if an indecomposable summand of $M$
lies in $\Sigma'$, then either $L$ or $N$ lies in $\Sigma'$. 
Here we may
 assume that $Q$ is a finite quiver with no
oriented cycle and that $T\in\ts{mod}\, kQ$ is a tilting module such that $B=\ts{End}_{kQ}(T)$.
Any module $X\in\ts{mod}\, B$ defines the $\Sigma$-module $\ts{Hom}_{B}(\Sigma,X)$
which, as a functor, assigns the vector space $\ts{Hom}_B(E,X)$ to
the object $E$ of $\Sigma$.
By the above properties of $B$, the map $x\mapsto \ts{Hom}_{kQ}\left(T,D(kQe_x)\right)$
defines an isomorphism of $k$-categories $kQ\xrightarrow{\sim} \Sigma$. We
denote by $\Gamma_{\leqslant \Sigma}$ the full subquiver of $\Gamma$ generated by
all the predecessors of $\Sigma$ in $\Gamma$.

\paragraph{The covering of the left section $\Sigma$}$\ $

Let $\widetilde{\Sigma}$ be the full subcategory of $k(\widetilde{\Gamma})$ whose objects are
the $x\in k(\widetilde{\Gamma})$ such that $p(x)\in\Sigma$. Therefore
$p\colon k(\widetilde{\Gamma})\to \ts{ind}\,\Gamma$ induces a covering functor
$p\colon\widetilde{\Sigma}\to \Sigma$. Note that $\widetilde\Sigma$
and $\widetilde{\Gamma}_{\leqslant\widetilde{\Sigma}}$ are
stable under $G$, as subquivers of $\widetilde{\Gamma}$.
Since $\Sigma$ is hereditary, so is $\widetilde{\Sigma}$. Therefore we have
$\widetilde{\Sigma}=k\widetilde Q$ for some quiver $\widetilde Q$. In particular, the
isomorphism $kQ\xrightarrow{\sim} \Sigma$ and the covering functor
$p\colon \widetilde{\Sigma}\to \Sigma$ induce a covering functor $q\colon k\widetilde Q\to
kQ$.

\paragraph{The covering functor of $B$}$\ $

Since $\pi$ and $p$ coincide on vertices, $\pi$ induces a Galois
covering of quivers $\pi\colon \widetilde Q\to Q$ with group
$G$. We write $\pi\colon k\widetilde Q\to kQ$ for the induced
Galois covering with group $G$. Note that $\widetilde Q$ is a disjoint union of copies of the
universal cover of $Q$ because $\widetilde{\Gamma}$ is simply
connected. Also, thanks to the Galois covering $\pi\colon\widetilde Q\to Q$
there is an action of $G$ on $\ts{mod}\, k\widetilde Q$.
Let $T=T_1\oplus\cdots\oplus T_n$ be such that $T_1,\ldots,T_n$ are indecomposable
and $\widetilde B$ be the full subcategory of $\ts{mod}\, k\widetilde Q$ with objects
the $^g\widetilde T_i$ (with $i\in\{1,\ldots,n\}$,
$g\in G$, see \ref{lem:lift}).
\begin{lem}
\label{lem2}
The $k$-category $\widetilde B$ is locally bounded. The push-down functor
$q_{\lambda}\colon \ts{mod}\, k\widetilde Q\to \ts{mod}\, kQ$ induces a covering functor:
\begin{equation}
  \begin{array}{crclc}
    F\colon & \widetilde B & \to & B&\\
    & ^g\widetilde T_i & \mapsto & T_i=q_{\lambda}(\,^g\widetilde T_i)&.
  \end{array}\notag
\end{equation}
Moreover, if $p\colon k(\widetilde{\Gamma})\to\ts{ind}\,\Gamma$ is a Galois covering with group
$\pi_1(\Gamma)$, then so is $F$.
\end{lem}
\noindent{\textbf{Proof:}} We apply the results of the preceding
section to the covering functor $q\colon k\widetilde Q\to kQ$ and the Galois
covering $\pi\colon k\widetilde Q\to kQ$. The first assertion follows from
\ref{lem1} and \ref{prop:cov}, and the second from 
\ref{prop:cov}. The last assertion was proved in
\cite[Lem. 2.2]{lemeur1}.\sq

We also have a Galois covering $\widetilde B\to B$ induced
by the push-down $\pi_{\lambda}\colon \ts{mod}\, k\widetilde Q\to\ts{mod}\, kQ$ (see
\cite[Lem. 2.2]{lemeur1}). In particular, the covering functor
$F\colon\widetilde B\to B$ and the Galois covering $\widetilde B\to B$ coincide on
objects. Therefore we may apply the results of the preceding section
to $F$.

In the sequel, we write $\widetilde T$ for the $k\widetilde Q$-module
$\bigoplus\{\,^g\widetilde T_i\ |\ i\in\{1,\ldots,n\}, g\in G\}$.
Although $\widetilde T$ is not necessarily finite dimensional, it
follows from \ref{prop:cov} that it induces a well-defined functor:
\begin{equation}
\ts{Hom}_{k\widetilde Q}(\widetilde T,-)\colon \ts{mod}\, k\widetilde Q\to \ts{mod}\, \widetilde B\ .\notag
\end{equation}
 More precisely,
if $X\in\ts{mod}\, k\widetilde Q$, then $\ts{Hom}_{k\widetilde Q}(\widetilde T,X)$ is the $\widetilde B$-module
defined by $^g\widetilde T_i\mapsto \ts{Hom}_{k\widetilde Q}(\,^g\widetilde T_i,X)$.
In particular, an object $x$ in $\widetilde{\Sigma}=k\widetilde Q$ defines the injective
$k\widetilde Q$-module $D(k\widetilde Q(x,-))$ which gives rise to the $\widetilde B$-module
$\ts{Hom}_{k\widetilde Q}(\widetilde T,D(k\widetilde Q(x,-)))$. Therefore every $\widetilde B$-module
$X$ defines a $\widetilde{\Sigma}$-module:
\begin{equation}
  \begin{array}{rclc}
    \widetilde{\Sigma}^{op} & \to & \ts{mod}\, k\\
x & \mapsto & \ts{Hom}_{\widetilde B}(\ts{Hom}_{k\widetilde Q}(\widetilde T,D(k\widetilde Q(x,-))),X)&.
  \end{array}\notag
\end{equation}
For reasons that will become clear later, this module is denoted by
$\ts{Hom}_{\widetilde B}(\widetilde{\Sigma},X)$.
In this way, we get a functor $\ts{Hom}_{\widetilde B}(\widetilde{\Sigma},-)\colon \ts{mod}\,\widetilde B\to \ts{mod}\,\widetilde{\Sigma}$.
We need the following result for later reference.
\begin{lem}
\label{lem3}
  The following diagram commutes up to isomorphism of functors:
  \begin{equation}
    \xymatrix{
\ts{mod}\, k\widetilde Q \ar[rrr]^{\ts{Hom}_{k\widetilde Q}(\widetilde T,-)} \ar[d]_{q_{\lambda}} &&& \ts{mod}\,
\widetilde B \ar[d]^{F_{\lambda}} \ar[rrr]^{\ts{Hom}_{\widetilde B}(\widetilde{\Sigma},-)}
&&& \ts{mod}\, \widetilde{\Sigma} \ar[d]^{p_{\lambda}}\\
\ts{mod}\, kQ \ar[rrr]_{\ts{Hom}_{kQ}(T,-)} &&& \ts{mod}\, B \ar[rrr]_{\ts{Hom}_{B}(\Sigma,-)}
&&& \ts{mod}\, \Sigma&.}\notag
  \end{equation}
 Moreover:
\begin{enumerate}
\item[(a)] The two top horizontal arrows are $G$-equivariant.
\item[(b)] If $\theta\colon \ts{mod}\, kQ\to \ts{mod}\, \Sigma$ (or $\widetilde{\theta}\colon \ts{mod}\,
  k\widetilde Q\to \ts{mod}\,\widetilde{\Sigma}$) denotes the composition of the two bottom (or
  top) horizontal arrows, then it induces an equivalence from the
  full subcategory of injective $kQ$-modules (or injective $k\widetilde
  Q$-modules) to the full subcategory of projective $\Sigma$-modules (or
  projective $\widetilde{\Sigma}$-modules, respectively).
\item[(c)] Let $\widetilde{\alpha}\colon\widetilde I\to \widetilde J$ be a surjective morphism between
  injective $k\widetilde Q$-modules. Let
  $\alpha\colon I\to J$ be equal to $q_{\lambda}(\widetilde{\alpha})$.
 Then $F_{\lambda}$ maps the connecting morphism $\ts{Hom}_{k\widetilde
  Q}(\widetilde T,\widetilde J)\to \ts{Ext}_{k\widetilde
  Q}^1(\widetilde T,\ts{Ker} \,\widetilde{\alpha})$ to the connecting morphism $\ts{Hom}_{kQ}(T,J)\to
\ts{Ext}_{kQ}^1(T,\ts{Ker} \,\alpha)$.
\end{enumerate}
\end{lem}
\noindent\textbf{Proof:}
The commutativity of the diagram is an easy exercise on covering functors, and left
to the reader.

(a) This  follows from a direct computation.

(b)  By tilting theory,  $\theta$ induces an
equivalence (see \cite[Chap. VIII Thm. 3.5]{ass}): 
\begin{equation}
  \begin{array}{crclc}
    \Phi\colon  &\ts{inj}\, kQ & \to & \ts{proj}\,\Sigma & \\
    &I & \mapsto & \ts{Hom}_B(\Sigma,\ts{Hom}_{kQ}(T,I)) & .
  \end{array}\notag
\end{equation}
Let $I\in\ts{inj}\, k\widetilde Q$. Then $p_{\lambda}\widetilde{\theta}(I)=\theta
q_{\lambda}(I)$. Moreover, $q_{\lambda}$ maps indecomposable injective
$k\widetilde Q$-modules to indecomposable injective $k Q$-modules,
because so does $\pi_{\lambda}\colon \ts{mod}\,k\widetilde{Q}\to \ts{mod}\,kQ$ (see
\ref{prop:cov}). So $p_{\lambda}\widetilde{\theta}(I)$ is
indecomposable projective, and therefore so is $\widetilde{\theta}(I)$ (see
\cite[3.2]{bog}). Consequently, $\widetilde{\theta}$ induces the following functor:
\begin{equation}
  \begin{array}{crclc}
    \Psi\colon  &\ts{inj}\, k\widetilde Q & \to & \ts{proj}\,\widetilde{\Sigma} & \\
    &I & \mapsto & \ts{Hom}_{\widetilde B}(\widetilde{\Sigma},\ts{Hom}_{k\widetilde Q}(\widetilde T,I)) & .
  \end{array}\notag  
\end{equation}
So we have a commutative diagram:
\begin{equation}
  \xymatrix{
\ts{inj}\, k\widetilde Q \ar[r]^{\Psi} \ar[d]_{q_{\lambda}} & \ts{proj}\,\widetilde{\Sigma}
\ar[d]^{p_{\lambda}}\\
\ts{inj}\, kQ \ar[r]_{\Phi}^{\sim} & \ts{proj}\,\Sigma &.
}\notag
\end{equation}
 In this diagram, $p_{\lambda}$, $q_{\lambda}$ and $\Phi$ are
faithful. Hence, so is $\Psi$.
Let $I,J\in\ts{inj}\, k\widetilde Q$ and $f\colon \Psi(I)\to
  \Psi(J)$. Let $h\colon q_{\lambda}I\to q_{\lambda} J$ be
  such that $\Phi(h)=p_{\lambda}(f)$. Using \ref{prop:cov},
  we have $h=\sum\limits_{g\in G}q_{\lambda}(h_g)$, where
  $(h_g)_g\in \bigoplus\limits_{g\in G} \ts{Hom}_{k\widetilde
    Q}(\widetilde I,\,^g\widetilde J)$. So $p_{\lambda}(f)=\sum\limits_{g\in
    G}p_{\lambda}\Psi(h_g)$. Using \ref{prop:cov} again, we
  deduce that $f=\Psi(h_1)$. So $\Psi$ is full.
Finally, we know from the preceding section that
  $q_{\lambda}\colon\ts{inj}\, k\widetilde Q\to\ts{inj}\, kQ$ is dense. Also, so is
  $p_{\lambda}\colon \ts{proj}\,\widetilde{\Sigma}\to\ts{proj}\,\Sigma$ (see \cite[3.2]{bog}, for instance). Since $\Phi$ is an
  equivalence, we deduce that $\Psi$ is dense.
Therefore $\Psi$ is an equivalence.

(c) The push-down functors $q_{\lambda}$ and $F_{\lambda}$ are exact. So we have a commutative
diagram up to isomorphism of functors:
\begin{equation}
  \xymatrix{
\mathcal D^b(\ts{mod}\, k\widetilde Q) \ar[d]_{q_{\lambda}} \ar[rr]^{\textbf{R}\ts{Hom}_{k\widetilde Q}(\widetilde T,-)}
&& \mathcal D^b(\ts{mod}\, \widetilde B) \ar[d]^{F_{\lambda}} \\
\mathcal D^b(\ts{mod}\, kQ) \ar[rr]_{\textbf{R}\ts{Hom}_{kQ}(T,-)} && \mathcal D^b(\ts{mod}\, B) &.
}\notag
\end{equation}
The statement follows from this diagram.
\sq

We wish to construct a functor $\varphi\colon k(\widetilde{\Gamma})\to \ts{mod}\,\widetilde
B$. We proceed in several steps:
\begin{enumerate}
\item Define a functor $\varphi_0\colon k(\widetilde{\Gamma}_{\leqslant \widetilde{\Sigma}})\to
  \ts{mod}\,\widetilde B$ where $k(\widetilde{\Gamma}_{\leqslant\widetilde{\Sigma}})$ denotes the full subcategory
  of $k(\widetilde{\Gamma})$ with objects the vertices in $\widetilde{\Gamma}_{\leqslant\widetilde{\Sigma}}$.
\item Define $\varphi$ on objects, so that it coincides with
  $\varphi_0$ on predecessors of $\widetilde{\Sigma}$.
\item Define $\varphi$ on morphisms, so that it extends $\varphi_0$.
\end{enumerate}

\paragraph{The functor $\varphi_0\colon k(\widetilde{\Gamma}_{\leqslant \widetilde{\Sigma}})\to
  \ts{mod}\,\widetilde B$}\label{ss1.2}$\ $

We first prove the following
lemma. In the case of a Galois covering whose group acts
freely on indecomposables, a corresponding
result was proved in \cite[3.6]{gabriel}. We know that
$p\colon\widetilde{\Sigma}\to\Sigma$ and $F\colon\widetilde B\to B$ are covering functors, and
that the latter coincides on objects with a Galois covering $\widetilde B\to B$
with group $G$. Finally, if $X\in\ts{ind}\, B$ is a summand
of a tilting $B$-module, then $\widetilde X\in\ts{ind}\,\widetilde B$ is such that
$p_{\lambda}(\widetilde X)\simeq X$ (see \ref{prop:cov}). 
\begin{lem}
  \label{lem:tau}
Let $X\in\Gamma_{\leqslant \Sigma}$ and $g_0\in G$. If $\widetilde u\colon \widetilde E\to
\,^{g_0}\widetilde X$ is right minimal almost split, then so is $p_{\lambda}\widetilde
u\colon p_{\lambda}\widetilde E\to X$. Consequently, $p_{\lambda}\tau_{\widetilde B}(\,^{g_0}\widetilde
X)\simeq \tau_BX$ if $X$ is not projective.
\end{lem}
\noindent{\textbf{Proof:}} Notice that $X$ is an indecomposable
summand of some tilting $B$-module. So we may apply the results of
\ref{prop:cov}. If $X$ is projective, the assertion follows
from \cite[3.2]{bog}. So we assume that $X$, and therefore $^{g_0}\widetilde X$, are not
projective. Let $u\colon E\to X$ be right minimal almost split and
$E=E_1\oplus\cdots\oplus E_t$ be such that $E_1,\ldots,E_t$ are indecomposable.
 Since $\Gamma_{\leqslant\Sigma}$ is acyclic (see 
\cite[2.2]{assem}), we have
$\ts{Ext}^1_B(E,\tau_BX)=0$. Also, the linear map
$\bigoplus\limits_{g\in G}\ts{Hom}_{\widetilde B}(\,^g\widetilde E_i,\,^{g_0}\widetilde
X)\to \ts{Hom}_B(E_i,X)$ is bijective, for every $i$ (see
\ref{prop:cov}). Therefore we apply
\ref{lem:lift} to the exact sequence
$0\to \tau_BX\to E\xrightarrow{u}X\to 0$: There exist
$g_1,\ldots,g_t\in G$ and  morphisms $\widetilde
u_i\colon\,^{g_i}\widetilde E_i\to \,^{g_0}\widetilde X$ ($i\in\{1,\ldots,t\}$)
fitting into a commutative diagram whose  vertical arrow on the left
is an isomorphism:
\begin{equation}
  \xymatrix{
E \ar[rr]^{p_{\lambda}[\widetilde u_1,\ldots,\widetilde u_t]} \ar[d]_{\sim} && X
  \ar@{=}[d]\\
E \ar[rr]_{u} && X&.
}\notag
\end{equation}
We identify $u$ and $p_{\lambda}[\widetilde u_1,\ldots,\widetilde u_t]$ via this
diagram. Let $i\in\{1,\ldots,n\}$. Then $\widetilde u_i\colon\,^{g_i}\widetilde E_i\to
\,^{g_0}\widetilde X$ is not a retraction because $^{g_i}\widetilde E_i$ and $^{g_0}\widetilde
X$ are non-isomorphic indecomposable modules. So $\widetilde u_i$ factors
through $\widetilde u$, for every $i$. Applying $p_{\lambda}$ to each
factorisation shows that $u$ factors through $p_{\lambda}(\widetilde u)$. On
the other hand, $p_{\lambda}\widetilde u$ is not a retraction because $X$ is
not a direct summand of $E$. So $p_{\lambda}(\widetilde u)$ factors through
$u$. The right minimality of $u$ implies that the morphism $u$ is a
direct summand of $p_{\lambda}(\widetilde u)$. Finally, the following
equality follows from \ref{lem:dimtau}:
\begin{equation}
  \dim \ts{Ker}  u=\dim \tau_BX=\dim\tau_Bp_{\lambda}\,^{g_0}\widetilde X=\dim
  p_{\lambda}\tau_{\widetilde B}\,^{g_0}\widetilde X=\dim\ts{Ker} 
  p_{\lambda}(\widetilde u)\ .\notag
\end{equation}
So $p_{\lambda}(\widetilde u)$ and $u$ are isomorphic, and $p_{\lambda}(\widetilde u)$ is
 right minimal almost split.\sq

Using the preceding lemma, we construct a functor
$\varphi_0\colon k(\widetilde{\Gamma}_{\leqslant \widetilde{\Sigma}})\to \ts{ind}\,\widetilde B$.
\begin{lem}
  \label{lem8}
There exists a full and faithful functor, $G$-equivariant on vertices,
$ \varphi_0\colon k(\widetilde{\Gamma}_{\leqslant \widetilde{\Sigma}})\to  \ts{ind}\, \widetilde B$.
This functor maps arrows in $\widetilde{\Gamma}_{\leqslant \widetilde{\Sigma}}$ to irreducible
maps, and meshes to almost split sequences. Moreover, it commutes with the translations and  extends the
canonical functor $\widetilde{\Sigma} \to \ts{ind}\, \widetilde B$
    defined on the objects by $x\mapsto  \ts{Hom}_{k\widetilde Q}(\widetilde T,D(k\widetilde Q(x,-))$.
Finally, the following diagram is commutative up to isomorphism of
functors:
\begin{equation}
  \xymatrix{
k(\widetilde{\Gamma}_{\leqslant \widetilde{\Sigma}})
\ar@{->}[r]^{\varphi_0} \ar@{->}[d]_p & \ts{ind}\, \widetilde
B \ar@{^(->}[r] & \ts{mod}\,\widetilde{B} \ar@{->}[d]^{F_{\lambda}}\\
\ts{ind}\,(\Gamma_{\leqslant \Sigma}) \ar@{^(->}[r] & \ts{ind}\, B
\ar@{^(->}[r] & \ts{mod}\,B.
}\notag
\end{equation}
\end{lem}
\noindent{\textbf{Proof:}} 
\textbf{Step $1$:} Clearly there is a functor $\varphi_0\colon
\widetilde{\Sigma}\to\ts{mod}\,\widetilde B$ given by $\widetilde
x\mapsto \ts{Hom}_{k\widetilde Q}(\widetilde T,D(k\widetilde
Q(\widetilde x,-)))$.
Note that $\widetilde{\Sigma}$ (or $\Sigma$) is naturally
equivalent to the full subcategory of $\ts{mod}\, k\widetilde Q$ (or $\ts{mod}\, kQ$,
respectively) consisting  of the indecomposable
injective modules. Therefore \ref{lem3} shows that this functor
is full and faithful, and that the following diagram commutes up to
isomorphism:
\begin{equation}
  \xymatrix{
\widetilde{\Sigma} \ar[r]^{\varphi_0} \ar[d]_{p} &
\ts{ind}\,\widetilde B \ar@{^(->}[r] &\ts{mod}\,\widetilde{B} \ar[d]^{F_{\lambda}}\\
\Sigma \ar@{^(->}[r] & \ts{ind}\, B \ar@{^(->}[r] & \ts{mod}\,B.
}\notag
\end{equation}
Note that $\varphi_0(M)$ is indecomposable for every $M$ because so is $F_{\lambda}\varphi_0(M)=p(M)$.
The functor $\varphi_0\colon\widetilde{\Sigma}\to\ts{ind}\,\widetilde B$ is
$G$-equivariant on vertices: Indeed, for every
$g\in G$, and every $\widetilde x\in\widetilde Q_0$, we have:
\begin{equation}
\varphi_0(g\widetilde
  x)=\ts{Hom}_{k\widetilde Q}(\widetilde T,D(k\widetilde Q({g\widetilde x},-)))=\ts{Hom}_{k\widetilde
    Q}(\widetilde T,\,^gD(k\widetilde Q({\widetilde
    x},-)))=\,^g\ts{Hom}_{k\widetilde Q}(\widetilde T,D(k\widetilde
  Q({\widetilde
    x},-)))=\,^g\varphi_0(\widetilde x)\ \ .\notag
\end{equation}

\textbf{Step $2$:}  If $M\in
k(\widetilde{\Gamma}_{\leqslant \widetilde{\Sigma}})$, there exists a unique $n\in\mathbb{N}$ such
that $\tau^{-n}M\in\widetilde{\Sigma}$. Let $\varphi_0(M)$ be the
$\widetilde B$-module:
\begin{equation}
  \varphi_0(M)=\tau_{\widetilde B}^n\varphi_0(\tau^{-n}M)\ \ .\notag
\end{equation}
It follows from \ref{lem:tau} that $F_{\lambda}\varphi_0(M)=p(M)$. Also
$\varphi_0(\,^gM)=\,^g\varphi_0(M)$ for every $g\in G$ and for every
vertex $M$
because $\tau$ commutes with the action of $G$. 

\textbf{Step $3$:} In order to define $\varphi_0$ on morphisms, we
construct inductively a sequence of $G$-invariant
left sections $\widetilde{\Sigma}_i$ of $\widetilde{\Gamma}$ such that
 $\widetilde{\Sigma}_0=\widetilde{\Sigma}$,
such that $\widetilde{\Sigma}_{i+1}\backslash\widetilde{\Sigma}_i$ consists of the $G$-orbit of a
  vertex, and such that, if 
$\bigcup\limits_{t=1}^i\widetilde{\Sigma}_t$ denotes the full subcategory of the
path category $k\widetilde{\Gamma}$
whose vertices are given by those of $\widetilde{\Sigma}_0,\ldots,\widetilde{\Sigma}_i$, then
$k\widetilde{\Gamma}_{\leqslant \widetilde{\Sigma}}=\bigcup\limits_{i\geqslant 0}\widetilde{\Sigma}_i$.
Each inductive step defines a functor $\varphi_0\colon
\bigcup\limits_{t=1}^i\widetilde{\Sigma}_t\to \ts{ind}\,\widetilde B$ which maps arrows to
irreducible maps and extends the
construction of the two preceding steps. 
This functor makes the following diagram commute:
  \begin{equation}
    \xymatrix{\bigcup\limits_{t=1}^i\widetilde{\Sigma}_t \ar@{->}[d]
      \ar@{->}[r]^{\varphi_0} & \ts{ind}\, \widetilde B \ar@{^(->}[r]
      & \ts{mod}\,\widetilde B\ar@{->}[d]^{F_{\lambda}}\\
\ts{ind}\, (\Gamma_{\leqslant \Sigma}) \ar@{^(->}[r] & \ts{ind}\, B
\ar@{^(->}[r]& \ts{mod}\,B,}\notag
  \end{equation}
where the vertical arrow on the left
is induced by $p$.
Assume that $\varphi_0\colon\bigcup\limits_{t=1}^i\widetilde{\Sigma}_t\to\ts{ind}\,\widetilde B$
has been defined for some $i\geqslant 0$. Since
$\widetilde{\Sigma}_i$ is acyclic, it has a sink. Assume that all sinks are
projective. First assume that $P$ is a projective sink, and let
$\widetilde{\Sigma}_{i+1}$ be equal to $\widetilde{\Sigma}_i\backslash\{\,^gP\ |\ g\in G\}$; then
$\widetilde{\Sigma}_{i+1}$ is a left section of $\Gamma$, and there is a unique
$\varphi_0\colon\bigcup\limits_{t=0}^{i+1}\widetilde{\Sigma}_t\to\ts{ind}\,\widetilde B$
satisfying the required conditions. Now assume that there is a non
projective sink $M$ in
$\widetilde{\Sigma}_i$. Then there exists a mesh in $\widetilde{\Gamma}$:
\begin{equation}
\xymatrix@R=2pt{  
&N_1\ar[rd]^{u_1} &\\
\tau M \ar[ru]^{v_1} \ar[rd]_{v_s} & \vdots & M & ,\\
&N_s \ar[ru]_{u_s} &
}\notag
\end{equation}
and $M,N_1,\ldots,N_s\in\widetilde{\Sigma}_i$ because $M\in\widetilde{\Sigma}_i$ is a sink. In
particular, $\varphi_0(u_j)$ is defined, and 
$F_{\lambda}\,\varphi_0(u_j)=p(u_j)_{|B}$ for every $i$. For simplicity, we
write $\varphi_0(u)=\begin{bmatrix}\varphi_0(u_1) & \ldots & \varphi_0(u_s)\end{bmatrix}$
and $p(u)=\begin{bmatrix}p(u_1) & \ldots&p(u_s)\end{bmatrix}$. Then
  $\varphi_0(u)$  is right minimal almost split in $\ts{mod}\, \widetilde 
B$: Indeed, there exists a right minimal almost split morphism
  $L\xrightarrow{w} \varphi_0(M)$.
Since $\varphi_0(u)\colon \bigoplus\limits_j \varphi_0(N_j)\to
\varphi_0(M)$ is not a retraction (because each $\varphi_0(u_j)$ is an
irreducible morphism, by the induction hypothesis), there
exists a morphism $w'\colon\bigoplus\limits_{j} \varphi_0(N_j)\to L$
such that $\varphi_0(u)=ww'$; applying $F_{\lambda}$, we have
$F_{\lambda}\varphi_0(u)=F_{\lambda}(w)F_{\lambda}(w')$; but now
$F_{\lambda}\varphi_0(u)=p(u)_{|B}$ is right minimal almost split by
construction, and so is $F_{\lambda}(w)$ (see \ref{lem:tau}); hence,
$F_{\lambda}(w')$ is an isomorphism and
therefore so is $w'$ because $F_{\lambda}$ is exact. We let
$\begin{bmatrix}\varphi_0(v_1) & \ldots & \varphi_0(v_s)\end{bmatrix}^t\colon\varphi_0(\tau
M)\to
\bigoplus\limits_{j=1}^s\varphi_0(N_j)$ be the kernel of
$\varphi_0(u)$. For simplicity, we set
$\varphi_0(v)=\begin{bmatrix}\varphi_0(v_1)&\ldots &
  \varphi_0(v_s)\end{bmatrix}^t$
and $p(v)=\begin{bmatrix}p(v_1)& \ldots & p(v_s)\end{bmatrix}^t$. We let
$\widetilde{\Sigma}_{i+1}=\left(\widetilde{\Sigma}_i\backslash\{\,^gM\ |\ g\in
  G\}\right)\bigcup \{\,^g\tau M\ |\ g\in G\}$. Clearly,
$\widetilde{\Sigma}_{i+1}$ is a left section. We now show that we may assume
$\varphi_0(v)$ to be taken such that
$F_{\lambda}\varphi_0(v)=p(u)_{|B}$. Indeed, the commutative diagram
with exact rows:
\begin{equation}
  \xymatrix{
0 \ar@{->}[r]  & F_{\lambda}\varphi_0(\tau M)
\ar@{->}[rr]^{F_{\lambda}\varphi_0(v)} && F_{\lambda}\varphi_0\left(\bigoplus\limits_{j=1}^sN_j\right) \ar@{=}[d]
\ar@{->}[rr]^{F_{\lambda}\varphi_0(u)} && F_{\lambda}\varphi_0(M) \ar@{->}[r] \ar@{=}[d]&
0 & \\
0 \ar@{->}[r] & p(\tau M)_{|B} \ar@{->}[rr]_{p(v)_{|B}} &&
p\left(\bigoplus\limits_{j=1}^sN_j\right)_{|B}
\ar@{->}[rr]_{p(u)_{|B}} && p(M)_{|B} \ar@{->}[r] & 0 &
}\notag
\end{equation}
gives an isomorphism $F_{\lambda}\varphi_0(\tau M)\to
p(\tau M)_{|B}$ making the left square commute. Since
$F_{\lambda}\varphi_0(\tau M)=p(\tau M)_{|B}$ is a brick (because it
belongs to $\Gamma_{\leqslant \Sigma}$), this isomorphism is the
multiplication by a non-zero constant $c$. Hence, $p(v)_{|B}=c\,
F_{\lambda}\varphi_0(v)$. Replacing $\varphi_0(v)$ by $c\,\varphi_0(v)$
does the trick. 
 Thus, we have defined
 $\varphi_0\colon\bigcup\limits_{t=1}^{i+1}\widetilde{\Sigma}_t\to \ts{ind}\,\widetilde B$.
  Clearly, the required conditions are satisfied. This induction
 gives a functor $\varphi_0\colon k\widetilde{\Gamma}_{\leqslant
   \widetilde{\Sigma}}\to\ts{ind}\,\widetilde B$ mapping arrows to irreducible maps and meshes
 to almost split sequences, and such that the following diagram commutes:
 \begin{equation}
   \xymatrix{
k\widetilde{\Gamma}_{\leqslant\widetilde{\Sigma}} \ar[d] \ar[r]^{\varphi_0} & \ts{ind}\,\widetilde B
\ar@{^(->}[r] &\ts{mod}\,\widetilde{B}\ar[d]^{F_{\lambda}}\\
\ts{ind}\,\Gamma_{\leqslant \Sigma} \ar@{^(->}[r] & \ts{ind}\, B \ar@{^(->}[r]& \ts{mod}\,B,
}\notag
 \end{equation}
where the vertical arrow on the left is induced by $p$. Since
$F_{\lambda}$ is faithful, $\varphi_0$ induces a functor $\varphi_0\colon
k(\widetilde{\Gamma}_{\leqslant\widetilde{\Sigma}})\to\ts{ind}\,\widetilde B$. It is 
now clear that this functor satisfies the conditions of the lemma.
\sq

It was shown in \cite[3.2]{assem} that the existence of a left section $\Sigma$
in an Auslander-Reiten component $\Gamma$ implies that $\Gamma_{\leqslant \Sigma}$
is generalised standard. We now prove that it is standard.
\begin{cor}
  \label{cor:standard-left-section}
Let $A$ be a finite dimensional $k$-algebra and $\Gamma$ be a  component
of $\Gamma(\ts{mod}\, A)$ having  a left section $\Sigma$. Then $\Gamma_{\leqslant \Sigma}$ is standard.
\end{cor}
\noindent{\textbf{Proof:}} Let $B=A/\ts{Ann}\,\Sigma$. Then $B$ is a product of
tilted algebras and the components of $\Sigma$ form complete slices of the
connecting components of the connected
 components of $B$. Let $\Gamma'$ be the union of the components of
$\Gamma(\ts{mod}\, B)$ intersecting $\Sigma$. The
arguments of the
proof of \ref{lem8} show that there exists  a full and faithful functor
$k(\Gamma'_{\leqslant \Sigma})\to\ts{ind}\, \Gamma'_{\leqslant \Sigma}$ extending the identity on
vertices. So $\Gamma_{\leqslant \Sigma}=\Gamma_{\leqslant \Sigma}'$ is standard.\sq

\begin{ex}
Let $A$ be the algebra given
by the quiver
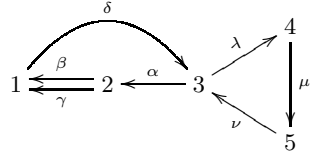
\begin{equation}
\xymatrix@R=10pt{
   &   &   & 4 \ar[dd]^{\mu}\\
1 \ar@/^2pc/[rr]^{\delta}& 2 \ar@<2pt>[l]^{\gamma} \ar@<-2pt>[l]_{\beta}
& 3 \ar[l]_{\alpha} \ar[ru]^{\lambda}& \\
  &   &   & 5 \ar[lu]^{\nu}
}\notag
\end{equation}
and the potential $W=\delta\beta\alpha+\nu\mu\lambda$ (or,
equivalently, by the relations $\beta\alpha=0$, $\delta\beta=0$,
$\alpha\delta=0$, $\mu\lambda=0$, $\nu\mu=0$ and $\lambda\nu=0$).
Then $A$ is a  
cluster-tilted algebra since it is the relation-extension (in the
sense of \cite{ABS08a}) of the tilted algebra of type 
 $\widetilde{\mathbb{A}}$ given by the quiver
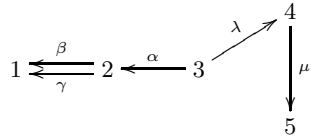
\begin{equation}
\xymatrix@R=10pt{
   &   &   & 4 \ar[dd]^{\mu}\\
1 & 2 \ar@<2pt>[l]^{\gamma} \ar@<-2pt>[l]_{\beta}
& 3 \ar[l]_{\alpha} \ar[ru]^{\lambda}& \\
  &   &   & 5 
}\notag
\end{equation}
bound by $\beta\alpha=0$ and $\mu\lambda=0$.
 The transjective component $\Gamma$ of $\Gamma(\ts{mod}\, A)$ is of the form
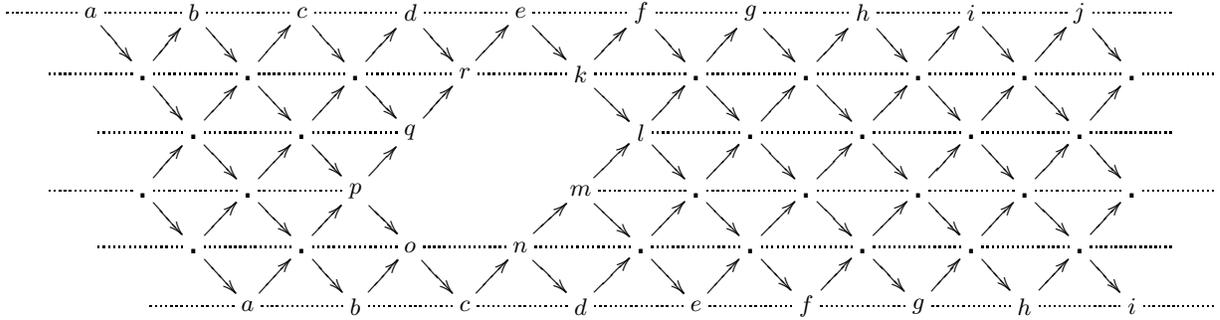
\begin{equation}
   \xymatrix@=10pt{
\ar@{.}[rr]&& a \ar@{.}[rr] \ar@{->}[rd]&&  b \ar@{.}[rr]
\ar@{->}[rd]&&  c 
\ar@{.}[rr] \ar@{->}[rd]&&  d \ar@{.}[rr] \ar@{->}[rd]&&
 e \ar@{.}[rr] \ar@{->}[rd]&& 
f \ar@{.}[rr] \ar@{->}[rd]&&  g \ar@{.}[rr]
\ar@{->}[rd]&&  h 
\ar@{.}[rr] \ar@{->}[rd]&&  i \ar@{.}[rr] \ar@{->}[rd]&&
j \ar@{->}[rd] \ar@{.}[rr] &&
\\
&\ar@{.}[rr] &&\centerdot \ar@{.}[rr]  \ar@{->}[rd] \ar@{->}[ru]&&  \centerdot \ar@{.}[rr]
\ar@{->}[rd]\ar@{->}[ru] &&  \centerdot 
\ar@{.}[rr] \ar@{->}[rd]\ar@{->}[ru] &&  r \ar@{.}[rr] \ar@{->}[ru] &&
k \ar@{.}[rr] \ar@{->}[rd]\ar@{->}[ru] && 
\centerdot \ar@{.}[rr] \ar@{->}[rd]\ar@{->}[ru] &&  \centerdot \ar@{.}[rr]
\ar@{->}[rd]\ar@{->}[ru] &&  \centerdot 
\ar@{.}[rr] \ar@{->}[rd]\ar@{->}[ru] &&  \centerdot \ar@{.}[rr] \ar@{->}[rd]\ar@{->}[ru] &&
\centerdot \ar@{.}[rr] &&
\\
 && \ar@{.}[rr] &&  \centerdot
\ar@{.}[rr] \ar@{->}[rd] \ar@{->}[ru] &&  \centerdot 
\ar@{.}[rr] \ar@{->}[rd] \ar@{->}[ru] &&  q
 \ar@{->}[ru] &&   && 
l \ar@{.}[rr] \ar@{->}[rd] \ar@{->}[ru] &&  \centerdot
\ar@{.}[rr] \ar@{->}[rd] \ar@{->}[ru] &&  \centerdot 
\ar@{.}[rr] \ar@{->}[rd] \ar@{->}[ru] &&  \centerdot \ar@{.}[rr]
\ar@{->}[rd] \ar@{->}[ru] && 
\centerdot \ar@{->}[rd] \ar@{->}[ru] \ar@{.}[rr] &&
\\
&\ar@{.}[rr] &&\centerdot \ar@{.}[rr] \ar@{->}[rd] \ar@{->}[ru] &&  \centerdot
\ar@{.}[rr] \ar@{->}[rd] \ar@{->}[ru] &&  p 
 \ar@{->}[rd] \ar@{->}[ru] &&    &&  m \ar@{.}[rr] \ar@{->}[ru]
 \ar@{->}[rd] &&  
\centerdot \ar@{.}[rr] \ar@{->}[rd] \ar@{->}[ru] &&  \centerdot
\ar@{.}[rr] \ar@{->}[rd] \ar@{->}[ru] &&  \centerdot 
\ar@{.}[rr] \ar@{->}[rd] \ar@{->}[ru] &&  \centerdot \ar@{.}[rr] 
\ar@{->}[rd] \ar@{->}[ru] &&  \centerdot \ar@{.}[rr] &&
\\
&& \ar@{.}[rr]  &&  \centerdot
\ar@{.}[rr] \ar@{->}[rd] \ar@{->}[ru] &&  \centerdot 
\ar@{.}[rr] \ar@{->}[rd] \ar@{->}[ru] &&  o \ar@{.}[rr]
\ar@{->}[rd]  &&  n \ar@{.}[rr] \ar@{->}[rd]
\ar@{->}[ru] && 
\centerdot \ar@{.}[rr] \ar@{->}[rd] \ar@{->}[ru] &&  \centerdot
\ar@{.}[rr] \ar@{->}[rd] \ar@{->}[ru] &&  \centerdot 
\ar@{.}[rr] \ar@{->}[rd] \ar@{->}[ru] &&  \centerdot \ar@{.}[rr]
\ar@{->}[rd] \ar@{->}[ru] && 
\centerdot \ar@{->}[rd] \ar@{->}[ru] \ar@{.}[rr] &&
\\
& && \ar@{.}[rr] &&  a \ar@{.}[rr]
\ar@{->}[ru] &&  b 
\ar@{.}[rr] \ar@{->}[ru] &&  c \ar@{.}[rr] \ar@{->}[ru] &&
d \ar@{.}[rr] \ar@{->}[ru] && 
e \ar@{.}[rr] \ar@{->}[ru] &&  f \ar@{.}[rr]
\ar@{->}[ru] &&  g 
\ar@{.}[rr] \ar@{->}[ru] &&  h \ar@{.}[rr] \ar@{->}[ru] &&
i \ar@{.}[rr] &&
}\notag
\end{equation}
where  vertices with the same label are identified. Then $\Gamma$ admits
a left section $\Sigma=\{e,r,q,p,o,c\}$ and $B=A/\ts{Ann}\,\Sigma$ is the
algebra given by the quiver:
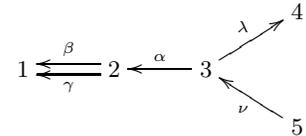
\begin{equation}
\xymatrix@R=10pt{
   &   &   & 4 \\
1 & 2 \ar@<2pt>[l]^{\gamma} \ar@<-2pt>[l]_{\beta}
& 3 \ar[l]_{\alpha} \ar[ru]^{\lambda}& 
\\
  &   &   & 5 \ar[lu]^{\nu}
}\notag  
\end{equation}
with the inherited relations. As we have seen, $\Gamma_{\leqslant \Sigma}$ is
standard (and generalised standard) while  $\Gamma$ itself is
not.
\end{ex}

The following corollary seems to be well-known. However we have been
unable to find a reference.
\begin{cor}
\label{rem:tilted_standard}
Let $B$ be a tilted algebra and $\Gamma$
 a connecting component of $B$. Then $\Gamma$ is standard.
\end{cor}
\noindent{\textbf{Proof:}} If $B$ is concealed, this follows
from \cite[2.4 (11) p. 80]{R84}. Assume that $B$ is not concealed. So $\Gamma$ is the
unique connecting component of $B$. Let $\Sigma$ be a complete slice in $\Gamma$.
As observed in \ref{cor:standard-left-section}, we have a full and
faithful functor $k(\Gamma_{\leqslant\Sigma})\to \ts{ind}\,\Gamma$ extending the
identity on vertices. A dual
construction extends this functor to a full and faithful functor
$k(\Gamma)\to\ts{ind}\,\Gamma$ extending the identity  on vertices. So $\Gamma$ is standard.
\sq

From now on, we identify $\widetilde{\Sigma}$ to a full subcategory of
$\ts{mod}\,\widetilde B$ by means of $\varphi_0$.

\paragraph{Construction of $\varphi$ on objects}$\ $

We prove that for any $M\in\widetilde{\Gamma}$, there exists $\varphi(M)\in \ts{mod}\, \widetilde
B$ whose image under $F_{\lambda}\colon \ts{mod}\, \widetilde B\to \ts{mod}\, B$
coincides with $p(M)_{|B}$, in such a way that
$\varphi(\,^gM)=\,^g\varphi(M)$, for every $g\in G$. 
We define $\mathcal L_{\Sigma}$ to be the full subcategory of $\ts{ind}\, B$ which consists
of the predecessors of the complete slice $\Sigma$. Also a \emph{minimal
  $\ts{add}\,\mathcal L_{\Sigma}$-presentation} of a module $R$ is a sequence of morphisms
$E_1\to E_2\to R$
where the morphism on the right is a minimal
$\ts{add}\,\mathcal L_{\Sigma}$-approximation and the one on the left is a minimal
$\ts{add}\,\mathcal L_{\Sigma}$-approximation of its kernel. Before constructing
$\varphi(M)$, we prove some lemmata.
\begin{lem}
  \label{lem11}
Let $R\in \ts{mod}\, B$ be a module with no direct summand in $\mathcal L_{\Sigma}$. There
exists an exact sequence in $\ts{mod}\, B$, which
is a minimal $\ts{add}\, \mathcal L_{\Sigma}$-presentation:
\begin{equation}
  0\to E_1\to E_2\to R\to 0\tag{$\star$}
\end{equation}
with $E_1,E_2\in \ts{add}\, \Sigma$. Moreover, the functor $\ts{Hom}_{kQ}(T,-)$
induces a bijection between the class of all such exact sequences, and
the class of minimal injective copresentations:
\begin{equation}
0\to \ts{Tor}_1^{B}(R,T)\to I_1\to I_2\to 0\ \  .\notag
\end{equation}
Finally, there is an isomorphism in $\ts{mod}\, B$:
\begin{equation}
  R\simeq \ts{Ext}_{kQ}^1(T,\ts{Tor}_1^{B}(R,T))\ \ .\notag
\end{equation}
\end{lem}
\noindent{\textbf{Proof:}}  
Let $\mathcal X(T)$ be the torsion class induced by $T$ in $\ts{mod}\, B$. So $R$
lies in $\mathcal X(T)$ and has no direct summand in $\Sigma$. Therefore $R$ is the
epimorphic image of a module in $\ts{add}\,\Sigma$. The first assertion then
follows from  \cite[2.2, (d)]{apt}.

 Let $f\colon I_1\to I_2$ be the morphism between injective
$kQ$-modules such that $\ts{Hom}_{kQ}(T,f)$ is equal to the morphism
$E_1\to E_2$ in $(\star)$. Because of the Brenner-Butler Theorem (see
\cite[Chap. VI, Thm. 3.8, p.207]{ass}), the functor $-\underset{B}{\otimes}T$ applied to
$(\star)$ yields an injective copresentation in $\ts{mod}\, kQ$:
\begin{equation}
  0\to \ts{Tor}_1^{B}(R,T)\to I_1\to I_2\to 0\ \ .\notag
\end{equation}
The minimality of this copresentation follows from the minimality of
$E_2\to R$. With these arguments, it is straightforward to check
that there is a well-defined bijection which carries the equivalence
class of the exact sequence $0\to E_1\to E_2\to R\to 0$ to the
equivalence class of the exact
sequence $0\to \ts{Tor}^B_1(R,T)\to I_1\to I_2\to 0$.

The last assertion follows from the Brenner-Butler
Theorem and the fact that $R\in\mathcal X(T)$.\sq

\begin{lem}
\label{lem:ls}
$\ts{add}\,\mathcal L_{\Sigma}$ is contravariantly finite in $\ts{mod}\, A$. Therefore if 
  $X\in\Gamma\backslash\mathcal L_{\Sigma}$, then $X_{|B}$ lies in  the torsion class induced by $T$ in $\ts{mod}\, B$.
\end{lem}
\noindent{\textbf{Proof:}} By \cite[Thm. B]{assem}, the algebra $B$ is the
endomorphism algebra of the indecomposable projective $A$-modules
 in $\mathcal L_{\Sigma}$. In particular, a projective $B$-module is projective
as an $A$-module so the projective dimensions in $\ts{mod}\, A$ and in
$\ts{mod}\, B$ coincide on $\mathcal L_{\Sigma}$. Also, by \cite[Thm. B]{assem}, all modules in $\mathcal L_{\Sigma}$ have
projective dimension at most
one as $B$-modules. Therefore $\mathcal L_{\Sigma}\subseteq\mathcal L_A$. Moreover,
$\bigoplus\Sigma$ is sincere as a $B$-module. Hence, \cite[8.2]{assem} implies that
$\ts{add}\,\mathcal L_{\Sigma}$ is contravariantly finite in $\ts{mod}\, A$.
Let $X\in\Gamma\backslash\mathcal L_{\Sigma}$. Let $P\twoheadrightarrow X_{|B}$ be a projective cover
in $\ts{mod}\, B$. As
noticed above, we have $P\in\ts{add}\,\mathcal L_{\Sigma}$. Therefore $P\twoheadrightarrow X_{|B}$ factors
through $\ts{add}\,\Sigma$. Thus, $X_{|B}$ lies in the torsion class.\sq

\begin{lem}
\label{lem12}
  There exists a map $\varphi\colon \widetilde{\Gamma}_o\to \ts{mod}\, \widetilde
  B$ extending $\varphi_0$, and such that
  $F_{\lambda}(\varphi(M))=p(M)_{|B}$, for every $M\in\widetilde{\Gamma}$. Moreover,
  $\varphi(\,^gM)=\,^g\varphi(M)$ for every $g\in G$ and
  $M\in\widetilde{\Gamma}$.
\end{lem}
\noindent{\textbf{Proof:}} Note that $\varphi$ is already defined on
$\widetilde{\Gamma}_{\leqslant \widetilde{\Sigma}}$ because of \ref{lem8}. Let
$M\in\widetilde{\Gamma}\backslash\widetilde{\Gamma}_{\leqslant \widetilde{\Sigma}}$. Then
$p(M)\in\Gamma\backslash\Gamma_{\leqslant \Sigma}=\Gamma\backslash\mathcal L_{\Sigma}$. By
\ref{lem:ls}, the module $p(M)_{|B}$ lies in the torsion class induced
by $T$ in $\ts{mod}\, B$. 
 So there is a decomposition in $\ts{mod}\, B$:
\begin{equation}
  p(M)_{|B}=R\oplus E\ \ ,\notag
\end{equation}
where $E\in \ts{add}\, \Sigma$ and $R$ has no indecomposable summand in $\mathcal L_{\Sigma}$. Also, fix a
decomposition in $\ts{mod}\, \widetilde{\Sigma}$:
\begin{equation}
  k(\widetilde{\Gamma})(\widetilde{\Sigma},M)=\widetilde R\oplus \widetilde P\ \ ,\notag
\end{equation}
where $\widetilde P$ is projective and maximal for this property. Let $\widetilde
E\in\ts{add}\,\widetilde{\Sigma}$ be such that $\widetilde P=k(\widetilde{\Gamma})(\widetilde{\Sigma},\widetilde E)$.

We claim that $p_{\lambda}\colon\ts{mod}\,\widetilde{\Sigma}\to\ts{mod}\,\Sigma$  maps $\widetilde R$ and
$\widetilde P$ to $\ts{Hom}_B(\Sigma,R)$
and $\ts{Hom}_B(\Sigma,E)$ respectively. Indeed, since  $p\colon k(\widetilde{\Gamma})\to
\ts{ind}\, \Gamma$ is a covering functor inducing $p\colon\widetilde{\Sigma}\to\Sigma$, the
image of  $k(\widetilde{\Gamma})(\widetilde{\Sigma},M)$ under
$p_{\lambda}\colon\ts{mod}\,\widetilde{\Sigma}\to\ts{mod}\,\Sigma$ is
$\ts{Hom}_A(\Sigma,p(M))\cong \ts{Hom}_B(\Sigma,p(M)_{|B})$
(functorially in $M$). Moreover, the decomposition
$p(M)_{|B}=R\oplus E$ in $\ts{mod}\, B$ gives a decomposition
$\ts{Hom}_B(\Sigma,p(M))=\ts{Hom}_B(\Sigma,R)\oplus \ts{Hom}_B(\Sigma,E)$ in
$\ts{mod}\, \Sigma$ where $\ts{Hom}_B(\Sigma,E)$ is projective and
$\ts{Hom}_B(\Sigma,R)$ has no non-zero projective direct summand. The claim then
follows from \cite[3.2]{bog}.

In order to prove that $R$ is the image of a $\widetilde B$-module under
$F_{\lambda}$, we consider a minimal projective presentation in $\ts{mod}\,\widetilde{\Sigma}$:
\begin{equation}
  0\to \widetilde P_1\to\widetilde P_2\to\widetilde R\to 0\ .\notag 
\end{equation}
Then there exists a morphism $\widetilde f\colon\widetilde I_1\to\widetilde I_2$ between
injective $k\widetilde Q$-modules such that the morphism $\widetilde P_1\to\widetilde P_2$
equals $\widetilde{\theta}(\widetilde f)$ (here $\widetilde{\theta}$ is as in \ref{lem3}).
 Let $f\colon
I_1\to I_2$ be the image of $\widetilde f$ under $q_{\lambda}\colon \ts{mod}\, k\widetilde Q\to
\ts{mod}\, kQ$. Hence, the image of $\ts{Ker} \widetilde f$ under
$q_{\lambda}\colon \ts{mod}\, k\widetilde Q\to \ts{mod}\, kQ$ is $\ts{Ker}  f$. Let $P_1\to P_2$ be the image of
$\widetilde{\theta}(\widetilde f)$ under $p_{\lambda}\colon\ts{mod}\,\widetilde{\Sigma}\to\ts{mod}\,\Sigma$. Therefore the
commutativity of the diagram in \ref{lem3} and the fact that
$\ts{Hom}_1(\Sigma,R)$ is the image of $\widetilde{R}$ under
$p_{\lambda}\colon\ts{mod}\,\widetilde{\Sigma}\to\ts{mod}\,\Sigma$ gives
a minimal projective presentation in $\ts{mod}\, \Sigma$:
\begin{equation}
  0\to P_1\to P_2\to \ts{Hom}_B(\Sigma,R)\to 0\ \ .\notag
\end{equation}
On the other hand, \ref{lem3} shows that $P_1\to P_2$ is equal
to the following morphism in $\ts{mod}\, \Sigma$:
\begin{equation}
  \ts{Hom}_B(\Sigma,\ts{Hom}_{kQ}(T,I_1))\xrightarrow{\ts{Hom}_B(\Sigma,\ts{Hom}_{kQ}(T,f))}
  \ts{Hom}_B(\Sigma,\ts{Hom}_{kQ}(T,I_2))\ \ .\notag
\end{equation}
Therefore we have a minimal $\ts{add}\, \mathcal L_{\Sigma}$-presentation:
\begin{equation}
  \ts{Hom}_{kQ}(T,I_1)\xrightarrow{\ts{Hom}_{kQ}(T,f)}\ts{Hom}_{kQ}(T,I_2)\to R\ \ .\notag
\end{equation}
Because of \ref{lem11}, the sequence $0\to \ts{Hom}_{kQ}(T,I_1)\to \ts{Hom}_{kQ}(T,I_2)\to
R\to 0$ is exact and
 $\ts{Ker} f=\ts{Tor}_1^{B}(R,T)$. In other words, 
$q_{\lambda}\colon \ts{mod}\, k\widetilde Q\to \ts{mod}\, kQ$ maps $\ts{Ker} \widetilde f$ to
$\ts{Tor}_1^{B}(R,T)$.  Using \ref{lem11} and the last diagram in the proof of
\ref{lem3}, we get $F_{\lambda}(\ts{Ext}_{k\widetilde Q}^1(\widetilde T,\ts{Ker} \widetilde f))=R$.

We give an explicit construction of $\varphi$. Let
$M\in k(\widetilde{\Gamma})$. We fix a minimal projective presentation in
$\ts{mod}\, \widetilde{\Sigma}$:
\begin{equation}
0\to  \widetilde P_1\xrightarrow{\widetilde u} \widetilde P_2\to \widetilde R\to 0\ \ ,\notag
\end{equation}
and injective $k\widetilde Q$-modules $\widetilde I_1$ and $\widetilde I_2$, together with a
morphism $\widetilde f\colon \widetilde I_1\to \widetilde I_2$ such that $\widetilde u=\widetilde{\theta}(\widetilde
f)$. Then we let $\varphi(M)$ be the
following $\widetilde B$-module:
\begin{equation}
  \varphi(M)=\varphi_0(\widetilde E)\oplus \ts{Ext}_{k\widetilde Q}^1(\widetilde T ,\ts{Ker} \widetilde
  f)\ \ ,\notag
\end{equation}
where $\varphi_0(\widetilde E)=\varphi_0(\widetilde E_1)\oplus\cdots\oplus \varphi_0(\widetilde E_s)$
if $\widetilde E=\widetilde E_1\oplus\cdots\oplus\widetilde E_s$ with $\widetilde E_1,\ldots,\widetilde E_s\in\widetilde{\Sigma}$.
This finishes the construction of the map $\varphi\colon
\widetilde{\Gamma}_o\to\ts{mod}\,\widetilde B$. We now prove the $G$-equivariance
property. Let $M\in k(\widetilde{\Gamma})$ be a vertex and let $g\in G$. We keep the
above notation $\widetilde R,\,\widetilde E,\, etc.$ introduced for $M$, and we
adopt the dashed notation $\widetilde R',\,\widetilde E',\, etc.$ for the corresponding
objects associated to $^gM$. We have
$k(\widetilde{\Gamma})(\widetilde{\Sigma},\,^gM)=\,^gk(\widetilde{\Gamma})(\widetilde{\Sigma},M)$.
Indeed, the $\widetilde{\Sigma}$-modules $k(\widetilde{\Gamma})(\widetilde{\Sigma},\,^gM)$ and
$^gk(\widetilde{\Gamma})(\widetilde{\Sigma},M)$ are given by
the functors $X\mapsto k(\widetilde{\Gamma})(X,\,^gM)$ and
$X\mapsto k(\widetilde{\Gamma})(\,^{g^{-1}}X,M)$ from $\widetilde{\Sigma}^{op}$ to $\ts{mod}\, k$, respectively. These
two functors coincide because $G$ acts on $k(\widetilde{\Gamma})$.
 Hence,
$\widetilde E'=\,^g\widetilde E$ and $\widetilde R'=\,^g\widetilde R$. Therefore any minimal
projective presentation of $\widetilde R'$ in $\ts{mod}\,\widetilde{\Sigma}$ is obtained from a
minimal projective presentation of $\widetilde R$ by applying $g$. Since,
moreover, $\widetilde{\theta}$ is $G$-equivariant (see
\ref{lem3}), we deduce that $\widetilde f'=\,^g\widetilde f$. Finally, the
$G$-action on $\ts{mod}\, k\widetilde Q$ implies, as above, that $\ts{Ext}^1_{k\widetilde Q}(\widetilde
T,\ts{Ker}  \,^g\widetilde f)=\ts{Ext}^1_{k\widetilde Q}(\widetilde
T,\,^g\ts{Ker}  \widetilde f)=\,^g\ts{Ext}^1_{k\widetilde Q}(\widetilde
T,\ts{Ker}  \widetilde f)$. From the construction of $\varphi$, we get $\varphi(\,^gM)=\,^g\varphi(M)$. \sq

\paragraph{Construction of $\varphi$ on morphisms}$\ $

We complete the construction of $\varphi$ by proving the following lemma.
\begin{lem}
  \label{lem14}
Let $u\colon M\to N$ be a morphism in $k(\widetilde{\Gamma})$. Then there exists a
unique morphism $\varphi(u)\colon \varphi(M)\to \varphi(N)$ in
$\ts{mod}\, \widetilde B$, such that $F_{\lambda}(\varphi(u))=p(u)_{|B}$.
\end{lem}
\noindent{\textbf{Proof:}} Since $F_{\lambda}$ is exact, it is
faithful so the morphism
$\varphi(u)$ is unique. We prove its existence. By \ref{lem8}, we
have constructed $\varphi(u)=\varphi_0(u)$ in case
$N\in\widetilde{\Gamma}_{\leqslant\widetilde{\Sigma}}$. So we may assume that $N\in
\widetilde{\Gamma}\backslash\widetilde{\Gamma}_{\leqslant \widetilde{\Sigma}}$. 
Since any path in $\widetilde{\Gamma}$ from a vertex in $\widetilde{\Gamma}_{\leqslant\widetilde{\Sigma}}$ to $N$ has a vertex in $\widetilde{\Sigma}$, we may
also assume that $M\in\left(\widetilde{\Gamma}\backslash\widetilde{\Gamma}_{\leqslant
    \widetilde{\Sigma}}\right)\bigcup \widetilde{\Sigma}$. The functor
$\varphi_0\colon k(\widetilde{\Gamma}_{\leqslant \widetilde{\Sigma}})\to \ts{mod}\, \widetilde B$ naturally
extends to a unique functor $\varphi_0\colon \ts{add}\,(k(\widetilde{\Gamma}_{\leqslant
  \widetilde{\Sigma}}))\to \ts{mod}\, \widetilde B$, such that the following diagram
commutes:
\begin{equation}
  \xymatrix{
\ts{add}\,(k(\widetilde{\Gamma}_{\leqslant
  \widetilde{\Sigma}})) \ar@{->}[rr]^{\varphi_0} \ar@{->}[d]_{\ts{add}\, p} && \ts{mod}\, \widetilde
B \ar@{->}[d]^{F_{\lambda}} & \\
\ts{add}\,(\ts{ind}\, \Gamma_{\leqslant \Sigma}) \ar@{^(->}[rr] && \ts{mod}\, B & .
}\notag
\end{equation}
We fix decompositions in $\ts{mod}\, \widetilde{\Sigma}$ as in the proof of
\ref{lem12}:
\begin{equation}
  k(\widetilde{\Gamma})(\widetilde{\Sigma},M)=\widetilde P\oplus\widetilde R,\ \ \text{and}\ \
  k(\widetilde{\Gamma})(\widetilde{\Sigma},N)=\widetilde P'\oplus\widetilde R'\ \ ,\notag
\end{equation}
where $\widetilde P,\widetilde P'$ are projective and $\widetilde R,\widetilde R'$ have no non-zero projective
direct summand. We let $\widetilde E,\widetilde E'\in\ts{add}\,\widetilde{\Sigma}$ be such that $\widetilde
P=k(\widetilde{\Gamma})(\widetilde{\Sigma},\widetilde E)$ and $\widetilde P'=k(\widetilde{\Gamma})(\widetilde{\Sigma},\widetilde E')$,
respectively. Therefore the morphism $k(\widetilde{\Gamma})(\widetilde{\Sigma},u)$ can be
written as:
\begin{equation}
  k(\widetilde{\Gamma})(\widetilde{\Sigma},u)=\begin{bmatrix} \widetilde u_1 & 0\\
    \widetilde u_2 & \widetilde
    u_3\end{bmatrix} \colon k(\widetilde{\Gamma})(\widetilde{\Sigma},\widetilde
  E)\oplus\widetilde R\to k(\widetilde{\Gamma})(\widetilde{\Sigma},\widetilde
  E')\oplus \widetilde R'\
  \ .\notag
\end{equation}
Similarly, we fix decompositions in $\ts{mod}\, B$:
\begin{equation}
  p(M)_{|B}=E\oplus R,\ \ p(N)_{|B}=E'\oplus R'\ \ ,\notag
\end{equation}
where $E,E'\in \ts{add}\, \Sigma$, and $R,R'$ have no direct summand in $\Sigma$. As
above, the morphism $p(u)_{|B}$ decomposes as:
\begin{equation}
  p(u)_{|B}=\begin{bmatrix} u_1 & 0\\ u_2&u_3\end{bmatrix} \colon
  E\oplus R\to E'\oplus R'\ \ .\notag
\end{equation}
Recall from the proof of \ref{lem12} that 
$p_{\lambda}\colon\ts{mod}\, \widetilde{\Sigma}\to \ts{mod}\, \Sigma$ maps $k(\widetilde{\Gamma})(\widetilde{\Sigma},\widetilde E)$, $\widetilde R$, $k(\widetilde{\Gamma})(\widetilde{\Sigma},\widetilde E')$ and $\widetilde R'$ to
$\ts{Hom}_B(\Sigma,E)$, $\ts{Hom}_B(\Sigma,R)$, $\ts{Hom}_B(\Sigma,E')$ and
$\ts{Hom}_B(\Sigma,R')$, respectively. As a consequence, it maps $\widetilde u_i$ to
$\ts{Hom}_B(\Sigma,u_i)$, for every $i$. As in
the proof of \ref{lem12}, we have morphisms $\widetilde f\colon \widetilde I_1\to \widetilde I_2$ and
$\widetilde f'\colon \widetilde I_1'\to \widetilde I_2'$ between
injective $k\widetilde Q$-modules and minimal projective
presentations in $\ts{mod}\, \widetilde{\Sigma}$:
\begin{equation}
   \widetilde{\theta}(\widetilde I_1)\xrightarrow{\widetilde{\theta}(\widetilde f)}\widetilde{\theta}(\widetilde I_2)
  \to \widetilde R\to 0\ \  
\text{and}\ \ \widetilde{\theta}(\widetilde I_1')\xrightarrow{\widetilde{\theta}(\widetilde f')}\widetilde{\theta}(\widetilde I_2')
  \xrightarrow{\widetilde v} \widetilde R'\to 0\ \ .\notag
\end{equation}
With these notations, we have:
\begin{equation}
  \varphi(M)=\varphi_0(\widetilde E)\oplus \ts{Ext}_{k\widetilde Q}^1(\widetilde{T},\ts{Ker} \widetilde
  f)\ \ 
\text{and} \  \varphi(N)=\varphi_0(\widetilde E')\oplus \ts{Ext}_{k\widetilde
  Q}^1(\widetilde{T},\ts{Ker} \widetilde f')\ \ .\notag
\end{equation}
Also, if $M\in \widetilde{\Sigma}$, then $\widetilde f=0$, so
$\varphi(M)=\varphi_0(\widetilde E)$.

It suffices to prove that each of $u_1,u_2,u_3$ is the image under
$F_{\lambda}$ of some morphism $\varphi_0(\widetilde E)\to \varphi_0(\widetilde E')$,
$\varphi_0(\widetilde E)\to \ts{Ext}_{k\widetilde
  Q}^1(\widetilde T,\ts{Ker} \widetilde f')$ and $\ts{Ext}_{k\widetilde Q}^1(\widetilde T,\ts{Ker} \widetilde
f)\to \ts{Ext}_{k\widetilde Q}^1(\widetilde T,\ts{Ker} \widetilde f')$, respectively. Clearly,
$u_1\colon k(\widetilde{\Gamma})(\widetilde{\Sigma},\widetilde E)\to
k(\widetilde{\Gamma})(\widetilde{\Sigma},\widetilde E')$ is induced by
a
morphism $\widetilde E\to \widetilde E'$ in $\ts{add}\,\widetilde{\Sigma}$. This and
\ref{lem8} imply that $u_1$ is the image under $F_{\lambda}$ of a
morphism $\varphi_0(\widetilde E)\to\varphi_0(\widetilde E')$. We
now
prove that $u_2$ is the image under $F_{\lambda}$ of a morphism
$\varphi_0(\widetilde E)\to\ts{Ext}_{k\widetilde Q}^1(\widetilde T,\ts{Ker} \widetilde f')$. Let $f'\colon I_1'\to I_2'$ be the image of
$\widetilde f'\colon\widetilde I_1'\to\widetilde I_2'$ under $q_{\lambda}\colon\ts{mod}\, k\widetilde Q\to \ts{mod}\,
kQ$. Therefore we have a minimal projective
presentation in $\ts{mod}\, \Sigma$ (see \ref{lem3} and the proof of \ref{lem12}):
\begin{equation}
  0\to \theta(I_1')\xrightarrow{\theta(f')}\theta(I_2')\to
  \ts{Hom}_B(\Sigma,R')\to 0 \ ,\notag
\end{equation}
together with a minimal injective copresentation in $\ts{mod}\, kQ$:
\begin{equation}
0\to \ts{Tor}_1^{B}(R',T)\to I_1'\xrightarrow{f'}I_2'\to 0\ \ .\notag  
\end{equation}
Recall that $\ts{Tor}_1^{B}(R',T)$ is equal to the image of $\ts{Ker} \widetilde f'$
  under $q_{\lambda}\colon\ts{mod}\, k\widetilde Q\to \ts{mod}\, kQ$.  Therefore we have an
  exact sequence in $\ts{mod}\, B$, which is also a minimal
  $\ts{add}\,\mathcal L_{\Sigma}$-presentation:
  \begin{equation}
    0\to \ts{Hom}_{kQ}(T,I_1')\to \ts{Hom}_{kQ}(T,I_2')\xrightarrow{v}R'\to 0\ \ ,\notag
  \end{equation}
where $v$ is such that $\ts{Hom}_B(\Sigma,v)$ is the image of
$\widetilde v\colon \widetilde{\theta}(\widetilde I_2')\to\widetilde R'$ under 
$p_{\lambda}\colon\ts{mod}\, \widetilde{\Sigma}\to \ts{mod}\, \Sigma$ (see the diagram in \ref{lem3}).
The projective cover $\widetilde v$ of $\widetilde R'$ in $\ts{mod}\, \widetilde{\Sigma}$ yields a
morphism $\widetilde{\delta}\colon\widetilde I\to \widetilde I_2'$ in $\ts{mod}\, k\widetilde Q$, where $\widetilde I$
is the injective $k\widetilde Q$-module such that $\widetilde P=\widetilde{\theta}(\widetilde I)$, and
such that the following diagram of $\ts{mod}\, \widetilde{\Sigma}$ commutes:
\begin{equation}
  \xymatrix{
&\widetilde P \ar@{->}[ld]_{\widetilde{\theta}(\widetilde{\delta})} \ar@{->}[d]^{\widetilde u_2} \\
\widetilde{\theta}(\widetilde I_2') \ar@{->}[r]_{\widetilde v} & \widetilde R' & \ \ .
}\notag
\end{equation}
Therefore if $\delta\colon I\to I_2'$ denotes the image of
$\widetilde{\delta}\colon \widetilde I\to \widetilde I_2'$ under
$q_{\lambda}\colon\ts{mod}\, k\widetilde Q\to \ts{mod}\, kQ$, then $\ts{Hom}_B(\Sigma,u_2)$ equals the
composition $\theta(I)\xrightarrow{\theta(\delta)}
\theta(I_2')\xrightarrow{\ts{Hom}_B(\Sigma,v)} \ts{Hom}_B(\Sigma,R')$. This is an
equality of morphisms in $\ts{mod}\,\Sigma$, hence, of morphisms between
contravariant functors from $\ts{add}\,\Sigma$ to $\ts{mod}\, k$. Applying
 this equality to $E$ yields that $u_2$ equals the
composition $E\xrightarrow{\ts{Hom}_{kQ}(T,\delta)}
\ts{Hom}_{kQ}(T,I_2')\xrightarrow{v} R'$. On the other hand, the morphism
  $\ts{Hom}_{kQ}(T,I_2')\xrightarrow{v} R'=\ts{Ext}_{kQ}^1(T,\ts{Ker}
  f')$
is the connecting morphism of the sequence resulting from the
application of $\ts{Hom}_{kQ}(T,-)$ to the exact sequence $0\to \ts{Ker}  f'\to
I_1'\xrightarrow{f'} I_2'\to 0$. Therefore \ref{lem3} implies
that $v$  equals the image under $F_{\lambda}$ of the 
connecting morphism of the sequence resulting from the application of
$\ts{Hom}_{k\widetilde Q}(\widetilde T,-)$ to the exact sequence $0\to \ts{Ker} \widetilde f'\to
\widetilde I_1'\to \widetilde I_2'\to 0$. Consequently, $u_2$ equals the image
under $F_{\lambda}$ of the 
composition
$\varphi_0(\widetilde E)\xrightarrow{\ts{Hom}_{k\widetilde Q}(\widetilde T,\delta)}
\varphi_0(\ts{Hom}_{kQ}(\widetilde T,\widetilde I_2'))\to \ts{Ext}_{k\widetilde Q}^1(\widetilde T,\ts{Ker} \widetilde f ')$.\\

It remains to prove that $u_3\colon R\to R'$  equals  the
image under $F_{\lambda}$ of a morphism $\ts{Ext}_{k\widetilde Q}^1(\widetilde T,\ts{Ker} \widetilde
f)\to \ts{Ext}_{k\widetilde Q}^1(\widetilde T,\ts{Ker} \widetilde
f')$ in $\ts{mod}\, \widetilde B$. Using the projective
presentations of $\widetilde R$ and $\widetilde R'$, we find
morphisms $\widetilde{\alpha}\colon \widetilde I_2\to \widetilde I_2'$ and $\widetilde{\beta}\colon \widetilde
I_1\to \widetilde I_1'$ such that the following diagram commutes:
\begin{equation}
  \xymatrix{
0 \ar@{->}[r] &  \widetilde\theta(\widetilde I_1) \ar@{->}[r]^{\widetilde{\theta}(\widetilde f)}
\ar@{->}[d]_{\widetilde{\theta}(\widetilde{\beta})} & \widetilde{\theta}(\widetilde I_2) \ar@{->}[r]
\ar@{->}[d]_{\widetilde{\theta}(\widetilde{\alpha})}& \widetilde R 
\ar@{->}[r] \ar@{->}[d]_{\widetilde u_3}& 0\\
0 \ar@{->}[r] & \widetilde{\theta}(\widetilde I_1')
\ar@{->}[r]_{\widetilde{\theta}(\widetilde f')} &
\widetilde{\theta}(\widetilde I_2')
\ar@{->}[r] & \widetilde R' \ar@{->}[r] & 0 &\ \ .
}\notag
\end{equation}
Therefore there exists a morphism $\widetilde{\gamma}\colon \ts{Ker} \widetilde f\to \ts{Ker} \widetilde
f'$ making the following diagram in $\ts{mod}\, k\widetilde Q$ commute:
\begin{equation}
  \xymatrix{
0\ar@{->}[r] & \ts{Ker} \widetilde f \ar@{->}[r] \ar@{->}[d]_{\widetilde{\gamma}} & \widetilde I_1
\ar@{->}[r]^{\widetilde f} \ar@{->}[d]_{\widetilde{\beta}} & \widetilde I_2
\ar@{->}[d]_{\widetilde{\alpha}} \ar@{->}[r] & 0\\
0\ar@{->}[r] & \ts{Ker} \widetilde f' \ar@{->}[r] & \widetilde I_1'
\ar@{->}[r]_{\widetilde f'} & \widetilde
I_2' \ar@{->}[r] & 0 & \ \ .
}\notag
\end{equation}
We claim that the image of $\ts{Ext}_{k\widetilde
  Q}^1(\widetilde T,\widetilde{\gamma})\colon \ts{Ext}_{k\widetilde Q}^1(\widetilde T,\ts{Ker} \widetilde
f)\to \ts{Ext}_{k\widetilde Q}^1(\widetilde T,\ts{Ker} \widetilde f')$
under $F_{\lambda}$  equals
 $u_3$. Indeed, let $\alpha,\beta,\gamma$ be the respective images of
$\widetilde{\alpha},\widetilde{\beta},\widetilde{\gamma}$ under $q_{\lambda}\colon\ts{mod}\, k\widetilde
Q\to \ts{mod}\, kQ$. Then the image of $\ts{Ext}_{k\widetilde
  Q}^1(\widetilde T,\widetilde{\gamma})$ under $F_{\lambda}$ is equal to (see
\ref{lem3}):
\begin{equation}
  \ts{Ext}_{kQ}^1(T,\gamma)\colon \ts{Ext}_{k Q}^1(T,\ts{Ker} 
f)\to \ts{Ext}_{k Q}^1(T,\ts{Ker}  f')\ \ .\notag
\end{equation}
 On the other hand, we have two commutative diagrams in $\ts{mod}\, kQ$ and
 $\ts{mod}\, \widetilde B$ respectively:
 \begin{equation}
\xymatrix{
0 \ar@{->}[r] & \ts{Ker}  f = \ts{Tor}_1^{B}(R,T) \ar@{->}[r]
\ar@{->}[d]_{\gamma} & I_1 \ar@{->}[r] \ar@{->}[d]_{\beta} &
I_2\ar@{->}[r] \ar@{->}[d]_{\alpha} & 0\\
0\ar@{->}[r] & \ts{Ker}  f'=\ts{Tor}_1^{B}(R',T) \ar@{->}[r] &I_1' \ar@{->}[r] & I_2'
\ar@{->}[r] & 0 & ,\ \text{and}\\
0\ar@{->}[r] & \ts{Hom}_{kQ}(T,I_1)\ar@{->}[r]  \ar@{->}[d]_{\ts{Hom}_{kQ}(T,\beta)} & \ts{Hom}_{kQ}(T,I_2)
\ar@{->}[r] \ar@{->}[d]_{\ts{Hom}_{kQ}(T,\alpha)} & R \ar@{->}[r]
\ar@{->}[d]_{u_3} & 0 \\
0 \ar@{->}[r] & \ts{Hom}_{kQ}(T,I_1') \ar@{->}[r] & \ts{Hom}_{kQ}(T,I_2') \ar@{->}[r] & R' \ar@{->}[r] &
0 & \ \ ,
}\notag
 \end{equation}
from which it is straightforward to check that $u_3\colon R\to R'$
coincides with $\ts{Ext}^1_{kQ}(T,\gamma)$. Thus, $u_3$ is equal to the
image under $F_{\lambda}$ of the morphism $\ts{Ext}^1_{k\widetilde
  Q}(\widetilde T,\widetilde{\gamma})\colon \ts{Ext}^1_{k\widetilde
  Q}(\widetilde T,\ts{Ker} \widetilde f)\to \ts{Ext}^1_{k\widetilde
  Q}(\widetilde T,\ts{Ker} \widetilde f')$. This completes the proof.\sq

We summarise our results in the following theorem.
\begin{thm}
  \label{prop2}
Let $A$ be a finite dimensional $k$-algebra and $\Gamma$ be a component of
$\Gamma(\ts{mod}\, A)$ containing a left section $\Sigma$. Let
$B=A/\ts{Ann}\,\Sigma$ and $\pi\colon \widetilde{\Gamma}\to\Gamma$ be 
a Galois covering with group $G$ of translation quivers such that
there exists a well-behaved functor
$p\colon k(\widetilde{\Gamma})\to\ts{ind}\,\Gamma$. Then there exists a covering 
$F\colon\widetilde B\to B$ with $\widetilde B$ locally bounded and a functor
$\varphi\colon k(\widetilde{\Gamma})\to\ts{mod}\,\widetilde B$ which is $G$-equivariant on
vertices and makes the following diagram commute:
\begin{equation}
  \xymatrix{
k(\widetilde{\Gamma}) \ar@{->}[rr]^{\varphi} \ar@{->}[d]_p && \ts{mod}\, \widetilde B
\ar@{->}[d]^{F_{\lambda}}\\
\ts{ind}\, \Gamma \ar@{->}[rr]_{\ts{Hom}_A(B,-)} && \ts{mod}\, B & \ \ .
}\notag
\end{equation}
\end{thm}
\noindent{\textbf{Proof:}} The functor $\varphi$ is constructed as
above. The $G$-equivariance on vertices
follows from \ref{lem12} and \ref{lem14}.\sq

\begin{cor}
  \label{cor1.5.1}
If $p\colon k(\widetilde{\Gamma})\to \ts{ind}\, \Gamma$ is a Galois covering (with respect to
the action of
$G$ on $k(\widetilde{\Gamma})$), then the functor
$\varphi\colon k(\widetilde{\Gamma})\to \ts{mod}\, \widetilde B$
of \ref{prop2} is $G$-equivariant.
\end{cor}
\noindent{\textbf{Proof:}} We already know that $\varphi$ is
$G$-equivariant on objects. Also $F\colon\widetilde B\to B$
is a Galois covering with group $G$ (see \ref{lem2}). Let $f\colon M\to N$ be a morphism
in $k(\widetilde{\Gamma})$, and $g\in G$. Then $\varphi(\,^gf)\colon
\varphi(\,^gM)\to \varphi(\,^gN)$ and $^g\varphi(f)\colon
\,^g\varphi(M)\to \,^g\varphi(N)$ are two morphisms in $\ts{mod}\, \widetilde
B$ such that
$F_{\lambda}(\varphi(\,^gf))=p(\,^gf)_{|B}=p(f)_{|B}=F_{\lambda}(\,^g\varphi(f))$
(recall that $F_{\lambda}=F_{\lambda}\circ g$ for every $g\in G$ because it is the
push-down functor of a Galois covering with group $G$). We deduce that
$\varphi(\,^gf)=\,^g\varphi(f)$.\sq

\section{The main theorem}$\ $
\label{sec:s5}

In this section we prove Theorem~\ref{thm2}. Assume
that $A$
is laura  with connecting components. We
 use the following notation:
\begin{enumerate}
\item[-] $\Gamma$ is the connecting component of $\Gamma(\ts{mod}\, A)$ (if $A$ is
  concealed we choose $\Gamma$ to be the unique postprojective component), and $\pi\colon
  \widetilde{\Gamma}\to \Gamma$ is a Galois covering with group $G$ of
  translation quivers such that 
there exists a well-behaved covering functor  $p\colon k(\widetilde{\Gamma})\to
\ts{ind}\, \Gamma$.  If $\Gamma$ is standard, we assume that $p$ equals
  the composition of $k(\pi)\colon k(\widetilde{\Gamma})\to k(\Gamma)$ with some
  isomorphism $k(\Gamma)\xrightarrow{\sim} \ts{ind}\, \Gamma$, so that $p$ is a
  Galois covering with group $G$.
\item[-] $\Sigma$ is the full subcategory of $\ts{ind}\, \Gamma$ whose objects are the
  $\ts{Ext}$-injective objects in $\mathcal L_A$.
\item[-] $B$ is the left support of $A$, that is, $B$ is the
  endomorphism algebra of the direct sum of the indecomposable
  projective modules lying on $\mathcal L_A$ (see Section~\ref{sec:s1}).
\end{enumerate}
Because of \cite[4.4, 5.1]{act}, the algebra $B$ is a product of tilted
algebras. Without loss of generality, we assume that:
\begin{enumerate}
\item[-] $B=\ts{End}_{kQ}(T)$, where $T=T_1\oplus\cdots\oplus T_n$ is a
  multiplicity-free tilting $kQ$-module ($T_i\in\ts{ind}\, kQ$).
\item[-] $\Sigma$ is the full subcategory of $\ts{mod}\, B$ with objects the
  modules of the form $\ts{Hom}_{kQ}\left(T,D(kQe_x)\right)$, $x\in Q_0$.
\end{enumerate}

It follows from \cite[2.1 Ex. b]{assem} that $\Sigma$ is a left section of $\Gamma$. So
we may apply \ref{prop2}.
 The proof of Theorem~\ref{thm2} is
done in the following steps: We first construct a locally bounded
$k$-category $\widetilde A$ endowed with a free
$G$-action in case $A$ is standard; then we construct a
covering functor $F\colon\widetilde A\to A$ extending the functor $F\colon\widetilde
B\to B$ of \ref{prop2} and satisfying the conditions of the
theorem; we also construct a functor $\Phi\colon k(\widetilde{\Gamma})\to\ts{mod}\,\widetilde A$
which extends the functor $\varphi\colon k(\widetilde{\Gamma})\to\ts{mod}\,\widetilde B$ of
\ref{prop2}; and finally we prove Theorem~\ref{thm2}.

\paragraph{The category $\widetilde A$}$\ $

 We need some
 notation. Let $C$ be the full subcategory of $\ts{ind}\, A$
with objects the indecomposable projective $A$-modules not in
$\mathcal L_A$. So $C$ is a full subcategory of $\ts{ind}\,\Gamma$. Let $\widetilde C$ be the
full subcategory $p^{-1}(C)$, so that $p$ induces a covering functor
$\widetilde C\to C$. If $A$ is standard and $p$ is Galois  with
group $G$, then $p\colon \widetilde C\to 
C$ is a Galois covering with group $G$. For every $x\in
\widetilde B_o$, let $\widetilde P_x$ be the corresponding
indecomposable projective $\widetilde B$-module. Also, $P_x\in
\ts{mod}\, B$ denotes the indecomposable projective $B$-module
associated to an object $x\in B_o$. We define the $\widetilde C-\widetilde
B$-bimodule $\widetilde M$ to be the functor $\widetilde C\times\widetilde B^{op}\to\ts{mod}\, k$
such that for every $\widetilde P\in\widetilde C_o$
and $x\in \widetilde B_o$
\begin{equation}
  _{\widetilde P}\widetilde M_{x}=\ts{Hom}_{\widetilde B}(\widetilde P_x,\varphi(\widetilde P))\ \ ,\notag
\end{equation}
with obvious actions of $\widetilde C$ (using $\varphi$) and $\widetilde B$.
The following lemma defines $\widetilde A$ and its $G$-action in case
$A$ is standard.
\begin{lem}
  \label{lem18}
Let $\widetilde A=\begin{bmatrix}\widetilde B & 0\\ \widetilde M &
  \widetilde C\end{bmatrix}$. Then $\widetilde A$ is locally bounded and
$G$ acts freely on $\widetilde A$ if $A$ is standard.
\end{lem}
\noindent{\textbf{Proof:}} 
We know that $\widetilde B$ and $\widetilde C$ are locally bounded. Let $P\in\widetilde C_o$. We
have the bijection of \ref{prop:cov}:
\begin{equation}
  \bigoplus\limits_{\widetilde x\in \widetilde B_o}\ts{Hom}_{\widetilde B}(\widetilde
  P_{\widetilde x},\varphi(\widetilde P)) \xrightarrow{\sim} \bigoplus\limits_{x\in
    B_o} \ts{Hom}_{B}(P_x,p(\widetilde P)_{|B})\ \ .\tag{i}
\end{equation}
Since the right-hand side is finite dimensional, then
so is $\bigoplus\limits_{\widetilde x\in \widetilde B_o}\ _{\widetilde P}\widetilde M_{\widetilde x}$, for every $\widetilde
P\in\widetilde C_o$.

Now let $P\in\widetilde C_o$, let $\widetilde x\in \widetilde B_o$, and let us prove that
$\bigoplus\limits_{p(P')=p(P)}\ts{Hom}_{\widetilde B}(\widetilde P_{\widetilde x},\varphi(P'))$
is finite dimensional. By definition of $p$, we have
$p^{-1}(p(P))=\{\,^gP\ |\ g\in G\}$. Also, we know from
\ref{prop2} that $\varphi$ is $G$-equivariant on
objects. Therefore:
\begin{equation}
  \bigoplus\limits_{p(P')=p(P)}\ts{Hom}_{\widetilde B}(\widetilde
  P_{\widetilde x},\varphi(P'))=\bigoplus\limits_{g\in G}\ts{Hom}_{\widetilde B}(\widetilde
  P_{\widetilde x},\,^g\varphi(P))=\bigoplus\limits_{g\in
    G}\ts{Hom}_{\widetilde{B}}(\,^{g^{-1}}\widetilde{P}_{\widetilde
    x},\varphi(P))\tag{ii}
\end{equation}
where the last equality follows from the $G$-action on
$\ts{mod}\,\widetilde{B}$. Applying
\ref{prop:cov} to  $\widetilde P_{\widetilde x}$ yields a bijection of vector
spaces:
\begin{equation}
  \bigoplus\limits_{g\in G}\ts{Hom}_{\widetilde B}(\,^{g^{-1}}\widetilde
  P_{\widetilde x},\varphi(P))\simeq \ts{Hom}_B(P_{F(\widetilde x)},F_{\lambda}\varphi(P))\ .\tag{iii}
\end{equation}
From (ii) and (iii) we infer that $\bigoplus\limits_{p(P')=p(P)}\ts{Hom}_{\widetilde B}(\widetilde
  P_{\widetilde x},\varphi(P'))$ is finite dimensional for every $\widetilde x\in \widetilde B_o$
  and  $P\in\widetilde C_o$.
 This shows that $\widetilde A$ is locally bounded.

Assume now that $A$ is standard and that $p\colon k(\widetilde{\Gamma})\to \ts{ind}\, \Gamma$ is a
Galois covering with group $G$. We define a free
$G$-action on $\widetilde A$. We already have a free
$G$-action on $\widetilde B$ and on $\widetilde C$. Also, for every
$\widetilde x\in \widetilde B_o$, $\widetilde P\in\widetilde C_o$ and $g\in
G$, we have an isomorphism of vector spaces:
\begin{equation}
  _{\widetilde P}\widetilde M_{\widetilde x}=\ts{Hom}_{\widetilde B}(\widetilde P_{\widetilde x},\varphi(P))\xrightarrow{\sim} \
  _{^g\widetilde P}\widetilde M_{g\widetilde x}=\ts{Hom}_{\widetilde B}(\,^g\widetilde
  P_{\widetilde x},\varphi(\,^g\widetilde P))\ \tag{$\star$}
\end{equation}
given by the $G$-action on $\ts{mod}\, \widetilde B$ (recall that
$\varphi$ is $G$-equivariant on objects, and that $\widetilde
P_{g\widetilde x}=\,^g\widetilde P_{\widetilde x}$). 
We define the action of $g$ on morphisms of $\widetilde A$ lying in $\widetilde M$
using this isomorphism. Since $G$ acts on $\ts{mod}\,\widetilde B$,
this defines a $G$-action on $\widetilde A$, that is, $g(vu)=g(v)g(u)$
whenever $u$ and $v$ are composable in $\widetilde A$. Moreover,
$G$ acts freely on objects in $\widetilde B$ and in $\widetilde C$. So we have a
free $G$-action on $\widetilde A$.\sq

\paragraph{The functor $F\colon\widetilde A\to A$}$\ $

\begin{lem}
  \label{lem20}
There exists a covering functor $F\colon \widetilde A\to A$ extending
$F\colon\widetilde B\to B$. If moreover $A$
is standard, then $F$ can be taken to be Galois with group $G$.
\end{lem}
\noindent{\textbf{Proof:}} Note that $A=\begin{bmatrix} B & 0\\ M & C\end{bmatrix}$
where $M$ is the $C-B$-bimodule such that
$_PM_x=\ts{Hom}_{B}(P_x,P_{|B})$ for every $P\in C_o$
and $x\in B_o$.
Let us define $F\colon \widetilde A\to A$ as follows:
\begin{enumerate}
\item[-] $F_{|\widetilde B}$ coincides with the functor
  $F\colon\widetilde B\to B$.
\item[-] $F_{|\widetilde C}$ coincides with $p\colon \widetilde C\to C$.
\item[-] Let $x\in \widetilde B_o$ and $\widetilde P\in\widetilde C_o$, then
  $F\colon\ _{\widetilde P}\widetilde M_x\to \ _{F(\widetilde P)}M_{F(x)}$ is the
  following map induced by $F_{\lambda}$:
  \begin{equation}
    \ts{Hom}_{\widetilde B}(\widetilde P_x,\varphi(\widetilde P))\to
    \ts{Hom}_{B}(P_{F(x)},p(\widetilde P)_{|B})\ \ .\notag
  \end{equation}
\end{enumerate}
Since $F_{\lambda}\colon\ts{mod}\,\widetilde B\to\ts{mod}\, B$ is a functor and
$F_{\lambda}\varphi=p(-)_{|B}$ (see \ref{prop2}), we have
defined a functor $F\colon \widetilde A\to A$.
We prove that $F\colon\widetilde A\to A$
is a covering functor. Since $F \colon \widetilde B\to B$ and
$p\colon \widetilde C\to C$ are covering functors, the bijections (i),
(ii) and
(iii) in the proof of \ref{lem18} show that
for any $\widetilde{a}\in \widetilde B_o$ and any
$\widetilde{P}\in\widetilde C_o$, the two following maps induced by $F_{\lambda}$ are
isomorphisms:
\begin{align}
  & \bigoplus\limits_{F(\widetilde x)=F(\widetilde{a})}\ts{Hom}_{\widetilde B}(\widetilde
    P_{\widetilde x},\varphi(\widetilde{P}))\to \ts{Hom}_{B}(P_{F(\widetilde{a})},
    p(\widetilde{P})_{|B})\ \ ,\notag\\
  & \bigoplus\limits_{p(\widetilde Q)=p(\widetilde{P})}\ts{Hom}_{\widetilde B}(\widetilde
    P_{\widetilde{a}},\varphi(\widetilde Q))\to \ts{Hom}_{B}(P_{F(\widetilde{a})},
    p(\widetilde{P})_{|B})\ \ .\notag
\end{align}
So $F$ is a covering functor.
 Assume now that $A$ is standard. We may suppose
that $p$ is a Galois covering with group $G$. By \ref{lem18},
there is a free $G$-action on $\widetilde A$.
Moreover, $F\colon \widetilde B\to B$, and therefore $F_{\lambda}\colon\ts{mod}\,\widetilde
  B\to\ts{mod}\, B$, are $G$-equivariant, and so is $p\colon \widetilde C\to
  C$, because it restricts the Galois covering $p\colon
  k(\widetilde{\Gamma})\to\ts{ind}\,\Gamma$. Therefore $F\colon\widetilde A\to A$ is $G$-equivariant.
Finally, the fibres of $F\colon\widetilde A\to A$ on objects are the
  $G$-orbits in $\widetilde A_o$ because $F\colon\widetilde B\to B$ and
  $p\colon \widetilde C\to C$ are Galois coverings.
Since $F\colon\widetilde A\to A$ is a covering functor, this implies that
it is also a Galois covering with group $G$ (see
for instance the proof of \cite[Prop. 6.1.37]{lemeur_thesis}).\sq

\paragraph{The functor $\Phi\colon k(\widetilde{\Gamma})\to \ts{mod}\, \widetilde A$}$\ $

 We can write an $\widetilde A$-module as a triple $(K,L,f)$ where $K\in \ts{mod}\, \widetilde B$, $L\in
\ts{mod}\, \widetilde C$ and $f\colon L\underset{\widetilde C}{\otimes}\widetilde M\to K$
is a morphism of $\widetilde B$-modules.  Let  $\psi\colon k(\widetilde{\Gamma})\to \ts{mod}\, \widetilde C$ be the
functor
  $\psi\colon X\mapsto k(\widetilde{\Gamma})(\widetilde C,X)$.
Clearly, it is $G$-equivariant.
Let $L\in k(\widetilde{\Gamma})$. Then $\psi(L)\underset{\widetilde C}{\otimes}\widetilde M$ is
the $\widetilde B$-module whose value at $x\in \widetilde
B_o$ equals:
\begin{equation}
  \left(\psi(L)\underset{\widetilde C}{\otimes} \widetilde
    M\right)(x)=\left(
\bigoplus\limits_{\widetilde P\in\widetilde C_o}k(\widetilde{\Gamma})(\widetilde P,L)\underset{k}{\otimes}
\ts{Hom}_{\widetilde B}\left(\widetilde P_x,\varphi(\widetilde P)\right)
\right)/N\ \ ,\notag
\end{equation}
where $N$ is the following subspace:
\begin{equation}
  N=\left<ff'\otimes u\ -\ f\otimes\varphi(f')u\ \left|\
      \begin{array}{l}
f\in k(\widetilde{\Gamma})(\widetilde P,L),
f'\in k(\widetilde{\Gamma})(\widetilde P',\widetilde P),
u\in \ts{Hom}_{\widetilde B}\left(\widetilde P_x,\varphi(\widetilde P')\right),\
\text{for every}\ \widetilde P,\widetilde P'\in\widetilde C_o
\end{array}\right.\right>\ \ .\notag
\end{equation}

For every $x\in \widetilde B_o$ and
$\widetilde P\in\widetilde C_o$, we have a $k$-linear map:
\begin{equation}
\begin{array}{crcll}
  \eta_{L,x,P}\colon &k(\widetilde{\Gamma})(\widetilde P,L)\underset{k}{\otimes} \ts{Hom}_{\widetilde
   B}\left(\widetilde P_x,\varphi(\widetilde P)\right) & \to
  & \ts{Hom}_{\widetilde B}\left( \widetilde P_x,\varphi(L)\right)=\varphi(L)(x)\\
& f\otimes u  & \mapsto & \varphi(f)\,u&.
\end{array}\notag
\end{equation}

It is not difficult to check that the family of maps
$\left(\eta_{L,x,\widetilde P}\right)_{L,x,\widetilde P}$ defines a functorial morphism:
\begin{equation}
\eta\colon \psi(-)\underset{\widetilde C}{\otimes}\widetilde M\to \varphi\ \ .\notag
\end{equation}
Moreover, if $\varphi$ is $G$-equivariant, then so is
$\eta$. We let $\Phi\colon k(\widetilde{\Gamma})\to \ts{mod}\,
\widetilde A$ be
the following functor:
\begin{equation}
  \Phi\colon L\longmapsto (\varphi(L),\psi(L),\eta_L)\ \ .\notag
\end{equation}

\paragraph{The main theorem}$\ $

\begin{thm}
\label{thm3}
Let $A$ be laura with connecting component $\Gamma$. Let
$\pi\colon \widetilde{\Gamma}\to\Gamma$ be a Galois covering with
group $G$ such that there exists a well-behaved covering functor
$p\colon k(\widetilde{\Gamma})\to \ts{ind}\,\Gamma$. Then there
exist a covering functor $F\colon\widetilde A\to A$ where $\widetilde
A$ is connected and locally bounded, and a
commutative diagram:
\begin{equation}
  \xymatrix{
k(\widetilde{\Gamma}) \ar@{->}[r]^{\Phi} \ar@{->}[d]_p & \ts{mod}\, \widetilde A\ar@{->}[d]^{F_{\lambda}}\\
\ts{ind}\, \Gamma \ar@{^(->}[r] & \ts{mod}\, A &\ \ ,
}\notag
\end{equation}
where $\Phi$ is
faithful. If, moreover, $A$ is standard, then $F$ and $p$ may be
assumed to be Galois coverings with group $G$, and $\Phi$ is
then $G$-equivariant and full.
\end{thm}
\noindent{\textbf{Proof:}} 
The commutativity of the above diagram follows from the one of
\ref{prop2} and from that of the diagram:
\begin{equation}
  \xymatrix{
k(\widetilde{\Gamma}) \ar@{->}[rrr]^{\psi} \ar@{->}[d]_p &&& \ts{mod}\, \widetilde C \ar@{->}[d]^{p_{\lambda}}\\
\ts{ind}\, \Gamma \ar@{->}[rrr]_{X \mapsto \ts{Hom}_A(C,X)} &&& \ts{mod}\, C & \ \ ,
}\notag
\end{equation}
Since $F_{\lambda}\Phi=p$ and $p$ is faithful, then $\Phi$ is faithful.
Therefore $\Phi(k(\widetilde{\Gamma}))$ is contained in a connected
component $\Omega$ of $\ts{mod}\, \widetilde A$. 

We now prove that $\widetilde A$ is connected. Let $x\in \widetilde A_o$ and $Q_x$ be
the corresponding indecomposable projective $\widetilde A$-module. If $\widetilde
F_{\lambda}Q_x\in C_o$, then, by construction, $Q_x$ lies in the image of
$\Phi$, so that $Q_x\in \Omega$. If $F_{\lambda}Q_x\not\in C_o$, then
$F(x)\in B_o$ and $x\in \widetilde B_o$. In
this case, there is a non-zero morphism $u\colon P_{F(x)}=\widetilde
F_{\lambda}Q_x\to E$ in $\ts{mod}\, B$, where $E\in \Sigma$. Fix $\widetilde E\in
p^{-1}(E)$ so that $F_{\lambda}\Phi(\widetilde E)=E$. 
Since $u$ is non-zero, \ref{prop:cov} implies that there
is a non-zero morphism $Q_x\to\,^g\varphi(\widetilde E)=\Phi(\,^g\widetilde E)$ in
$\ts{mod}\,\widetilde B$ (recall that $\varphi$ is $G$-equivariant on vertices). So
$Q_x\in \Omega$, and
$\Omega$ contains all the indecomposable projective
$\widetilde A$-modules. This proves that $\widetilde A$ is connected.

It remains to prove that if $A$ is standard, then $\Phi$ is full,
$G$-equivariant, and $F$ is Galois with group
$G$. In case $A$ is standard, we suppose that $p\colon
k(\widetilde{\Gamma})\to \ts{ind}\, \Gamma$ is Galois  with group
$G$. Therefore $\varphi$ is $G$-equivariant (see
\ref{cor1.5.1}) and so is $\eta$. Hence, $\Phi$ is
$G$-equivariant. Also, $F$ is Galois because
of \ref{lem20}. We prove that $\Phi$ is full. Given a
morphism $f\colon\Phi(L)\to\Phi(N)$, there exists
$(f_g)_g\in\bigoplus\limits_{g\in G}\ts{Hom}_{k(\Gamma)}(L,\,^gN)$
such that $F_{\lambda}(f)=\sum\limits_g p(f_g)$ (because $p$ is
Galois). So $F_{\lambda}(f-\Phi(f_1))-\sum_{g\neq 1}F_{\lambda}(\Phi(f_g))=0$. Since
$F$ is Galois with group $G$ and since $\Phi$ is
$G$-equivariant, we get $f=\Phi(f_1)$. So $\Phi$ is full and
the theorem is proved.\sq

The following example of a non-standard representation-finite algebra
due to Riedtmann shows that $F$ needs not be a Galois covering.
\begin{ex}
\label{ex2}
 Assume that $\ts{char}(k)=2$ and $A$ is given by the
  bound quiver (see \cite[\S 7, Ex. 14 bis]{bog} and \cite{R83}):
   \begin{equation}
     \xymatrix{
x\ar@<2pt>@{<-}[r]^{\sigma} \ar@<-2pt>[r]_{\delta}& y
\ar@/_/@(ur,dr)^{\rho} 
},
\ \ \rho^4=0,\ \rho^2=\delta\sigma,\
     \sigma\delta=\sigma\rho\delta\ .\notag
   \end{equation}
Then $A$ is representation-finite and not standard, with the following
Auslander-Reiten quiver:
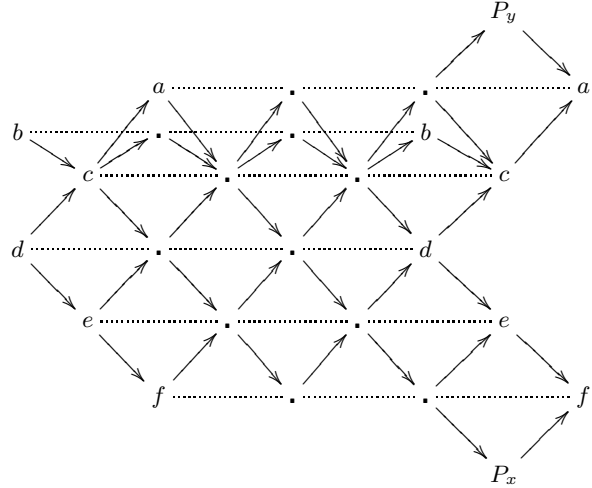
\begin{equation}
  \xymatrix@=5pt{
&& && && && && && && P_y \ar[rrdd]\\
\\
&& && a \ar[rrdd] \ar@{.}[rrrr] && && \centerdot \ar[rrdd] \ar@{.}[rrrr] && &&
\centerdot \ar[rruu] \ar@{.}[rrrr] \ar[rrdd] && && a\\
b \ar@{.}[rrrr] \ar[rrd] && && \centerdot \ar@{.}[rrrr] \ar[rrd] && &&
\centerdot \ar@{.}[rrrr] \ar[rrd] && && b \ar[rrd] \\
&& c \ar[rru] \ar[rruu] \ar@{.}[rrrr] \ar[rrdd] && && \centerdot
\ar@{.}[rrrr] \ar[rru] \ar[rruu] \ar[rrdd] && && \centerdot
\ar@{.}[rrrr] \ar[rru] \ar[rruu] \ar[rrdd] && && c \ar[rruu]\\
\\
d \ar[rruu] \ar[rrdd] \ar@{.}[rrrr] && && \centerdot \ar[rruu]
\ar[rrdd] \ar@{.}[rrrr] && && \centerdot \ar[rruu] \ar[rrdd]
\ar@{.}[rrrr] && && d \ar[rruu] \ar[rrdd]\\
\\
&& e \ar[rruu] \ar[rrdd] \ar@{.}[rrrr] && && \centerdot \ar[rruu]
\ar[rrdd] \ar@{.}[rrrr] && && \centerdot \ar[rruu] \ar[rrdd]
\ar@{.}[rrrr] && && e  \ar[rrdd]\\
\\
&& &&  f \ar[rruu]  \ar@{.}[rrrr] && && \centerdot \ar[rruu]
 \ar@{.}[rrrr] && && \centerdot \ar[rruu] \ar[rrdd]
\ar@{.}[rrrr] && && f \\
\\
&& && && && && && && P_x \ar[rruu] &&
}\notag
\end{equation}
where the two copies of $a$, $b$, $c$, $d$, $e$ and $f$, respectively,
are identified. 
In this case, there exists a well-behaved covering functor associated to the
universal cover $\widetilde{\Gamma}$ of $\Gamma(\ts{mod}\,A)$ (which
is equal to the generic covering).
Here, $G=\pi_1(\Gamma)\simeq \mathbb{Z}$
and $\widetilde A$ is the locally bounded $k$-category, given by
the following bound quiver:
\begin{equation}
  \begin{array}{c}
  \xymatrix{
{\hdots} & y_{i-1} \ar[r] \ar[rd] & y_i  \ar[r] \ar[rd] & y_{i+1} \ar[r]
\ar[rd] & y_{i+2} & \hdots
\\
{\hdots} & x_{i-1} \ar[ru]& x_i \ar[ru] & x_{i+1} \ar[ru] & x_{i+2} & \hdots
}
\\
\\
 \delta_{i+1}\sigma_i=\rho_{i+1}\rho_i,\ \sigma_{i+1}\delta_i=0\
 ,\text{for all $i$},    
  \end{array}\notag
\end{equation}
where $\sigma_i$, $\delta_i$ and $\rho_i$ denote the arrows $y_i\to
x_{i+1}$, $x_i\to y_{i+1}$, and $y_i\to y_{i+1}$, respectively.
Now the covering functor $F\colon\widetilde A\to A$ is as follows:
\begin{enumerate}
\item $F(\rho_i)=\rho$ for every $i$,
\item $F(\sigma_i)=\sigma$ for every $i\equiv 0,1\ \ts{mod}\, 4$,
\item $F(\sigma_i)=\sigma+\sigma\rho$ for every $i\equiv 2,3\ \ts{mod}\, 4$,
\item $F(\delta_i)=\delta$ for every $i\equiv 1,3\ \ts{mod}\, 4$,
\item $F(\delta_i)=\delta+\rho\delta$, for every $i\equiv 0,2\ \ts{mod}\,
  4$.
\end{enumerate}
Obviously, $F$ is a covering functor which is not Galois. Actually,
one can easily check that $A$ is simply connected, that is, the
fundamental group (in the sense of \cite{martinezvilla_delapena}) of
any presentation of $A$ is trivial.
Hence, $A$ has no proper
Galois covering by a locally bounded and connected $k$-category.
\end{ex}

The following corollary is a particular case of our main theorem. 
We
state it for later purposes.
  \begin{cor}
\label{cor2.2}
    Let $A$ be a standard laura algebra and let $\Gamma$ be a
    connecting component. There exists a Galois covering $
    F\colon\widetilde A\to A$ with group $\pi_1(\Gamma)$ and where
    $\widetilde A$ is connected and locally bounded, together with a commutative
    diagram:
    \begin{equation}
\xymatrix{
      k(\widetilde{\Gamma}) \ar@{->}[r]^{\Phi} \ar@{->}[d]_{k(\pi)} & \ts{mod}\,
      \widetilde A
    \ar@{->}[d]^{F_{\lambda}} \\
k(\Gamma) \ar@{^(->}[r] & \ts{mod}\, A&\ \ ,
}\notag
    \end{equation}
where $\pi\colon\widetilde{\Gamma}\to \Gamma$ is the universal cover and where $\Phi$
is full, faithful and $\pi_1(\Gamma)$-equivariant.
\hfill$\blacksquare$ 
 \end{cor}
\noindent{\textbf{Proof:}}
Since $\Gamma$ is a standard component, there exists an isomorphism of
categories $k(\Gamma)\to \ts{ind}\,\Gamma$ and the universal cover
$\pi\colon\widetilde{\Gamma}\to\Gamma$ induces a well-behaved functor
$k(\pi)\colon k(\widetilde{\Gamma})\to k(\Gamma)$ and therefore a
well-behaved covering functor
$k(\widetilde{\Gamma})\to\ts{ind}\Gamma$. We then apply \ref{thm3}.\sq

We pose the following problems.
\begin{pb}
  Does there exist a combinatorial characterisation of standardness
  for laura algebras (as happens for representation-finite algebras,
  see \cite{bgrs})?
\end{pb}

\begin{pb}
  Let $A$ be a left supported algebra. Is it possible to construct 
  coverings $\widetilde A\to A$ associated to the coverings of a
  component of $\Gamma(\ts{mod}\, A)$ containing the $\ts{Ext}$-injective modules
  of $\mathcal L_A$?
\end{pb}

\section{Galois coverings of the connecting component}$\ $
\label{sec:s6}

\begin{thm}
  \label{thm2.1}
Let $A$ be a standard laura algebra, and $p\colon \Gamma'\to \Gamma$ be a
Galois covering with group $G$ of a connecting component. Then
there exist a Galois covering $F'\colon A'\to A$ with group $G$,
where $A'$ is connected and locally bounded, and
 a commutative diagram:
\begin{equation}
  \xymatrix{
k(\Gamma') \ar@{->}[r]^{\Phi'} \ar@{->}[d]_{k(p)} & \ts{mod}\, A'
\ar@{->}[d]^{F_{\lambda}'} &\\
k(\Gamma) \ar@{^(->}[r]^j & \ts{mod}\, A &\ \ ,
}\notag
\end{equation}
where $\Phi'$ is full, faithful and $G$-equivariant.
\end{thm}
\noindent{\textbf{Proof:}} Since $A$ is standard, there exists a full
and faithful functor $j\colon k(\Gamma)\hookrightarrow \ts{ind}\, A$ with image
$\ts{ind}\, \Gamma$, which maps meshes to almost split sequences. Let
$\pi\colon \widetilde{\Gamma}\to \Gamma$ be the universal
cover. Then there exists a normal subgroup $H\vartriangleleft
\pi_1(\Gamma)$ such that $\widetilde{\Gamma}/H\simeq \Gamma'$ and $G\simeq \pi_1(\Gamma)/H$, and
such that under these identifications, the following diagram commutes:
\begin{equation}
  \xymatrix{
\widetilde{\Gamma} \ar@{->}[rd]^q \ar@{->}[dd]_{\pi} & \\
& \widetilde{\Gamma}/H=\Gamma' \ar@{->}[ld]^p\\
\Gamma & &\ \ ,
}\notag
\end{equation}
where $q$ is the projection. These identifications
imply that $p\colon \Gamma'\to \Gamma$ is induced by $\pi\colon\widetilde{\Gamma}\to \Gamma$ by
factoring out by $H$. By  \ref{cor2.2}, there exist a
Galois covering $F\colon\widetilde A\to A$ with group $\pi_1(\Gamma)$ and a
commutative diagram:
\begin{equation}
  \xymatrix{
k(\widetilde{\Gamma}) \ar@{->}[r]^{\Phi} \ar@{->}[d]_{k(\pi)} & \ts{mod}\, \widetilde A
\ar@{->}[d]^{F_{\lambda}} \\
k(\Gamma) \ar@{^(->}[r]^j & \ts{mod}\, A &\ \ ,
}\notag
\end{equation}
where $\Phi$ is full, faithful and
$\pi_1(\Gamma)$-equivariant.  Setting $A'=\widetilde A/H$, we deduce a Galois
covering $F'\colon A'\to A$ with group $G$ and where $A'$ is connected
and locally bounded, making
the following diagram commute:
\begin{equation}
  \xymatrix{
\widetilde A\ar@{->}[rd]^{F''} \ar@{->}[dd]_{F} &\\
& \widetilde A/H=A' \ar@{->}[ld]^{F'} \\
A & & \ \ ,
}\notag
\end{equation}
where $F''$ is the natural projection (and $F'$ is deduced from $F$
by factoring out by $H$). Therefore we have a commutative diagram of
solid arrows:
\begin{equation}
  \xymatrix{
k(\widetilde{\Gamma}) \ar@{->}[rrr]^{\Phi} \ar@{->}[dd]_{k(\pi)}
\ar@{->}[rd]^{k(q)}&&& \ts{mod}\, \widetilde A
\ar@{->}[dd]_(.3){F_{\lambda}} \ar@{->}[rd]^{F''_{\lambda}} &\\
&k(\Gamma') \ar@{->}[ld]^{k(p)} \ar@{.>}[rrr]^(.3){\Phi'} &&& \ts{mod}\, A'
\ar@{->}[ld]^{F'_{\lambda}}\\
k(\Gamma) \ar@{^(->}[rrr]^j &&& \ts{mod}\, A&&\ \ .
}\notag
\end{equation}

We prove the existence of the dotted arrow $\Phi'$ such that $\Phi'\
k(q)=F''_{\lambda}\ \Phi$. For this purpose, recall that $k(q)$ is a
Galois covering with group $H$. Hence, it suffices to prove that
$F''_{\lambda}\ \Phi$ is $H$-invariant. Indeed, we have
$F''_{\lambda}\Phi' h=F''_{\lambda}h\Phi'=F''_{\lambda}\Phi'$, for every $h\in
H$, because $\Phi$ is $\pi_1(\Gamma)$-equivariant and $F''$ is a Galois
covering with group $H$. Now, we prove that the whole
diagram commutes. We have:
\begin{equation}
  (F'_{\lambda}\ \Phi')\ k(q)=F'_{\lambda}F''_{\lambda}\ \Phi=F_{\lambda}\ \Phi=j\
  k(\pi)=j\ k(p)\ k(q)\ \ ,\notag
\end{equation}
hence, $F'_{\lambda}\ \Phi'=j\ k(p)$. We prove next that $\Phi'$ is full and
faithful. Let $f\colon X\to Y$ be a morphism in $k(\Gamma')$ such that
$\Phi'(f)=0$. Fix $\widetilde X,\widetilde Y\in k(\widetilde{\Gamma})$ such that $q(\widetilde X)=X$ and $q(\widetilde
Y)=Y$. Since $k(q)$ is Galois with group $H$, there exists
$(f_h)_{h\in H}\in \bigoplus\limits_{h\in H}k(\widetilde{\Gamma})(\widetilde X,\,^h\widetilde Y)$ such
that $\sum\limits_{h\in H} k(q)(f_h)=f$. The
commutativity of the diagram gives:
\begin{equation}
  0=\sum\limits_{h\in H} F''_{\lambda}(\Phi(f_h))\ \ ,\notag
\end{equation}
where $(\Phi(f_h))_{h\in H}\in \bigoplus\limits_{h\in
  H}\ts{Hom}_{\widetilde A}(\Phi(\widetilde X),\,^h\Phi(\widetilde Y))$ (recall that
$\Phi$ is $\pi_1(\Gamma)$-equivariant). Since $F''\colon \widetilde A\to A'$ is
Galois with group $H$, we deduce that $\Phi(f_h)=0$ for every
$h\in H$, so that $f_h=0$ for every $h\in H$, because $\Phi$ is
faithful. Thus, $f=\sum\limits_{h\in H}k(q)(f_h)=0$ and
$\Phi'$ is faithful. Let $X,Y\in k(\Gamma')$ and $u\colon \Phi'(X)\to
\Phi'(Y)$ be a morphism in $\ts{mod}\, A'$, and fix $\widetilde X,\widetilde Y\in k(\widetilde{\Gamma})$
as above. In particular, $\Phi'(X)=F''_{\lambda}(\Phi(\widetilde X))$ and
$\Phi'(Y)=F''_{\lambda}(\Phi(\widetilde Y))$. Therefore there exists
$(\widetilde u_h)_{h\in H}\in \bigoplus\limits_{h\in H} \ts{Hom}_{\widetilde A}(\Phi(\widetilde
X),\,^h\Phi(\widetilde Y))$ such that $u=\sum\limits_{h\in H} F''_{\lambda}(\widetilde
u_h)$. 
Since $\Phi$ is $\pi_1(\Gamma)$-equivariant, we have $\ts{Hom}_{\widetilde
  A}(\Phi(\widetilde X),\,^h\Phi(\widetilde Y))=\ts{Hom}_{\widetilde A}(\Phi(\widetilde
X),\Phi(\,^h\widetilde Y))$, for every $h\in H$. Since
$\Phi$ is full, there exists $(\widetilde f_h)_{h\in H}\in
\bigoplus\limits_{h\in H}k(\widetilde{\Gamma})(\widetilde X,\,^h\widetilde Y)$ such that $\widetilde
u_h=\Phi(\widetilde f_h)$ for every $h\in H$. Since $k(q)$ is Galois
with group $H$, we deduce that $\sum\limits_{h\in H} k(q)(\widetilde f_h)\in
k(\Gamma)(X,Y)$. Moreover, we have:
\begin{equation}
  \Phi'\left(\sum\limits_{h\in H} k(q)(\widetilde f_h)\right)=\sum\limits_{h\in H}
  F''_{\lambda}\Phi(\widetilde f_h)=\sum\limits_{h\in H} F''_{\lambda}\widetilde u_h=u\ \ ,\notag
\end{equation}
whence the fullness of $\Phi'$. To finish, it remains to prove that
$\Phi'$ is $G$-equivariant. Let $g\in G$ be the residual class of
$\sigma\in \pi_1(\Gamma)$ modulo $H$. We need to prove that $\Phi'\circ g=g\circ \Phi'$. We have
  $\Phi'\circ g\circ k(q)=\Phi'\circ k(q)\circ \sigma$, because $q\colon\widetilde{\Gamma}\to \Gamma'=\widetilde{\Gamma}/H$
  is the canonical projection. Hence, $\Phi'\circ g\circ k(q)=F''_{\lambda}\circ \sigma\circ
  \Phi$, because $F''_{\lambda}\circ \Phi=\Phi'\circ k(q)$, and $\Phi$ is
  $\pi_1(\Gamma)$-equivariant. Since $F''_{\lambda}\circ \sigma=g\circ F''_{\lambda}$ (because
  $F''$ is deduced from $F$ by factoring out by $H$), we have
  $\Phi'\circ g\circ k(q)=g\circ F''_{\lambda}\circ \Phi=g\circ
  \Phi'\circ k(q)$, and so
  $\Phi'\circ g=g\circ \Phi'$. The proof is complete.\sq

  \begin{cor}
    \label{cor2.3}
In the situation of Theorem~\ref{thm2.1}, the full subquiver
$\Omega$ of $\Gamma(\ts{mod}\, A')$
with vertex set equal to $\{X\in \ts{ind}\, A'\ |\ F'_{\lambda}X\in\Gamma\}$ is
a faithful and generalised standard component of
$\Gamma(\ts{mod}\, A')$, isomorphic, as a translation quiver, to
$\Gamma'$. Moreover, there exists a Galois covering of translation
quivers $\Gamma'\to \Gamma$ with group $G$ extending the map
$X\mapsto F'_{\lambda}X$.
  \end{cor}
\noindent{\textbf{Proof:}} 
Since $F'_{\lambda}\Phi'=j\,k(p)$, the module $\Phi'(X)$ is
indecomposable and lies in $\Omega$, for every $X\in\Gamma'$. On the
other hand if $X\in\Omega$,
there exists $ X'\in \Gamma'$ such that $F'_{\lambda}X=k(p)(X')$. Therefore
$F'_{\lambda}X=F'_{\lambda}\Phi'(X')$. Since $X$ and $\Phi'(X')$ are
indecomposable, there exists $g\in G$ such that
$X=\,^g\Phi'(X')=\Phi'(\,^gX')\in\Phi'(\Gamma')$. Thus, we have shown that:
\begin{enumerate}
\item[(i)] $\Omega$ coincides with the full subquiver of $\Gamma(\ts{mod}\, A')$ with
set of vertices $\{\Phi'(X)\ |\ X\in \Gamma'\}$.
\end{enumerate}
 Let
$X\xrightarrow{u}Y$ be an arrow in $\Gamma'$. Since $F'_{\lambda}\Phi'=j\,k(p)$,
then
$F'_{\lambda}\ \Phi'(u)$ is an irreducible morphism between indecomposable
$A$-modules. Using \cite[Lem. 2.1]{ws}, we deduce that $\Phi'(u)$ is
irreducible. This proves
that:
\begin{enumerate}
\item[(ii)] The full subquiver of $\Gamma(\ts{mod}\, A')$ with set of vertices
  $\{\Phi'(X)\ |\ X\in\Gamma'\}$ is contained in a connected component of $\Gamma(\ts{mod}\,
A')$.
\end{enumerate}
Combining (i), (ii) and \cite[Lem. 2.3]{ws}, we deduce that
$\Omega$ is a component of $\Gamma(\ts{mod}\, A')$. The same lemma shows that
$\Omega$ is faithful and generalised standard because so is $\Gamma$.

Let us prove that $\Phi'$ induces an isomorphism between $\Gamma'$ and
$\Omega$. 
 Since $q\colon\widetilde{\Gamma}\to\Gamma'$ is surjective on vertices and
$F''_{\lambda}\,\Phi=\Phi'\,k(q)$, then $X\in\Omega$ lies
in the image of $F''_{\lambda}$.
Also, $k(q)$ and $\Phi$ commute with the
translation, and so does $F''_{\lambda}$ (see \cite[Lem. 2.1]{ws}). Hence $\Phi'$ commutes with the
translation.
Finally $k(q)$ maps meshes to meshes, and $\Phi$ maps meshes to
  almost split sequences. So $\Phi'$ maps meshes to almost split
  sequences (see \cite[Lem. 2.2]{ws}).
Therefore there exists a morphism of translation quivers
$\Gamma'\to\Omega$ extending the map $X\mapsto \Phi'(X)$ on
vertices. Since it is a bijection on vertices, it is an
isomorphism $\Gamma'\xrightarrow{\sim}\Omega$.

 Finally, the stabiliser $G_X=\{g\in G\ |\ ^gX\simeq X\}$ of $X$ is
 trivial for every $X\in\Omega$, because $G$ acts freely on
 $\Gamma'$ and $\Phi'$ is $G$-equivariant. Therefore there exists a
Galois covering of translation quivers $\Omega\to \Gamma$ with group
$G$ and extending the map $X\mapsto F'_{\lambda}(X)$ (see \cite[3.6]{gabriel}).\sq

\begin{cor}
  \label{cor2.4}
In the situation of Theorem~\ref{thm2.1}, if $G$ is
finite, then $A'$ is a finite dimensional standard laura algebra.
\end{cor}
\noindent{\textbf{Proof:}} Since $G$ is finite, $A'$ is finite
dimensional. By the preceding corollary,   $\Gamma'$ is
generalised standard and faithful. Since $\Gamma$ has only finitely
many isomorphism classes of indecomposable modules lying on oriented
cycles, the same is true for $\Gamma'$. Therefore $\Gamma'$ is
quasi-directed and faithful. Applying \cite[3.1]{rs} (or \cite[Thm. 2]{S04a}) shows that $A'$ is 
a laura algebra with $\Gamma'$ as a connecting component. Finally,
the full and faithful functor $\Phi'\colon k(\Gamma')\to \ts{mod}\, A'$ with image
equal to $\ts{ind}\, \Gamma'$ shows that $\Gamma'$ is standard, that
is, $A'$ is standard.\sq

\begin{rem}
  The above corollary may be compared with \cite[Thm. 1.2]{alr} and
  \cite[Thm. 3]{ws}. Indeed, if $A'$ is a finite dimensional algebra
  endowed with the free action of a (necessarily finite) group $G$,
  then the category $A/G$ and the skew-group algebra $A[G]$ are
  Morita equivalent.
\end{rem}

We end this section with the following corollary:
\begin{cor}
  \label{cor2.5}
In the situation of Theorem~\ref{thm2.1}, if $G$ is
finite, then:
\begin{enumerate}
\item[(a)] $A$ is tame if and only if $A'$ is tame.
\item[(b)] $A$ is wild if and only if $A'$ is wild.
\end{enumerate}
\end{cor}
\noindent{\textbf{Proof:}} This follows from Theorem~\ref{thm2.1} and
from \cite[5.3, (b)]{ac}.\sq

  \begin{ex}
Consider the algebra $A$ of \ref{ex1}, (a). The
connecting component $\Gamma$ admits a Galois covering with group
$\mathbb{Z}/2\mathbb{Z}$ by the following translation quiver:
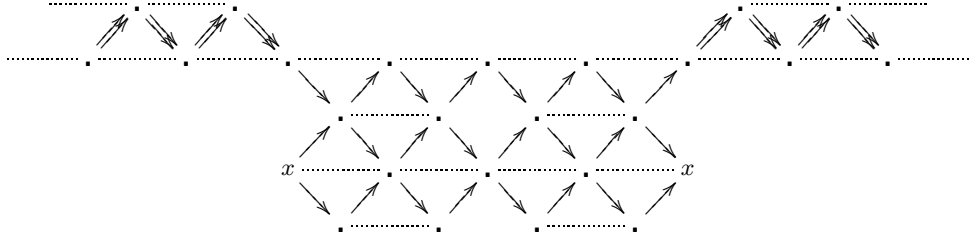
\begin{equation}
  \xymatrix@=10pt{
& \ar@{.}[rr]& &  \centerdot \ar@<1pt>[rd] \ar@<-1pt>[rd] \ar@{.}[rr]& & 
\centerdot \ar@<1pt>[rd] \ar@<-1pt>[rd] & & & & & & & & & &
\centerdot \ar@<1pt>[rd] \ar@<-1pt>[rd] \ar@{.}[rr] & & 
\centerdot \ar@<1pt>[rd] \ar@<-1pt>[rd] \ar@{.}[rr]& & &\\
\ar@{.}[rr] & & \centerdot \ar@<1pt>[ru] \ar@<-1pt>[ru] \ar@{.}[rr] & & 
\centerdot \ar@<1pt>[ru] \ar@<-1pt>[ru] \ar@{.}[rr] & &\centerdot \ar[rd] \ar@{.}[rr] & &
\centerdot \ar@{.}[rr] \ar[rd] & &
\centerdot \ar[rd] \ar@{.}[rr] &  & \centerdot \ar[rd] \ar@{.}[rr]& 
& 
\centerdot \ar@<1pt>[ru] \ar@<-1pt>[ru] \ar@{.}[rr] & &
\centerdot \ar@<1pt>[ru] \ar@<-1pt>[ru] \ar@{.}[rr] & &
 \centerdot \ar@{.}[rr]&&\\
&&&&&&& 
\centerdot \ar[ru] \ar@{.}[rr] \ar[rd]& &
\centerdot \ar[ru] \ar[rd] & & \centerdot \ar[ru] \ar[rd] \ar@{.}[rr]
& & \centerdot \ar[rd] \ar[ru] & & & & & &&\\
&&& & & & x \ar[ru] \ar@{.}[rr] \ar[rd]& & 
\centerdot \ar[ru] \ar[rd] \ar@{.}[rr] & & \centerdot \ar[ru] \ar[rd]
\ar@{.}[rr]  & & \centerdot  \ar[ru] \ar[rd] \ar@{.}[rr] & &x & &\\
&&& & & & &
\centerdot \ar[ru] \ar@{.}[rr] 
& & 
\centerdot \ar[ru] && \centerdot \ar[ru] \ar@{.}[rr] &&\centerdot \ar[ru] &&&&&&&
}\notag
\end{equation}
where the two copies of $x$ are identified. With our construction, we
get a Galois covering $F\colon A'\to A$ with group
$\mathbb{Z}/2\mathbb{Z}$, where $A'$ is the radical square zero
algebra with the following quiver:
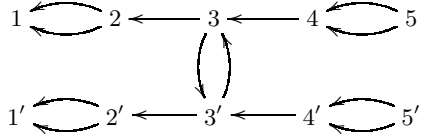
\begin{equation}
  \xymatrix{
 1 & 2 \ar@/_/[l] \ar@/^/[l]& 3 \ar[l] \ar@/_/[d] & 4 \ar[l]& 5
\ar@/_/[l] \ar@/^/[l]
\\
 1' & 2' \ar@/_/[l] \ar@/^/[l]& 3' \ar[l] \ar@/_/[u] & 4' \ar[l]& 5'
\ar@/_/[l] \ar@/^/[l]&.
}\notag
\end{equation}
Both $A$ and $A'$ are tame.
  \end{ex}

\section{Proof of Theorem~\ref{thm1}}$\ $
\label{sec:s7}

We recall the definition of the \emph{orbit graph} $\mathcal
O(\Gamma)$ (see \cite[4.2]{bog}).
Given a vertex $x\in\Gamma$, its $\tau$-orbit
$x^{\tau}$ is the set $\{y\in\Gamma\ |\ y=\tau^lx,\ \text{for
  some $l\in\mathbb{Z}$}\}$. Also, we fix a polarisation
$\sigma$ in $\Gamma$. The periodic components of $\Gamma$ are
defined as follows. Consider the full translation subquiver of $\Gamma$
with vertices the periodic vertices in $\Gamma$. To this subquiver, add a
new arrow $x\to \tau x$ for every vertex $x$. A \emph{periodic
  component} of $\Gamma$ is a connected component of the resulting quiver.
Then:
\begin{enumerate}
\item The vertices of $\mathcal O(\Gamma)$ are the periodic components of
$\Gamma$ and the $\tau$-orbits of the non-periodic vertices.
\item  For each
periodic component, there is a loop attached to the associated vertex
in $\mathcal O(\Gamma)$.
\item  Let $u^{\sigma}$ be the $\sigma$-orbit of an
arrow $u\colon x\to y$. If both $x$ and $y$ are non-periodic, then there is
an edge
 between $x^{\tau}$ and $y^{\tau}$. If $x$ (or $y$) is
non-periodic and $y$ (or $x$) is periodic, then there is an edge 
between
$x^{\tau}$ (or $y^{\tau}$) and the vertex associated to the periodic component
containing $y$ (or $x$, respectively). Otherwise, no arrow
is associated to $u^{\sigma}$.
\end{enumerate}
By \cite[4.2]{bog}, the fundamental group of the orbit graph $\mathcal O(\Gamma)$ is
isomorphic to $\pi_1(\Gamma)$.

 Throughout this section, we assume
that $A$ is standard laura, having $\Gamma$ as a connecting
component. We use the following lemmata:
\begin{lem}
  \label{lem3.1}
If $\mathcal{O}(\Gamma)$ is a tree, then $A$ is weakly shod.
\end{lem}
\noindent{\textbf{Proof:}} If $\mathcal{O}(\Gamma)$ is a tree, then $\Gamma$
is simply connected (see \cite[4.1 and 4.2]{bog}). In particular, $\Gamma$ has no
oriented cycle. Hence, $A$ is laura and its non semiregular
component (there is at most one) has no oriented cycles. So $A$
is weakly shod (\cite[2.5]{cl}).\sq

\begin{lem}
  \label{lem3.2}
Let $A$ be a product of laura algebras with connecting components. If
the orbit graph of any connecting
component is a tree, then $A$ is a product of simply connected
algebras and $\ts{HH}^1(A)=0$.
\end{lem}
\noindent{\textbf{Proof:}} This follows from the preceding lemma and
from \cite[Cor. 2]{ws}.\sq

We now prove Theorem~\ref{thm1} whose statement we recall for convenience.
\setcounter{Thm}{1}
\begin{Thm}
\label{thm2.2}
Let $A$ be a standard laura algebra, and $\Gamma$ its connecting component(s). The following are
equivalent:
\begin{enumerate}
\item[(a)] $A$ has no proper Galois covering, that is, $A$ is simply connected.
\item[(b)] $\ts{HH}^1(A)=0$.
\item[(c)] $\Gamma$ is simply connected.
\item[(d)] The orbit graph $\mathcal O(\Gamma)$  is a tree.
\end{enumerate}
Moreover, if these conditions are verified, then $A$ is weakly shod.
\end{Thm}
\noindent{\textbf{Proof:}} By \cite[4.1, 4.2]{bog} and the
above lemma, (c) and (d) are equivalent and imply (a) and
(b). If $A$ is simply connected, then \ref{cor2.2}
implies $\pi_1(\Gamma)=1$. So (a) implies
(c). Finally, assume that
$\ts{HH}^1(A)=0$. By \ref{cor2.2}, the algebra $A$ admits a Galois
covering with group $\pi_1(\Gamma)$. This group is free because of
\cite[4.2]{bog}. On the other hand, the rank of $\pi_1(\Gamma)$ is less
than or equal to $\dim\ts{HH}^1(A)$ because of
\cite[Cor. 3]{PS01}. Therefore $\pi_1(\Gamma)=1$. So (b) implies (c). Thus the conditions
are equivalent, and imply that $A$ is weakly shod by \ref{lem3.1}.\sq

We illustrate Theorem~\ref{thm2.2} on the following examples. In
particular, note that this theorem does not necessarily hold true if
one drops  standardness.
\begin{ex}
  \begin{enumerate}
  \item[(a)] Let $A$ be as in \ref{ex1}, (a). Then $A$ clearly
    admits a Galois covering with group a free group of rank $3$ by a
    locally bounded $k$-category. It is given by the universal cover
    of the underlying graph of the ordinary quiver. So $A$ is not
    simply connected. The orbit graph $\mathcal O(\Gamma)$ of the connecting
    component $\Gamma$ is
    as follows:
    
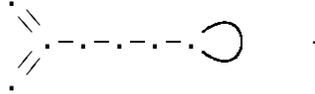
\begin{equation}
      \xymatrix@=5pt{
\centerdot \ar@<2pt>@{-}[rd] \ar@<-2pt>@{-}[rd]\\
& \centerdot \ar@{-}[r] & \centerdot \ar@{-}[r] & \centerdot
\ar@{-}[r] & \centerdot \ar@{-}[r]& \centerdot \ar@{-}@(ru,rd)&&&& .\\
\centerdot
\ar@<2pt>@{-}[ru] \ar@<-2pt>@{-}[ru]\\
}\notag
    \end{equation}
\vskip 12pt
Then $\pi_1(\Gamma)$  is free of rank $3$. A
straightforward computation gives $ \dim\ts{HH}^1(A)=7$ (see also
\cite[Thm. 1]{CS01}).
\item[(b)] Let $A$ be as in \ref{ex2}. As already noticed, $A$ is a
  simply connected representation-finite algebra. Also, it is not
  standard. The orbit graph of
  its Auslander-Reiten quiver is as follows:
  
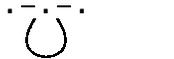
\begin{equation}
    \xymatrix@=5pt{
\centerdot \ar@{-}[r] &
\centerdot \ar@{-}@(ld,rd) \ar@{-}[r]&
\centerdot  &&&.\\
\null
}\notag
  \end{equation}
Finally, $A$ admits the following outer derivation, yielding a
non-zero element in $\ts{HH}^1(A)$ (see \cite[4.2]{bl}) 
\begin{equation}
  \begin{array}{crclc}
    d\colon & A &\to & A&\\
    & \sigma,\delta & \mapsto & 0&\\
    & \rho & \mapsto & \rho^3& .
  \end{array}\notag
\end{equation}
This example shows that Theorem~\ref{thm2.2} may fail if one drops
standardness. Note that the definition of simple connectedness
we use differs slightly from that used in\cite[4.3]{bl}: In \cite{bl},
as in \cite[\S 6]{bog}, a
representation-finite algebra is called simply connected if its
Auslander-Reiten quiver is simply connected.
\item[(c)] Let $A$ be given by the quiver:
  \begin{equation}
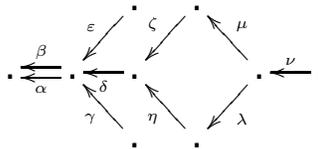
\begin{array}{cr}
      \xymatrix@=15pt{
& & \centerdot \ar[ld]_\varepsilon  & \centerdot\ar[ld]_\zeta \\
\centerdot & \centerdot \ar@<2pt>[l]^\alpha \ar@<-2pt>[l]_\beta  &\centerdot
\ar[l]^\delta & & \centerdot \ar[lu]_\mu \ar[ld]^\lambda & \centerdot \ar[l]_\nu\\
& & \centerdot \ar[lu]^\gamma  & \centerdot\ar[lu]^\eta 
}\notag\\
&
\end{array}
  \end{equation}
bound by $\beta\varepsilon=0,\,\alpha\gamma=0,\,\beta\delta=\alpha\delta,
\,\delta\zeta=0,\,\delta\eta=0,
\,\zeta\mu=\eta\lambda,\, 
\zeta\mu\nu=0$
Then $A$ is laura. Actually, it is right glued 
\cite[4.2]{ACLST05}. The orbit graph $\mathcal O(\Gamma)$ of its
connecting component $\Gamma$
is as follows:
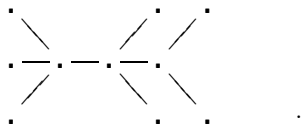
\begin{equation}
%       \xymatrix@=10pt{
% \centerdot \ar@{-}[rdd] & &\centerdot
% \\
% \centerdot \ar@{-}[rd] &&\centerdot
% \\ 
% & \centerdot  \ar@{-}[ruu]
% \ar@{-}[ru] \ar@{-}[rdd] \ar@{-}[rd] \ar@{-}[d]&
% \\
% \centerdot \ar@{-}[ru] & \centerdot& \centerdot
% \\
% \centerdot \ar@{-}[ruu] &&\centerdot& .
% }\notag
      \xymatrix@=10pt{
\centerdot \ar@{-}[rd]& & & \centerdot & \centerdot\\
\centerdot \ar@{-}[r]& \centerdot \ar@{-}[r]& \centerdot \ar@{-}[r]
\ar@{-}[ru] \ar@{-}[rd]& \centerdot \ar@{-}[ru] \ar@{-}[rd]\\
\centerdot \ar@{-}[ru]& & & \centerdot &  \centerdot& & .
}\notag
\end{equation}
It is a tree. Also, $A$ is simply connected, and it is not
hard to see that  $\ts{HH}^1(A)=0$ using, for instance, Happel's long exact sequence (see
\cite[5.3]{happel}).
  \end{enumerate}
\end{ex}

We end with the following problem.
\begin{pb}
  Let $A$ be a non-standard laura algebra. How can the vanishing of
  $\ts{HH}^1(A)$ be expressed in terms of topological properties of $A$?
\end{pb}

\section*{Acknowledgements}
$\ $

The first author gratefully acknowledges financial support from
the NSERC of Canada and the Universit\'e de Sherbrooke,  the second
 from Universidad San
Fransisco de Quito and NSERC of Canada  and  the third 
from the NSERC of Canada and
the \'Ecole normale sup\'erieure de Cachan. Part of the work presented
in this paper was done during visits of the second and third author at the
Universit\'e de Sherbrooke. They thank the
 algebra research group of Sherbrooke for their warm
hospitality.

\bibliographystyle{plain}
\bibliography{biblio}

\noindent Ibrahim Assem\\
\textit{e-mail:} ibrahim.assem@USherbrooke.ca\\
\textit{address:} D\'epartement de Math\'ematiques, Universit\'e
  de Sherbrooke, Sherbrooke, Qu\'ebec, Canada J1K 2R1\\

\noindent Juan Carlos Bustamante\\
\textit{e-mail:} jbustamante@ufsq.edu.ec\\
\textit{address:} Departamento de Matem\'aticas, Universidad
San Fransisco de Quito. Cumbay\'a - Quito,
  Ecuador\\

\noindent Patrick Le Meur\\
\textit{e-mail:} patrick.lemeur@cmla.ens-cachan.fr\\
\textit{address:} CMLA, ENS Cachan, CNRS, UniverSud,
61 Avenue du President Wilson, F-94230 Cachan, France
\end{document}